\documentclass[11pt]{article}
\usepackage{geometry,graphicx,latexsym,amssymb,verbatim}
\usepackage{multirow}
\usepackage{array}
\usepackage{multicol}
\usepackage{tikz}
\usetikzlibrary{calc}

\usepackage{color}

\def\vertex(#1){\put(#1){\circle*{2}}}
\def\vertexo(#1){\put(#1){\circle{2}}}
\def\vert(#1){\put(#1){\circle*{1.5}}}
\def\verto(#1){\put(#1){\circle{1.5}}}
\def\lab(#1)#2{\put(#1){\makebox(0,0)[c]{#2}}}
\newtheorem{thm}{Theorem}
\newtheorem{ob}[thm]{Observation}
\newtheorem{lem}[thm]{Lemma}
\newtheorem{prop}[thm]{Proposition}
\newtheorem{cor}[thm]{Corollary}

\newcommand{\dtd}{\gamma^d_t}
\newcommand{\ndtd}{\gamma^{d}_{nt}}

\newcommand{\gt}{\gamma_t}

\newcommand{\mod}{{\rm mod}}
\newcommand{\qed}{$\Box$}

\newcommand{\cB}{{\cal B}}

\newcommand{\cF}{{\cal F}}
\newcommand{\cG}{{\cal G}}
\newcommand{\cH}{{\cal H}}
\newcommand{\cL}{{\cal L}}

\newcommand{\cS}{{\cal S}}
\newcommand{\cC}{{\cal C}}
\newcommand{\cD}{{\cal D}}

\newcommand{\smallqed}{{\tiny ($\Box$)}}

\newenvironment{unnumbered}[1]{\trivlist
\item [\hskip \labelsep {\bf #1}]\ignorespaces\it}{\endtrivlist}

\newcommand{\1}{\vspace{0.1cm}}
\newcommand{\2}{\vspace{0.2cm}}

\begin{document}

\title{Graphs with Large Disjunctive Total Domination Number}
\author{$^1$Michael A. Henning\thanks{Research supported in part
by the South African National Research Foundation and the University of Johannesburg}  \, and $^{1,2}$Viroshan Naicker\thanks{Research supported in part by the University of Johannesburg, Rhodes University and the Claude Leon Foundation} \\
\\
$^1$Department of Pure and Applied Mathematics\\
University of Johannesburg \\
Auckland Park, 2006 South Africa \\
Email: mahenning@uj.ac.za\\
\\
$^2$Department of Mathematics \\
Rhodes University \\
Grahamstown, 6140 South Africa \\
Email: v.naicker@ru.ac.za
 }

\date{}
\maketitle

\begin{abstract}
Let $G$ be a graph with no isolated vertex. In this paper, we study a parameter that is a relaxation of arguably the most important domination parameter, namely the total domination number, $\gamma_t(G)$. A set $S$ of vertices in $G$ is a disjunctive total dominating set of $G$ if every vertex is adjacent to a vertex of $S$ or has at least two vertices in $S$ at distance~$2$ from it. The disjunctive total domination number, $\gamma^d_t(G)$, is the minimum cardinality of such a set. We observe that $\gamma^d_t(G) \le \gamma_t(G)$. Let $G$ be a connected graph on $n$ vertices with minimum degree $\delta$. It is known [J. Graph Theory 35 (2000), 21--45] that if $\delta \ge 2$ and $n \ge 11$, then $\gamma_t(G) \le 4n/7$. Further [J. Graph Theory 46 (2004), 207--210] if $\delta \ge 3$, then $\gamma_t(G) \le n/2$. We prove that if $\delta \ge 2$ and $n \ge 8$, then $\gamma^d_t(G) \le n/2$ and we characterize the extremal graphs.
\end{abstract}


\newpage
\section{Introduction}

In this paper we continue the study of disjunctive total domination in graphs introduced and studied by the authors in~\cite{HeN13}. As remarked in~\cite{HeN13}, a common issue in network design is to minimize the trade-off between resource allocation and redundancy. Key resources are usually expensive and cannot be allocated across an entire network, and, in addition, if there is a possibility of resource failure at a particular node, redundancy and backup requirements then become vital but require extra resources to be allocated. This problem has been addressed, in various guises, by using graphs as a model for the network and searching for vertex subsets which are `close' to the rest of the graph and satisfy pertinent redundancy criteria. A neural network that learns and attempts to optimally allocate resources as it grows has attracted considerable attention as is evidenced, for example, by the article of Platt~\cite{Pl91}.

Domination and, in particular, total domination are well studied topics in the graph theory literature which attempt a solution of this problem~(see, for example, \cite{hhs1,hhs2,He09,HeYe_book}.
Let $G$ be a graph that serves as a model of a network and let $G$ have vertex set $V$. On the one hand, for purposes of resource allocation, we select a set $D$ of vertices, called a \emph{dominating set}, of $G$ such that every vertex in $V \setminus D$ is adjacent to at least one vertex in $D$. On the other hand, for the purpose of extending the domination problem to include redundancy, we select a set $S$ of vertices, called a \emph{total dominating set}, abbreviated TD-set, of $G$ such that every vertex in $V$, including those in $S$, is adjacent to at least one vertex in $S$.
However, as remarked in~\cite{HeN13}, given the sheer scale of modern networks (see~\cite{Chung10}), many existing domination type structures are expensive to implement. Variations on the theme of dominating and total dominating sets studied to date tend to focus on adding restrictions which in turn raises their implementation costs.  As an alternative a relaxation of the domination number, called \emph{disjunctive domination}, was proposed and studied by Goddard et al.~\cite{GoHePi11}. This concept was recently extended in~\cite{HeN13} to a relaxation of total domination, called \emph{disjunctive total domination}. This new variant of total domination offers greater flexibility in the modelling of the network resource allocation problem while maintaining the redundancy and proximity features of the classical total domination parameter. In addition, as shown in \cite{HeN13} and by our main result below, there is a significant reduction in implementation cost over total domination in terms of the number of nodes of the network.

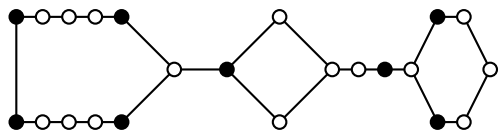
\begin{figure}[htb]
\tikzstyle{every node}=[circle, draw, fill=black!0, inner sep=0pt,minimum width=.18cm]
\begin{center}
\begin{tikzpicture}[thick,scale=.7]
  \draw(0,0) { 
    +(0.00,2.00) -- +(0.00,0.00)
    +(0.00,2.00) -- +(0.50,2.00)
    +(0.50,2.00) -- +(1.00,2.00)
    +(1.00,2.00) -- +(1.50,2.00)
    +(1.50,2.00) -- +(2.00,2.00)
    +(2.00,2.00) -- +(3.00,1.00)
    +(3.00,1.00) -- +(2.00,0.00)
    +(2.00,0.00) -- +(1.50,0.00)
    +(1.50,0.00) -- +(1.00,0.00)
    +(1.00,0.00) -- +(0.50,0.00)
    +(0.50,0.00) -- +(0.00,0.00)
    +(3.00,1.00) -- +(4.00,1.00)
    +(4.00,1.00) -- +(5.00,2.00)
    +(5.00,2.00) -- +(6.00,1.00)
    +(6.00,1.00) -- +(5.00,0.00)
    +(5.00,0.00) -- +(4.00,1.00)
   +(6.00,1.00) -- +(6.50,1.00)
    +(6.50,1.00) -- +(7.00,1.00)
    +(7.00,1.00) -- +(7.50,1.00)
    +(7.50,1.00) -- +(8.00,2.00)
    +(8.00,2.00) -- +(8.50,2.00)
    +(8.50,2.00) -- +(9.00,1.00)
    +(9.00,1.00) -- +(8.50,0.00)
    +(8.50,0.00) -- +(8.00,0.00)
    +(8.00,0.00) -- +(7.50,1.00)
    +(0.00,2.00) node[circle, draw,fill=black!100]{}
    +(0.50,2.00) node{}
    +(1.00,2.00) node{}
    +(1.50,2.00) node{}
    +(2.00,2.00) node[circle, draw,fill=black!100]{}
    +(2.00,0.00) node[circle, draw,fill=black!100]{}
    +(3.00,1.00) node{}
   +(4.00,1.00) node[circle, draw,fill=black!100]{}
    +(5.00,2.00) node{}
    +(6.00,1.00) node{}
    +(5.00,0.00) node{}
    +(6.50,1.00) node{}
    +(7.00,1.00) node[circle, draw,fill=black!100]{}
    +(7.50,1.00) node{}
    +(8.00,2.00) node[circle, draw,fill=black!100]{}
    +(8.50,2.00) node{}
    +(9.00,1.00) node{}
    +(8.00,0.00) node[circle, draw,fill=black!100]{}
    +(8.50,0.00) node{}
    +(1.50,0.00) node{}
    +(1.00,0.00) node{}
    +(0.50,0.00) node{}
    +(0.00,0.00) node[circle, draw,fill=black!100]{}
  };
\end{tikzpicture}
\end{center}
\vskip -0.6 cm \caption{A graph $G$ with $\dtd(G) = 8$ and $\gamma_t(G) = 11$.} \label{f:DTD5}
\end{figure}

\begin{figure}[htb]
\tikzstyle{every node}=[circle, draw, fill=black!0, inner sep=0pt,minimum width=.16cm]
\begin{center}
\begin{tikzpicture}[thick,scale=.5]
\draw(0,0){+(0,0)--+(0,1) +(0,1)--+(1,2) +(0,0)--+(1,-1)--+(2,0)--+(2,1)--+(1,2)
+(0,0) node[circle, draw, fill=black!100, inner sep=0pt,minimum width=.16cm]{} +(1,2) node[circle, draw, fill=black!100, inner sep=0pt,minimum width=.16cm]{} +(0,1) node{} +(2,1) node{} +(1,-1) node{} +(2,0) node[circle, draw, fill=black!100, inner sep=0pt,minimum width=.16cm]{}
};
\draw(5,0){+(0,0)--+(0,1) +(0,1)--+(1,2) +(0,0)--+(1,-1)--+(2,0)--+(2,1)--+(1,2)
+(0,0) node{} +(0,1) node[circle, draw, fill=black!100, inner sep=0pt,minimum width=.16cm]{} +(2,1) node[circle, draw, fill=black!100, inner sep=0pt,minimum width=.16cm]{} +(1,2) node[circle, draw, fill=black!100, inner sep=0pt,minimum width=.16cm]{} +(1,-1) node{} +(2,0) node{}
};
\draw(10,0){+(0,0)--+(0,1) +(0,1)--+(1,2) +(0,0)--+(1,-1)--+(2,0)--+(2,1)--+(1,2)
+(0,0) node[circle, draw, fill=black!100, inner sep=0pt,minimum width=.16cm]{} +(1,2) node[circle, draw, fill=black!100, inner sep=0pt,minimum width=.16cm]{} +(0,1) node{} +(2,1) node{} +(1,-1) node{} +(2,0) node[circle, draw, fill=black!100, inner sep=0pt,minimum width=.16cm]{} +(0,1)--+(-1,2) +(2,1)--+(3,2) +(3,2) node{} +(-1,2) node{}
};

\draw(15,0){+(0,0)--+(0,1) +(0,1)--+(1,2) +(0,0)--+(1,-1)--+(2,0)--+(2,1)--+(1,2)
+(0,0) node{} +(0,1) node[circle, draw, fill=black!100, inner sep=0pt,minimum width=.16cm]{} +(2,1) node[circle, draw, fill=black!100, inner sep=0pt,minimum width=.16cm]{} +(1,2) node[circle, draw, fill=black!100, inner sep=0pt,minimum width=.16cm]{} +(1,-1) node{} +(2,0) node{} +(1,1) node{}  +(1,0) node{} +(0,0)--+(1,1) +(1,0)--+(0,0) +(1,0)--+(2,0) +(1,1)--+(2,0)
};

\draw(1,-1)node[label=below:(a)]{};
\draw(6,-1)node[label=below:(b)]{};
\draw(11,-1)node[label=below:(c)]{};
\draw(16,-1)node[label=below:(d)]{};

\end{tikzpicture}
\end{center}
\vskip -0.6 cm \caption{Novel features of DTD-sets}\label{cyc:dtd}
\end{figure}
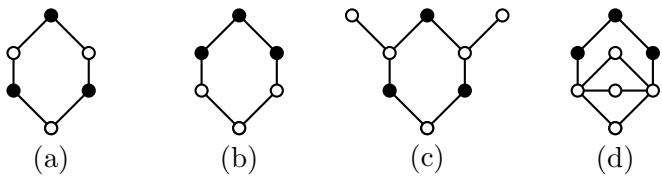

A set $S$ of vertices in $G$ is a \emph{disjunctive total dominating set}, abbreviated DTD-set, of $G$ if every vertex is adjacent to a vertex of $S$ or has at least two vertices in $S$ at distance~$2$ from it.
For example, the set of eight darkened vertices in the graph $G$ shown in Figure~\ref{f:DTD5} is a DTD-set of $G$.
We say that a vertex $v \in V$ is \emph{disjunctively totally dominated}, abbreviated DT-\emph{dominated}, by a set $S$, if $v$ has a neighbor in $S$ or if $v$ is at distance~$2$ from at least two vertices of $S$. Further if $v$ has a neighbor in $S$, we say $S$ \emph{totally dominates} the vertex $v$, while if $v$ is at distance~$2$ from at least two vertices of $S$, we say $S$ \emph{disjunctively dominates} the vertex $v$. The \emph{disjunctive total domination number}, $\dtd(G)$, is the minimum cardinality of a DTD-set in $G$. A DTD-set of cardinality $\dtd(G)$ is called a $\dtd(G)$-set.
Examples of two $\dtd(G)$-sets when $G$ is a cycle on six vertices are shown in Figure~\ref{cyc:dtd}(a) and Figure~\ref{cyc:dtd}(b). In each case the darkened vertices represent a $\dtd(G)$-set. We remark that a feature of the parameter is, depending on the graph, that a DTD-set may be chosen to be independent Figure~\ref{cyc:dtd}(a) or connected Figure~\ref{cyc:dtd}(b). The feature of choosing an independent DTD-set is essentially different from total domination in which independence is, by definition, lost. Further, if new vertices are added to the graph (nodes to a network) these may be added in a way that resources which are already allocated are, in a sense, `buffered' by their existing neighbors. For example, the graphs in Figure~\ref{cyc:dtd}(a) and Figure~\ref{cyc:dtd}(b) have been modified by adding two additional vertices to give Figure~\ref{cyc:dtd}(c) and Figure~\ref{cyc:dtd}(d), respectively, while the existing DTD-set is preserved. In each case the vertices need not be adjacent to a vertex in the DTD-set. This is a different feature, again, from total domination in which if a new vertex is joined to the graph by an edge either the existing TD-set must be enlarged and a new resource allocated to the network, or the new vertex has to be joined to a vertex in an existing TD-set which invites possible resource overcrowding by new users. 

\vspace{2mm}

The \emph{total domination number} of $G$, denoted by $\gt(G)$, is the minimum cardinality of a TD-set of $G$. Every TD-set is a DTD-set, implying the following observation.

\begin{ob}{\rm (\cite{HeN13})}
For every graph $G$ with no isolated vertex, $\dtd(G) \le \gt(G)$.
 \label{ob:Tdom}
\end{ob}

The known upper bounds on the total domination number of a graph $G$ in terms of its order~$n$ and small minimum degree~$\delta(G)$ are summarized in Table~1.

{\small
\[
\begin{array}{||ccrcllc||} \hline \hline %
& & & & & & \\
\delta(G) \ge 1 & \Rightarrow & \gt(G) & \le &
          \displaystyle{ \frac{2}{3} \, n } & \mbox{ if } n \ge 3 \mbox{ and $G$ is connected} & (\cite{CoDaHe80}) \\
& & & & & & \\
\delta(G) \ge 2 & \Rightarrow & \gt(G) & \le & \displaystyle{
        \frac{4}{7} \, n} & \mbox{ if } n \ge 11
        \mbox{ and $G$ is connected} & (\cite{He00}) \\
& & & & & &  \\
\delta(G) \ge 3 & \Rightarrow & \gt(G) &  \le & \displaystyle{ \frac{1}{2} \,  n } & & (\cite{Alfewy04,ChMc,Tu90}) \\
& & & & & & \\
\hline \hline
\end{array}
\]
}
\begin{center}
{Table~1:} Upper bounds on the total domination number of a graph $G$.
\end{center}

By Observation~\ref{ob:Tdom} and the result of \cite{CoDaHe80} shown in Table~1, if $G$ is a connected graph of order~$n \ge 3$, then $\dtd(G) \le 2n/3$.  The authors showed in~\cite{HeN13} that this upper bound on the disjunctive total domination number can be improved ever-so-slightly. Further, they characterized the concomitant extremal graphs.
\begin{thm}
\label{1:thm} {\rm (\cite{HeN13})}
If $G$ is a connected graph of order $n \ge 8$, then $\dtd(G) \le 2(n-1)/3$, and this bound is tight.
\end{thm}

In addition, the authors showed in \cite{HeN13} that the upper bound in Table~1 of $\gt(G)\le 4n/7$ when $G$ has order $n\ge 11$ and $\delta(G)\ge 2$ holds for the disjunctive total domination number when the minimum degree is relaxed from 2 to 1 and $G$ is restricted to the class of connected claw-free graphs of order $n>14$. Further, they characterized the concomitant extremal graphs.

\begin{thm}
\label{1cf:thm} {\rm (\cite{HeN13})}
If $G$ is a connected, claw-free, graph of order $n > 14$, then $\dtd(G) \le 4n/7$, and this bound is tight.
\end{thm}

In this paper, we show that if we restrict the minimum degree to at least~$2$, then the result of \cite{He00} shown in Table~1 on the total domination number can be improved significantly for the disjunctive total domination number. Perhaps surprisingly we show that the upper bound of one-half the order of a graph of \cite{Alfewy04,ChMc,Tu90} shown in Table~1 on the total domination number is an upper bound on the disjunctive  total domination number even if we relax the minimum degree from~$3$ to~$2$. More precisely, we shall prove the following result.

\begin{thm}
\label{t:main1}
If $G$ is a connected graph of order $n \ge 13$, then $\dtd(G) \le (n-1)/2$, and this bound is tight.
\end{thm}

Table~2 summarizes the upper bounds on the disjunctive total domination number of a graph $G$ in terms of its order~$n$ and minimum degree~$\delta(G)$.

{\small
\[
\begin{array}{||ccrcllc||} \hline \hline %
& & & & & & \\
\delta(G) \ge 1 & \Rightarrow & \dtd(G) & \le &
          \displaystyle{ \frac{2}{3} \, (n-1) } & \mbox{ if } n \ge 8 \mbox{ and $G$ is connected} & (\cite{HeN13}) \\
& & & & & & \\
\delta(G) \ge 1 & \Rightarrow & \dtd(G) & \le &
          \displaystyle{ \frac{4}{7} \, n } & \mbox{ if } n \ge 14 \mbox{ and $G$ is connected, claw-free} & (\cite{HeN13}) \\
& & & & & & \\
\delta(G) \ge 2 & \Rightarrow & \dtd(G) &  \le & \displaystyle{ \frac{1}{2} \,  (n-1) } & \mbox{ if } n \ge 13 \mbox{ and $G$ is connected}  & (Theorem~\ref{t:main1}) \\
& & & & & & \\
\hline \hline
\end{array}
\]
}
\begin{center}
{Table~2:} Upper bounds on the disjunctive total domination number of a graph $G$.
\end{center}

\subsection{Notation}

For notation and graph theory terminology, we in general follow \cite{hhs1}. Specifically, let $G=(V,E)$ be a graph with vertex set $V$ and edge set $E$. We denote the \emph{degree} of $v$ in $G$ by $d_G(v)$. A vertex of degree~$k$ is called a \emph{degree}-$k$ vertex. The maximum (minimum) degree among the vertices of $G$ is denoted by $\Delta(G)$ ($\delta(G)$, respectively).
For a set $S \subseteq V$, the subgraph induced by $S$ is denoted by $G[S]$, while the graph obtained from $G$ be removing all vertices in $S$ and their incident edges is denoted by $G - S$. For two vertices $u$ and $v$ in a connected graph $G$, the \emph{distance} $d_G(u,v)$ between $u$ and $v$ is the length of a shortest $u$--$v$ path in $G$.
%
The \emph{open neighborhood} of a vertex $v$ is the set $N_G(v) = \{u \in V \, | \, uv \in E\}$ and the \emph{closed neighborhood of $v$} is $N_G[v] = \{v\} \cup N_G(v)$. For a set $S\subseteq V$, its \emph{open neighborhood} is the set
$N_G(S) = \bigcup_{v \in S} N_G(v)$,
and its \emph{closed neighborhood} is the set $N_G[S] = N_G(S) \cup S$. If the graph $G$ is clear from the context, we simply write $N(v)$, $N[v]$, $N(S)$, $N[S]$, $d(v)$ and $d(u,v)$ rather than $N_{G}(v)$, $N_{G}[v]$, $N_G(S)$, $N_G[S]$, $d_{G}(v)$, and $d_{G}(u,v)$, respectively.

A \emph{cycle} and \emph{path} on $n$ vertices are denoted by $C_n$ and $P_n$, respectively. For $m \ge 3$ and $n \ge 1$, we denote by $L_{m,n}$ the graph obtained by joining with an edge a vertex in $C_m$ to an end-vertex of $P_n$. The graph $L_{m,n}$ is called a \emph{key}.


A $4$-\emph{subdivision} of a graph $G$ is the graph obtained from $G$ by subdividing an edge of $G$ four times. If $e$ is an edge of a graph $G$, we let $G_e$ denote the graph obtained from $G$ by subdividing the edge $e$ four times. If two graphs $G$ and $H$ are isomorphic, we write $G \cong H$. %
Let $G$ be a graph with $\delta(G) \ge 2$. We define a vertex $v$ of $G$ to be \emph{large} if $d_G(v) \ge 3$ and \emph{small} if $d_G(v) = 2$. A \emph{cycle edge} of $G$ is an edge that belongs to a cycle in $G$.
We call a vertex, $v$, in a graph $G$ a \emph{good-vertex} of $G$ if it belongs to some $\dtd(G)$-set; otherwise, we call $v$ a \emph{bad-vertex} of $G$. We call an edge, $e$, in $G$ a \emph{good-edge} of $G$ if there is a $\dtd(G)$-set which contains both vertices incident with $e$; otherwise, we call $e$ a \emph{bad-edge}.

\subsection{Special Families}
\label{s3:sec}

In this section we define several special families of graphs. Let
\[
\cC = \{C_{3}, C_{4}, C_{5}, C_{6}, C_{7}, C_{8}, C_{9}, C_{11}, C_{12}, C_{13}, C_{17}\}
\]
be a family of cycles.
A \emph{daisy} with $k\ge 2$ \emph{petals} is a connected graph that can be constructed from $k\ge 2$ disjoint cycles by identifying a set of $k$ vertices, one from each cycle into one vertex. If the cycles have lengths $n_{1}$, $n_{2}$,\ldots, $n_{k}$, we denote the daisy by $D(n_{1}, n_{2}, \ldots, n_{k})$. Let
$
\cD = \{D(3,3), D(4,4), D(3,7)\}
$
be the family of three daisies shown in Figure~\ref{f:cD}.

\begin{figure}[htb]
\tikzstyle{every node}=[circle, draw, fill=black!0, inner sep=0pt,minimum width=.16cm]
\begin{center}
\begin{tikzpicture}[thick,scale=.7]

\draw(-0.5,5) {+(0.50,0.50)--+(0.00,1.00)--+(0.00,0.00)--+(0.50,0.50) +(0.50,0.50)--+(1.00,1.00)--+(1.00,0.00)--+(0.50,0.50)
+(0.50,0.50) node{} +(0.00,1.00) node{} +(0.00,0.00) node{} +(1.00,0.00) node{} +(1.00,1.00) node{}
+(0.5,-0.75) node[rectangle, draw=white!0, fill=white!100]{$D(3,3)$}
};

\draw(8.5,5){+(0.50,0.50)--+(0.00,1.00)--+(0.00,0.00)--+(0.50,0.50) +(0.50,0.50)--+(1.00,1.00)--+(2.00,1.00)--+(2.00,0.00)--+(1.00,0.00)--+(0.50,0.50)
 +(0.50,0.50) node{} +(0.00,1.00) node{} +(0.00,0.00) node{} +(1.00,1.00) node{}  +(1.50,1.00) node{}  +(2.00,1.00) node{}  +(1.50,0.00) node{}  +(2.00,0.00) node{}  +(1.00,0.00) node{}
+(1,-0.75) node[rectangle, draw=white!0, fill=white!100]{$D(3,7)$}
};

\draw(3.5,5){+(1.00,0.50)--+(0.50,1.00)--+(0.00,0.50)--+(0.50,0.00)--+(1.00,0.50) +(1.00,0.50)--+(1.50,1.00)--+(2.00,0.50)--+(1.50,0.00)--+(1.00,0.50)
+(1.00,0.50) node{} +(0.50,1.00) node{} +(0.00,0.50) node{} +(0.50,0.00) node{}
+(1.50,1.00) node{} +(2.00,0.50) node{} +(1.50,0.00) node{}
+(1,-0.75) node[rectangle, draw=white!0, fill=white!100]{$D(4,4)$}
};


\end{tikzpicture}
\end{center}
\vskip -0.5 cm \caption{The family, $\cD$, of daisies.} \label{f:cD}
\end{figure}

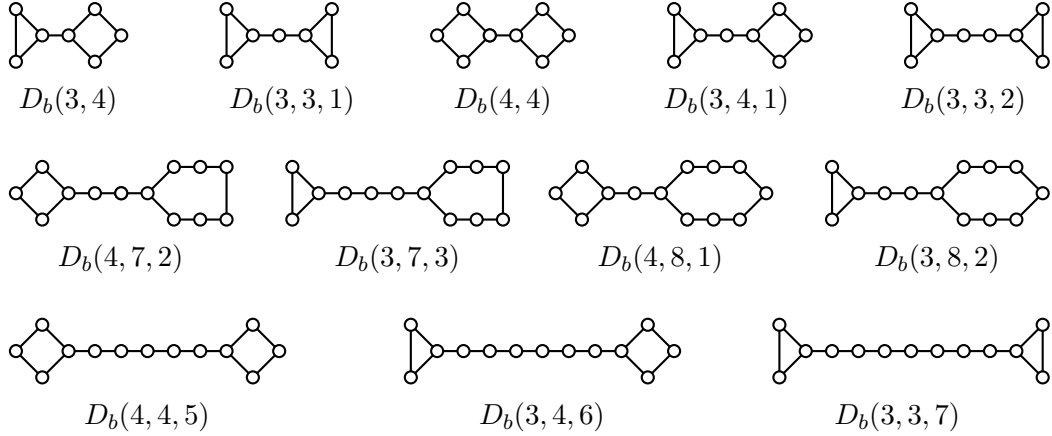
\begin{figure}[htb]
\tikzstyle{every node}=[circle, draw, fill=black!0, inner sep=0pt,minimum width=.16cm]
\begin{center}
\begin{tikzpicture}[thick,scale=.7]

\draw(0,15){+(0.50,0.50)--+(0.00,1.00)--+(0.00,0.00)--+(0.50,0.50) +(0.50,0.50)--+(1.00,0.50)--+(1.50,1.00)--+(2.00,0.50)--+(1.50,0.00)--+(1.00,0.50)
+(0.00,0.00) node{} +(0.00,1.00) node{} +(0.50,0.50) node{} +(1.00,0.50) node{} +(1.50,1.00) node{} +(2.00,0.50) node{} +(1.50,0.00) node{}
+(1.00,-0.75) node[rectangle, draw=white!0, fill=white!100]{$D_{b}(3,4)$}
};

\draw(8,15){+(0.00,0.50)--+(0.50,1.00)--+(1.00,0.50)--+(0.50,0.00)--+(0.00,0.50)
+(1.00,0.50)--+(1.50,0.50)--+(2.00,1.00)--+(2.50,0.50)--+(2.00,0.00)--+(1.50,0.50)
+(0.00,0.50) node{} +(0.50,1.00) node{} +(1.00,0.50) node{} +(0.50,0.00) node{} +(1.00,0.50) node{} +(2.00,1.00) node{} +(2.00,0.00) node{} +(2.50,0.50) node{} +(1.50,0.50) node {}
+(1.25,-0.75) node[rectangle, draw=white!0, fill=white!100]{$D_{b}(4,4)$}
};

\draw(4,15){+(0.50,0.50)--+(0.00,1.00)--+(0.00,0.00)--+(0.50,0.50) +(0.50,0.50)--+(1.00,0.50) +(1.00,0.50)--+(1.50,0.50)--+(2.00,1.00)--+(2.00,0.00)--+(1.50,0.50)
+(0.00,0.00) node{} +(0.00,1.00) node{} +(0.50,0.50) node{}
+(1.50,0.50) node{} +(2.00,1.00) node{} +(2.00,0.00) node{}
+(1.00,0.50) node{}
+(1.25,-0.75) node[rectangle, draw=white!0, fill=white!100]{$D_{b}(3,3,1)$}
};

\draw(17,15){+(0.50,0.50)--+(0.00,1.00)--+(0.00,0.00)--+(0.50,0.50) +(0.50,0.50)--+(2.00,0.50)--+(2.50,1.00)--+(2.50,0.00)--+(2.00,0.50)
+(0.00,0.00) node{} +(0.00,1.00) node{} +(0.50,0.50) node{}
+(2.00,0.50) node{} +(2.50,1.00) node{} +(2.50,0.00) node{}
+(1.00,0.50) node{} +(1.50,0.50) node{}
+(1.00,-0.75) node[rectangle, draw=white!0, fill=white!100]{$D_{b}(3,3,2)$}
};

\draw(12.5,15){+(0.50,0.50)--+(0.00,1.00)--+(0.00,0.00)--+(0.50,0.50) +(0.50,0.50)--+(1.50,0.50) +(1.50,0.50)--+(2.00,1.00)--+(2.50,0.50)--+(2.00,0.00)--+(1.50,0.50)
+(0.00,0.00) node{} +(0.00,1.00) node{} +(0.50,0.50) node{}
+(1.00,0.50) node{}
+(1.50,0.50) node{} +(2.00,1.00) node{} +(2.50,0.50) node{} +(2.00,0.00) node{}
+(1.00,-0.75) node[rectangle, draw=white!0, fill=white!100]{$D_{b}(3,4,1)$}
};

\draw(5.25,12){+(0.50,0.50)--+(0.00,1.00)--+(0.00,0.00)--+(0.50,0.50) +(0.50,0.50)--+(2.50,0.50)--+(3.00,1.00)--+(4.00,1.00)--+(4.00,0.00)--+(3.00,0.00)--+(2.50,0.50)
+(0.00,0.00) node{} +(0.00,1.00) node{} +(0.50,0.50) node{}
+(1.00,0.50) node{} +(1.50,0.50) node{} +(2.00,0.50) node{}
+(2.50,0.50) node{} +(3.00,1.00) node{} +(3.50,1.00) node{} +(4.00,1.00) node{} +(4.00,0.00) node{} +(3.50,0.00) node{} +(3.00,0.00) node{}
+(2.00,-0.75) node[rectangle, draw=white!0, fill=white!100]{$D_{b}(3,7,3)$}
};

\draw(0.0,12){+(1.00,0.50)--+(0.50,1.00)--+(0.00,0.50)--+(0.50,0.00)--+(1.00,0.50) +(1.00,0.50)--+(2.50,0.50) +(2.50,0.50)--+(3.00,1.00)--+(4.00,1.00)--+(4.00,0.00)--+(3.00,0.00)--+(2.50,0.50) 
+(1.00,0.50) node{} +(0.50,1.00) node{} +(0.00,0.50) node{} +(0.50,0.00) node{}
+(2.00,0.50) node{} +(4.00,1.00) node{} +(3.00,1.00) node{} +(3.50,1.00) node{} +(3.50,0.00) node{} +(3.00,0.00) node{} +(4.00,0.00) node{}
+(1.50,0.50) node{} +(2.00,0.50) node{} +(2.50,0.50) node{}
+(2.00,-0.75) node[rectangle, draw=white!0, fill=white!100]{$D_{b}(4,7,2)$}
};

\draw(0,9){+(1.00,0.50)--+(0.50,1.00)--+(0.00,0.50)--+(0.50,0.00)--+(1.00,0.50) +(1.00,0.50)--+(4.00,0.50) +(4.00,0.50)--+(4.50,1.00)--+(5.00,0.50)--+(4.50,0.00)--+(4.00,0.50)
+(1.00,0.50) node{} +(0.50,1.00) node{} +(0.00,0.50) node{} +(0.50,0.00) node{}
+(4.00,0.50) node{} +(4.50,1.00) node{} +(5.00,0.50) node{} +(4.50,0.00) node{}
+(1.50,0.50) node{} +(2.00,0.50) node{} +(2.50,0.50) node{} +(3.00,0.50) node{} +(3.50,0.50) node{}
+(2.50,-0.75) node[rectangle, draw=white!0, fill=white!100]{$D_{b}(4,4,5)$}
};

\draw(7.5,9){+(0.50,0.50)--+(0.00,1.00)--+(0.00,0.00)--+(0.50,0.50) +(0.50,0.50)--+(4.00,0.50) +(4.00,0.50)--+(4.50,1.00)--+(5.00,0.50)--+(4.50,0.00)--+(4.00,0.50)
+(0.00,0.00) node{} +(0.00,1.00) node{} +(0.50,0.50) node{} +(3.50,0.50) node{}
+(1.00,0.50) node{} +(1.50,0.50) node{} +(2.00,0.50) node{} +(2.50,0.50) node{} +(3.00,0.50) node{}
+(4.00,0.50) node{} +(4.50,1.00) node{} +(5.00,0.50) node{} +(4.50,0.00) node{}

+(2.50,-0.75) node[rectangle, draw=white!0, fill=white!100]{$D_{b}(3,4,6)$}
};

\draw(14.5,9){+(0.50,0.50)--+(0.00,1.00)--+(0.00,0.00)--+(0.50,0.50) +(0.50,0.50)--+(4.50,0.50) +(4.50,0.50)--+(5.00,1.00)--+(5.00,0.00)--+(4.50,0.50)
+(0.00,0.00) node{} +(0.00,1.00) node{} +(0.50,0.50) node{}
+(1.00,0.50) node{} +(1.50,0.50) node{} +(2.00,0.50) node{} +(2.50,0.50) node{} +(3.00,0.50) node{} +(3.50,0.50) node{} +(4.00,0.50) node{}
+(4.50,0.50) node{} +(5.00,0.00) node{} +(5.00,1.00) node{}
+(2.25,-0.75) node[rectangle, draw=white!0, fill=white!100]{$D_{b}(3,3,7)$}
};

\draw(10.25,12){+(1.00,0.50)--+(0.50,1.00)--+(0.00,0.50)--+(0.50,0.00)--+(1.00,0.50) +(1.00,0.50)--+(2.00,0.50) +(2.00,0.50)--+(2.50,1.00)--+(3.50,1.00)--+(4.00,0.50)--+(3.50,0.00)--+(2.50,0.00)--+(2.00,0.50)
+(1.00,0.50) node{} +(0.50,1.00) node{} +(0.00,0.50) node{} +(0.50,0.00) node{}
+(1.50,0.50) node{}
+(2.00,0.50) node{} +(2.50,1.00) node{} +(3.00,1.00) node{} +(3.50,1.00) node{} +(3.50,0.00) node{} +(4.00,0.50) node{} +(3.00,0.00) node{} +(2.50,0.00) node{}
+(2.00,-0.75) node[rectangle, draw=white!0, fill=white!100]{$D_{b}(4,8,1)$}
};

\draw(15.5,12){+(0.50,0.50)--+(0.00,1.00)--+(0.00,0.00)--+(0.50,0.50) +(0.50,0.50)--+(2.00,0.50) +(2.00,0.50)--+(2.50,1.00)--+(3.50,1.00)--+(4.00,0.50)--+(3.50,0.00)--+(2.50,0.00)--+(2.00,0.50)

+(0.50,0.50) node{} +(0.00,0.00) node{} +(0.00,1.00) node{} +(1.00,0.50) node{}
+(1.50,0.50) node{}
+(2.00,0.50) node{} +(2.50,1.00) node{} +(3.00,1.00) node{} +(3.50,1.00) node{} +(3.50,0.00) node{} +(4.00,0.50) node{} +(3.00,0.00) node{} +(2.50,0.00) node{}
+(2.00,-0.75) node[rectangle, draw=white!0, fill=white!100]{$D_{b}(3,8,2)$}
};

\end{tikzpicture}
\end{center}
\vskip -0.5 cm \caption{The family, $\cD_{b}$, of dumb-bell graphs.} \label{f:cDb}
\end{figure}

A \emph{dumb-bell} is a connected graph on $n=n_{1}+n_{2}+\ell$ vertices that can be constructed by joining a vertex of a cycle $C_{n_{1}}$ to a vertex of a cycle $C_{n_{2}}$ by an edge and subdividing this edge $\ell$ times. If $\ell = 0$ we denote the dumb-bell by $D_{b}(n_{1}, n_{2})$, and if $\ell \ge 1$ we denote the dumb-bell by $D_{b}(n_{1}, n_{2}, \ell)$. Let
\[
\begin{array}{lcl}
\cD_b & = & \{ D_b(3,4), D_b(3,3,1), D_b(4,4), D_b(3,4,1), D_b(3,3,2), D_{b}(4,7,2), \\
& & D_{b}(3,7,3), D_b(4,8,1), D_{b}(3,8,2), D_{b}(4,4,5), D_{b}(3,4,6), D_{b}(3,4,7) \}
\end{array}
\]
be the family of twelve dumb-bells shown in Figure \ref{f:cDb}.

\begin{figure}[htb]
\tikzstyle{every node}=[circle, draw, fill=black!0, inner sep=0pt,minimum width=.16cm]
\begin{center}
\begin{tikzpicture}[thick,scale=.7]

\draw(6.5,4) {+(0.00,0.00)--+(1.00,0.00)--+(0.50,0.50)--+(0.00,0.00) +(0.50,0.50)--+(0.50,1.50) +(0.00,0.00) node{} +(1.00,0.00) node{} +(0.50,0.50) node{} +(0.50,1.00) node{}
    +(0.50,1.50) node[circle, draw, fill=black!100, inner sep=0pt,minimum width=.16cm]{}
    +(0.5,-0.75) node[rectangle, draw=white!0, fill=white!100]{$U_{1}$}
    +(0.5,1.9) node[rectangle, draw=white!0, fill=white!100]{$v$}
};

\draw(9.5,4) {+(0.00,0.50)--+(0.50,0.00)--+(1.00,0.50)--+(0.50,1.00)--+(0.00,0.50) +(0.50,1.00)--+(0.50,1.50) +(0.00,0.50) node{} +(0.50,0.00) node{} +(1.00,0.50) node{} +(0.50,1.00) node{} +(0.50,1.50) node[circle, draw, fill=black!100, inner sep=0pt,minimum width=.16cm]{}
+(0.5,-0.75) node[rectangle, draw=white!0, fill=white!100]{$U_{2}$}
+(0.5,1.9) node[rectangle, draw=white!0, fill=white!100]{$v$}
};

\draw(12.5,4) {+(0.00,0.50)--+(0.50,0.00)--+(1.00,0.50)--+(0.50,1.00)--+(0.00,0.50) +(0.50,1.00)--+(0.50,1.50) +(0.00,0.50)--+(1.00,0.50) +(0.00,0.50) node{} +(0.50,0.00) node{} +(1.00,0.50) node{} +(0.50,1.00)  node{} +(0.50,1.50) node[circle, draw, fill=black!100, inner sep=0pt,minimum width=.16cm]{}
+(0.5,-0.75) node[rectangle, draw=white!0, fill=white!100]{$U_{3}$}
+(0.5,1.9) node[rectangle, draw=white!0, fill=white!100]{$v$}
};

\draw(1.5,0) {+(0.00,0.00)--+(1.00,0.00)--+(1.00,0.50)--+(0.50,1.00)--+(0.00,0.50)--+(0.00,0.00) +(0.00,0.00) node{} +(0.00,0.50) node{} +(1.00,0.50) node{} +(1.00,0.00) node{} +(0.50,1.00) node[circle, draw, fill=black!100, inner sep=0pt,minimum width=.16cm]{}
+(0.5,-0.75) node[rectangle, draw=white!0, fill=white!100]{$X_{1}$}
+(0.5,1.35) node[rectangle, draw=white!0, fill=white!100]{$v$}
};

\draw(4.5,0) {+(0.00,0.50)--+(0.50,0.00)--+(1.00,0.50)--+(1.00,1.00)--+(0.50,1.50)--+(0.00,1.00)--+(0.00,0.50) +(0.50,1.50)--+(0.50,2.00) +(0.00,1.00) node{} +(0.50,0.00) node{} +(0.00,0.50) node{} +(1.00,0.50) node{} +(1.00,1.00) node{} +(0.50,1.50) node{} +(0.50,2.00) node[circle, draw, fill=black!100, inner sep=0pt,minimum width=.16cm]{}
+(0.5,-0.75) node[rectangle, draw=white!0, fill=white!100]{$X_{2}$}
+(0.5,2.4) node[rectangle, draw=white!0, fill=white!100]{$v$}
};

\draw(7.5,0) {+(0.00,0.00)--+(1.00,0.00)--+(0.50,0.50)--+(0.00,0.00) +(0.50,0.50)--+(0.50,2.00) +(0.00,0.00) node{} +(1.00,0.00) node{} +(0.50,0.50) node{} +(0.50,1.50) node{} +(0.50,1.00) node{} +(0.50,1.50)--+(1.00,1.50) +(1.00,1.50)--+(1.50,2.00)--+(1.50,1.00)--+(1.00,1.50) +(1.00,1.50) node{} +(1.50,2.00) node{} +(1.50,1.00) node{} +(0.50,2.00) node[circle, draw, fill=black!100, inner sep=0pt,minimum width=.16cm]{}
+(0.5,-0.75) node[rectangle, draw=white!0, fill=white!100]{$X_{3}$}
+(0.5,2.4) node[rectangle, draw=white!0, fill=white!100]{$v$}
};

\draw(11,0) {+(0.00,0.00)--+(1.00,0.00)--+(0.50,0.50)--+(0.00,0.00) +(0.50,0.50)--+(0.50,2.00) +(0.00,0.00) node{} +(1.00,0.00) node{} +(0.50,0.50) node{} +(0.50,1.50) node{} +(0.50,1.00) node{} +(0.50,1.50)--+(1.00,2.00)--+(1.50,1.50)--+(1.00,1.00)--+(0.50,1.50) +(1.50,1.50) node{} +(1.00,2.00) node{} +(1.00,1.00) node{} +(0.50,2.00) node[circle, draw, fill=black!100, inner sep=0pt,minimum width=.16cm]{}
+(0.5,-0.75) node[rectangle, draw=white!0, fill=white!100]{$X_{4}$}
+(0.5,2.4) node[rectangle, draw=white!0, fill=white!100]{$v$}
};

\draw(14.5,0) {+(0.00,0.50)--+(0.50,0.00)--+(1.00,0.50)--+(0.50,1.00)--+(0.00,0.50) +(0.50,1.00)--+(0.50,2.00) +(0.00,0.50) node{} +(0.50,0.00) node{} +(1.00,0.50) node{} +(0.50,1.00) node{} +(0.50,1.50) node{} +(0.50,1.50)--+(1.00,1.50) +(1.00,1.50)--+(1.50,2.00)--+(1.50,1.00)--+(1.00,1.50) +(1.00,1.50) node{} +(1.50,2.00) node{} +(1.50,1.00) node{} +(0.50,2.00) node[circle, draw, fill=black!100, inner sep=0pt,minimum width=.16cm]{}
+(0.5,-0.75) node[rectangle, draw=white!0, fill=white!100]{$X_{5}$}
+(0.5,2.4) node[rectangle, draw=white!0, fill=white!100]{$v$}
};

\draw(17.5,0) {+(0.00,0.50)--+(0.50,0.00)--+(1.00,0.50)--+(0.50,1.00)--+(0.00,0.50) +(0.50,1.00)--+(0.50,2.00) +(0.00,0.50) node{} +(0.50,0.00) node{} +(1.00,0.50) node{} +(0.50,1.00) node{} +(0.50,1.50) node{} +(1.00,1.00) node{} +(1.50,1.50) node{} +(1.00,2.00) node{} +(0.50,2.00) node[circle, draw, fill=black!100, inner sep=0pt,minimum width=.16cm]{} +(0.50,1.50)--+(1.00,2.00)--+(1.50,1.50)--+(1.00,1.00)--+(0.50,1.50) +(1.50,1.50)
+(0.5,-0.75) node[rectangle, draw=white!0, fill=white!100]{$X_{6}$}
+(0.5,2.4) node[rectangle, draw=white!0, fill=white!100]{$v$}
};

\draw(4.5,-4) {+(0.00,0.00)--+(1.00,0.00)--+(1.00,0.50)--+(0.50,1.00)--+(0.00,0.50)--+(0.00,0.00) +(0.00,0.00) node{} +(0.00,0.50) node{} +(1.00,0.50) node{} +(1.00,0.00) node{} +(0.50,1.00) node[circle, draw, fill=black!100, inner sep=0pt,minimum width=.16cm]{} +(1.00,0.50)--+(0.00,0.50) +(0.5,-0.75) node[rectangle, draw=white!0, fill=white!100]{$X_{7}$}
+(0.5,1.4) node[rectangle, draw=white!0, fill=white!100]{$v$}
};

\draw(7.5,-4) {+(0.00,0.00)--+(1.00,0.00)--+(0.50,0.50)--+(0.00,0.00) +(0.50,0.50)--+(0.50,2.00) +(0.00,0.00) node{} +(1.00,0.00) node{} +(0.50,0.50) node{} +(0.50,1.50) node{} +(0.50,1.00) node{} +(0.50,1.50)--+(1.00,2.00)--+(1.50,1.50)--+(1.00,1.00)--+(0.50,1.50) +(1.50,1.50) node{} +(1.00,2.00) node{} +(1.00,1.00) node{} +(0.50,2.00) node[circle, draw, fill=black!100, inner sep=0pt,minimum width=.16cm]{} +(1.00,1.00)--+(1.00,2.00)
+(0.5,-0.75) node[rectangle, draw=white!0, fill=white!100]{$X_{8}$}
+(0.5,2.4) node[rectangle, draw=white!0, fill=white!100]{$v$}
};

\draw(11,-4) {+(0.00,0.50)--+(0.50,0.00)--+(1.00,0.50)--+(0.50,1.00)--+(0.00,0.50) +(0.50,1.00)--+(0.50,2.00) +(0.00,0.50) node{} +(0.50,0.00) node{} +(1.00,0.50) node{} +(0.50,1.00) node{} +(0.50,1.50) node{} +(0.50,1.50)--+(1.00,1.50) +(1.00,1.50)--+(1.50,2.00)--+(1.50,1.00)--+(1.00,1.50) +(1.00,1.50) node{} +(1.50,2.00) node{} +(1.50,1.00) node{} +(0.50,2.00) node[circle, draw, fill=black!100, inner sep=0pt,minimum width=.16cm]{} +(0.00,0.50)--+(1.00,0.50)
+(0.5,-0.75) node[rectangle, draw=white!0, fill=white!100]{$X_{9}$}
+(0.5,2.4) node[rectangle, draw=white!0, fill=white!100]{$v$}
};

\draw(14.5,-4) {+(0.00,0.50)--+(0.50,0.00)--+(1.00,0.50)--+(0.50,1.00)--+(0.00,0.50) +(0.50,1.00)--+(0.50,2.00) +(0.00,0.50) node{} +(0.50,0.00) node{} +(1.00,0.50) node{} +(0.50,1.00) node{} +(0.50,1.50) node{} +(1.00,1.00) node{} +(1.50,1.50) node{} +(1.00,2.00) node{} +(0.50,2.00) node[circle, draw, fill=black!100, inner sep=0pt,minimum width=.16cm]{} +(0.50,1.50)--+(1.00,2.00)--+(1.50,1.50)--+(1.00,1.00)--+(0.50,1.50) +(1.50,1.50) +(1.00,1.00)--+(1.00,2.00) +(0.00,0.50)--+(1.00,0.50)
+(0.5,-0.75) node[rectangle, draw=white!0, fill=white!100]{$X_{10}$}
+(0.5,2.4) node[rectangle, draw=white!0, fill=white!100]{$v$}
};

\end{tikzpicture}
\end{center}
\vskip -0.7 cm \caption{The graphs $U_1,U_2,U_3$ and $X_1,X_2,\ldots,X_{10}$.} \label{f:UV}
\end{figure}
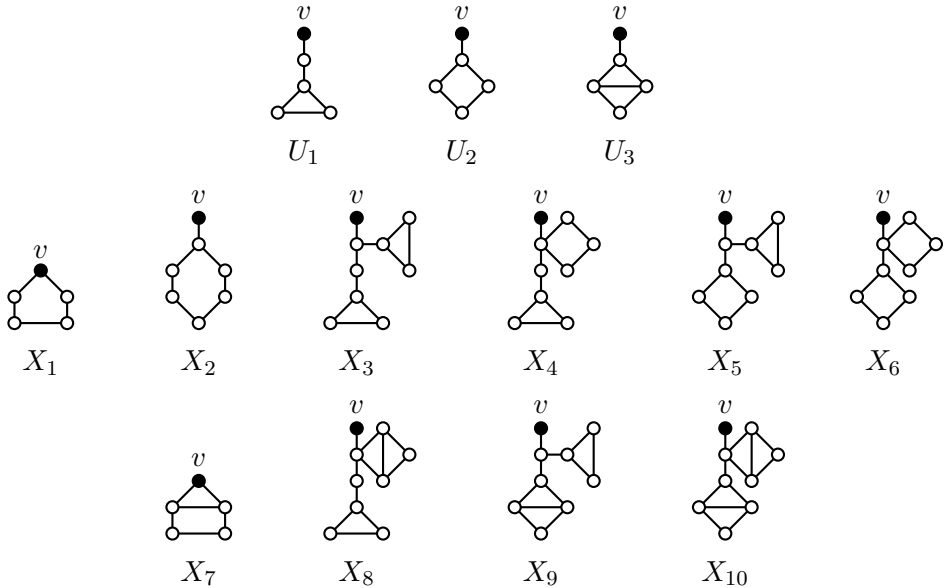

Let $U_{1}, U_{2}, U_3$ and $X_{1}, X_{2}\ldots X_{10}$ be the thirteen graphs shown in Figure~\ref{f:UV}. We define a \emph{unit} to be a graph that is isomorphic to the graph $U_i$ for some $i$, $1 \le i \le 3$, or the graph $X_j$ for some $j$, $1 \le j \le 10$. The darkened vertex, named $v$, in each unit in Figure~\ref{f:UV} we call the \emph{link vertex} of the unit.
For $i = 1,2,3$, we call a unit isomorphic to the graph $U_i$ a \emph{type}-$i$ \emph{unit}. A unit isomorphic to the graph $X_j$, $1 \le j \le 10$, we call a $X_j$-\emph{unit}.

\vskip -0.2cm


For $n = n_1 + n_2 \ge 2$, let $G = G_0(n_1,n_2)$ be the graph obtained from the disjoint union of $n_1$ units of type-$1$ and $n_2$ units of type-$2$ by identifying the $n$ link vertices, one from each unit, into one new vertex which we
call the \emph{identified vertex} of $G$. Let $\cG_0$ denote the family of all such graphs $G$. For $n = n_1 + n_2 + 1 \ge 2$ and for $i = 1,2,\ldots,6$, let $G = G_i(n_1,n_2)$ be the graph obtained from the disjoint union of $n_1$ units of type-$1$, $n_2$ units of type-$2$, and one $X_i$-unit by identifying the $n$ link vertices, one from each unit, into one new vertex which we call the \emph{identified vertex} of $G$. Let $\cG_i$ denote the family of all such graphs $G$. Let
\[
\cG = \bigcup_{i=0}^6 \cG_i.
\]

A graph in each family $\cG_i$, $0 \le i \le 6$, is shown in Figure~\ref{f:cG}, where the identified vertices are indicated in bold.

\begin{figure}[htb]
\tikzstyle{every node}=[circle, draw, fill=black!0, inner sep=0pt,minimum width=.16cm]
\begin{center}
\begin{tikzpicture}[thick,scale=.7]

\draw(5,20){


+(-2.00,1.50) node{} +(-2.00,1.00) node{} +(-2.50,0.50) node{} +(-1.50,0.50) node{}


+(1.50,2.50)--+(-2.00,1.50)--+(-2.00,1.00)--+(-2.50,0.50)--+(-1.50,0.50)--+(-2.00,1.00)


+(5.00,1.50) node{} +(5.00,1.00) node{} +(4.50,0.50) node{} +(5.50, 0.50) node{}


+(1.50,2.50)--+(5.00,1.50)--+(5.00,1.00)--+(4.50,0.50)--+(5.50,0.50)--+(5.00,1.00) 


+(1.50,2.50) node[circle, draw, fill=black!100, inner sep=0pt,minimum width=.16cm]{}


+(-0.50,1.00) node{} +(0.50,1.00) node{} +(0.00,0.50) node{} +(0.00,1.50) node{} 


+(3.00,0.50) node{} +(2.50,1.00) node{} +(3.50,1.00) node{} +(3.00,1.50) node{} 


+(1.50,1.00) node{} +(1.00,0.50) node{} +(2.00,0.50) node{} +(1.50,1.50) node{} 


+(1.50,2.50)--+(0.00,1.50) +(1.50,2.50)--+(1.50,1.50) +(1.50,2.50)--+(3.00,1.50)


+(1.50,1.50)--+(1.50,1.00)--+(1.00,0.50)--+(2.00,0.50)--+(1.50,1.00)


+(0.00,1.50)--+(-0.50,1.00)--+(0.00,0.50)--+(0.50,1.00)--+(0.00,1.50)


+(3.00,1.50)--+(2.50,1.00)--+(3.00,0.50)--+(3.50,1.00)--+(3.00,1.50)

+(1.5,-0.75) node[rectangle, draw=white!0, fill=white!100]{$G_{0}(3,2)\in\cG_{0}$}

};

\draw(-1,14){


+(1.50,2.50) node{} +(1.00,3.00) node{} +(2.00,3.00) node{} +(1.00,3.50) node{} +(2.00,3.50) node{}


+(1.50,2.50)--+(1.00,3.00)--+(1.00,3.50)--+(2.00,3.50)--+(2.00,3.00)--+(1.50,2.50)


+(1.50,2.50) node[circle, draw, fill=black!100, inner sep=0pt,minimum width=.16cm]{}


+(-0.50,1.00) node{} +(0.50,1.00) node{} +(0.00,0.50) node{} +(0.00,1.50) node{} 


+(3.00,0.50) node{} +(2.50,1.00) node{} +(3.50,1.00) node{} +(3.00,1.50) node{} 


+(1.50,1.00) node{} +(1.00,0.50) node{} +(2.00,0.50) node{} +(1.50,1.50) node{} 


+(1.50,2.50)--+(0.00,1.50) +(1.50,2.50)--+(1.50,1.50) +(1.50,2.50)--+(3.00,1.50)


+(1.50,1.50)--+(1.50,1.00)--+(1.00,0.50)--+(2.00,0.50)--+(1.50,1.00)


+(0.00,1.50)--+(-0.50,1.00)--+(0.00,0.50)--+(0.50,1.00)--+(0.00,1.50)


+(3.00,1.50)--+(2.50,1.00)--+(3.00,0.50)--+(3.50,1.00)--+(3.00,1.50)

+(1.5,-0.75) node[rectangle, draw=white!0, fill=white!100]{$G_1(1,2)\in\cG_1$}

};

\draw(5,14){


+(1.50,2.50) node{} +(1.00,3.00) node{} +(2.00,3.00) node{} +(1.00,3.50) node{} +(2.00,3.50) node{} +(1.50,4.00) node{}


+(1.50,2.50)--+(1.00,3.00)--+(1.00,3.50)--+(1.50,4.00)--+(2.00,3.50)--+(2.00,3.00)--+(1.50,2.50)


+(1.50,2.50) node {}
+(1.50,2.00) node[circle, draw, fill=black!100, inner sep=0pt,minimum width=.16cm]{}


+(-0.50,1.00) node{} +(0.50,1.00) node{} +(0.00,0.50) node{} +(0.00,1.50) node{} 


+(3.00,0.50) node{} +(2.50,1.00) node{} +(3.50,1.00) node{} +(3.00,1.50) node{} 


+(1.50,1.00) node{} +(1.00,0.50) node{} +(2.00,0.50) node{} +(1.50,1.50) node{} 


+(1.50,2.00)--+(0.00,1.50) +(1.50,2.50)--+(1.50,1.50) +(1.50,2.00)--+(3.00,1.50)


+(1.50,1.50)--+(1.50,1.00)--+(1.00,0.50)--+(2.00,0.50)--+(1.50,1.00)


+(0.00,1.50)--+(-0.50,1.00)--+(0.00,0.50)--+(0.50,1.00)--+(0.00,1.50)


+(3.00,1.50)--+(2.50,1.00)--+(3.00,0.50)--+(3.50,1.00)--+(3.00,1.50)

+(1.5,-0.75) node[rectangle, draw=white!0, fill=white!100]{$G_2(1,2)\in\cG_2$}

};

\draw(11,14){


+(1.50,3.00) node{} +(1.50,3.50) node{} +(1.00,4.00) node{} +(2.00, 4.00) node{}

+(2.00,3.00) node{} +(2.50,3.00) node{} +(3.00,3.50) node{} +(3.00,2.50) node{}


+(1.50,2.50)--+(1.50,3.50)--+(1.00,4.00)--+(2.00,4.00)--+(1.50,3.50)
+(1.50,3.00)--+(2.50,3.00)--+(3.00,3.50)--+(3.00,2.50)--+(2.50,3.00)


+(1.50,2.50) node[circle, draw, fill=black!100, inner sep=0pt,minimum width=.16cm]{}


+(-0.50,1.00) node{} +(0.50,1.00) node{} +(0.00,0.50) node{} +(0.00,1.50) node{} 


+(3.00,0.50) node{} +(2.50,1.00) node{} +(3.50,1.00) node{} +(3.00,1.50) node{} 


+(1.50,1.00) node{} +(1.00,0.50) node{} +(2.00,0.50) node{} +(1.50,1.50) node{} 


+(1.50,2.50)--+(0.00,1.50) +(1.50,2.50)--+(1.50,1.50) +(1.50,2.50)--+(3.00,1.50)


+(1.50,1.50)--+(1.50,1.00)--+(1.00,0.50)--+(2.00,0.50)--+(1.50,1.00)


+(0.00,1.50)--+(-0.50,1.00)--+(0.00,0.50)--+(0.50,1.00)--+(0.00,1.50)


+(3.00,1.50)--+(2.50,1.00)--+(3.00,0.50)--+(3.50,1.00)--+(3.00,1.50)

+(1.5,-0.75) node[rectangle, draw=white!0, fill=white!100]{$G_{3}(1,2)\in \cG_{3}$}

};

\draw(-1,8){

+(1.50,3.00) node{} +(1.50,4.00) node{} +(1.00,3.50) node{} +(2.00, 3.50) node{}

+(2.00,3.00) node{} +(2.50,3.00) node{} +(3.00,3.50) node{} +(3.00,2.50) node{}


+(1.50,2.50)--+(1.50,3.00)--+(1.00,3.50)--+(1.50,4.00)--+(2.00,3.50)--+(1.50,3.00)
+(1.50,3.00)--+(2.50,3.00)--+(3.00,3.50)--+(3.00,2.50)--+(2.50,3.00)


+(1.50,2.50) node[circle, draw, fill=black!100, inner sep=0pt,minimum width=.16cm]{}


+(-0.50,1.00) node{} +(0.50,1.00) node{} +(0.00,0.50) node{} +(0.00,1.50) node{} 


+(3.00,0.50) node{} +(2.50,1.00) node{} +(3.50,1.00) node{} +(3.00,1.50) node{} 


+(1.50,1.00) node{} +(1.00,0.50) node{} +(2.00,0.50) node{} +(1.50,1.50) node{} 


+(1.50,2.50)--+(0.00,1.50) +(1.50,2.50)--+(1.50,1.50) +(1.50,2.50)--+(3.00,1.50)


+(1.50,1.50)--+(1.50,1.00)--+(1.00,0.50)--+(2.00,0.50)--+(1.50,1.00)


+(0.00,1.50)--+(-0.50,1.00)--+(0.00,0.50)--+(0.50,1.00)--+(0.00,1.50)


+(3.00,1.50)--+(2.50,1.00)--+(3.00,0.50)--+(3.50,1.00)--+(3.00,1.50)

+(1.5,-0.75) node[rectangle, draw=white!0, fill=white!100]{$G_{4}(1,2)\in\cG_{4}$}

};

\draw(5,8){


+(1.50,3.00) node{} +(1.50,3.50) node{} +(1.00,4.00) node{} +(2.00, 4.00) node{}
+(2.00,3.00) node{} +(3.00,3.00) node{} +(2.50,3.50) node{} +(2.50,2.50) node{}


+(1.50,2.50)--+(1.50,3.50)--+(1.00,4.00)--+(2.00,4.00)--+(1.50,3.50)
+(1.50,3.00)--+(2.00,3.00)--+(2.50,2.50)--+(3.00,3.00)--+(2.50,3.50)--+(2.00,3.00)


+(1.50,2.50) node[circle, draw, fill=black!100, inner sep=0pt,minimum width=.16cm]{}


+(-0.50,1.00) node{} +(0.50,1.00) node{} +(0.00,0.50) node{} +(0.00,1.50) node{} 


+(3.00,0.50) node{} +(2.50,1.00) node{} +(3.50,1.00) node{} +(3.00,1.50) node{} 


+(1.50,1.00) node{} +(1.00,0.50) node{} +(2.00,0.50) node{} +(1.50,1.50) node{} 


+(1.50,2.50)--+(0.00,1.50) +(1.50,2.50)--+(1.50,1.50) +(1.50,2.50)--+(3.00,1.50)


+(1.50,1.50)--+(1.50,1.00)--+(1.00,0.50)--+(2.00,0.50)--+(1.50,1.00)


+(0.00,1.50)--+(-0.50,1.00)--+(0.00,0.50)--+(0.50,1.00)--+(0.00,1.50)


+(3.00,1.50)--+(2.50,1.00)--+(3.00,0.50)--+(3.50,1.00)--+(3.00,1.50)

+(1.5,-0.75) node[rectangle, draw=white!0, fill=white!100]{$G_{5}(1,2)\in\cG_{5}$}

};

\draw(11,8){

+(1.50,3.00) node{} +(1.50,4.00) node{} +(1.00,3.50) node{} +(2.00, 3.50) node{}
+(2.00,3.00) node{} +(3.00,3.00) node{} +(2.50,3.50) node{} +(2.50,2.50) node{}

+(1.50,2.50)--+(1.50,3.00)--+(1.00,3.50)--+(1.50,4.00)--+(2.00,3.50)--+(1.50,3.00)
+(1.50,3.00)--+(2.00,3.00)--+(2.50,2.50)--+(3.00,3.00)--+(2.50,3.50)--+(2.00,3.00)


+(1.50,2.50) node[circle, draw, fill=black!100, inner sep=0pt,minimum width=.16cm]{}


+(-0.50,1.00) node{} +(0.50,1.00) node{} +(0.00,0.50) node{} +(0.00,1.50) node{} 


+(3.00,0.50) node{} +(2.50,1.00) node{} +(3.50,1.00) node{} +(3.00,1.50) node{} 


+(1.50,1.00) node{} +(1.00,0.50) node{} +(2.00,0.50) node{} +(1.50,1.50) node{} 


+(1.50,2.50)--+(0.00,1.50) +(1.50,2.50)--+(1.50,1.50) +(1.50,2.50)--+(3.00,1.50)


+(1.50,1.50)--+(1.50,1.00)--+(1.00,0.50)--+(2.00,0.50)--+(1.50,1.00)


+(0.00,1.50)--+(-0.50,1.00)--+(0.00,0.50)--+(0.50,1.00)--+(0.00,1.50)


+(3.00,1.50)--+(2.50,1.00)--+(3.00,0.50)--+(3.50,1.00)--+(3.00,1.50)

+(1.5,-0.75) node[rectangle, draw=white!0, fill=white!100]{$G_{6}(1,2)\in \cG_{6}$}

};

\end{tikzpicture}
\end{center}
\vskip -0.5 cm \caption{Graphs in the family $\cG$.} \label{f:cG}
\end{figure}

For $n = n_1 + n_2 + n_3 \ge 2$, let $H = H_0(n_1,n_2,n_3)$ be the graph obtained from the disjoint union of $n_1$ units of type-$1$, $n_2$ units of type-$2$ and $n_3$ units of type-$3$ by identifying the $n$ link vertices, one from each unit, into one new vertex which we call the \emph{identified vertex} of $H$. Let $\cH_0$ denote the family of all such graphs $H$. For $n = n_1 + n_2 + n_3 + 1 \ge 2$ and for $i = 1,2,\ldots,10$, let $H = H_i(n_1,n_2,n_3)$ be the graph obtained from the disjoint union of $n_1$ units of type-$1$, $n_2$ units of type-$2$, $n_3$ units of type-$3$, and one $X_i$-unit by identifying the $n$ link vertices, one from each unit, into one new vertex which we call the \emph{identified vertex} of $H$. Let $\cH_i$ denote the family of all such graphs $H$. Let
\[
\cH = \bigcup_{i=0}^{10} \cH_i.
\]

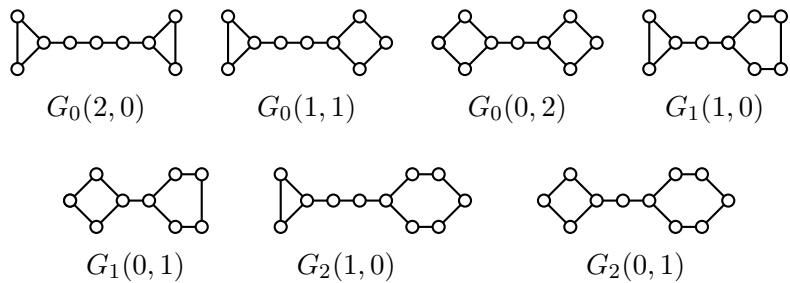
\begin{figure}[htb]
\tikzstyle{every node}=[circle, draw, fill=black!0, inner sep=0pt,minimum width=.16cm]
\begin{center}
\begin{tikzpicture}[thick,scale=.7]

\draw(0,18){+(0.00,0.00)--+(0.00,1.00)--+(0.50,0.50)--+(0.00,0.00) +(0.50,0.50)--+(2.50,0.50) +(2.50,0.50)--+(3.00,1.00)--+(3.00,0.00)--+(2.50,0.50)
+(0.00,0.00)node{} +(0.00,1.00)node{} +(0.50,0.50)node{} +(1.00,0.50)node{} +(1.50,0.50)node{} +(2.00,0.50)node{} +(2.50,0.50)node{} +(3.00,0.00)node{} +(3.00,1.00)node{}
+(1.5,-0.75) node[rectangle, draw=white!0, fill=white!100]{$G_{0}(2,0)$}
};

\draw(4,18){+(0.00,0.00)--+(0.00,1.00)--+(0.50,0.50)--+(0.00,0.00) +(0.50,0.50)--+(2.00,0.50) +(2.00,0.50)--+(2.50,1.00)--+(3.00,0.50)--+(2.50,0.00)--+(2.00,0.50)
+(0.00,0.00)node{} +(0.00,1.00)node{} +(0.50,0.50)node{} +(1.00,0.50)node{} +(1.50,0.50)node{} +(2.00,0.50)node{} +(2.50,1.00) node{} +(3.00,0.50) node{} +(2.50,0.00) node{}
+(1.5,-0.75) node[rectangle, draw=white!0, fill=white!100]{$G_{0}(1,1)$}
};

\draw(8,18){+(0.00,0.50)--+(0.50,1.00)--+(1.00,0.50)--+(0.50,0.00)--+(0.00,0.50) +(1.00,0.50)--+(2.00,0.50) +(2.00,0.50)--+(2.50,1.00)--+(3.00,0.50)--+(2.50,0.00)--+(2.00,0.50)
+(0.00,0.50) node{} +(0.50,1.00) node{} +(1.00,0.50) node{} +(0.50,0.00) node{} +(0.00,0.50) node{} +(1.50,0.50) node{} +(2.00,0.50) node{} +(3.00,0.50) node{} +(2.50,1.00) node{} +(2.50,0.00) node{}
+(1.5,-0.75) node[rectangle, draw=white!0, fill=white!100]{$G_{0}(0,2)$}
};

\draw(12,18){+(0.00,0.00)--+(0.00,1.00)--+(0.50,0.50)--+(0.00,0.00) +(0.50,0.50)--+(1.50,0.50) +(1.50,0.50)--+(2.00,1.00)--+(2.50,1.00)--+(2.50,0.00)--+(2.00,0.00)--+(1.50,0.50)
+(0.00,0.00)node{} +(0.00,1.00)node{} +(0.50,0.50)node{} +(1.00,0.50)node{} +(1.50,0.50)node{} +(2.00,1.00) node{} +(2.50,1.00) node{} +(2.50,0.00) node{} +(2.00,0.00) node{} +(1.50,0.50) node{}
+(1.25,-0.75) node[rectangle, draw=white!0, fill=white!100]{$G_1(1,0)$}
};

\draw(1,15){+(0.00,0.50)--+(0.50,1.00)--+(1.00,0.50)--+(0.50,0.00)--+(0.00,0.50) +(1.00,0.50)--+(1.50,0.50)--+(2.00,1.00)--+(2.50,1.00)--+(2.50,0.00)--+(2.00,0.00)--+(1.50,0.50)
+(0.00,0.50) node{} +(0.50,1.00) node{} +(1.00,0.50) node{} +(0.50,0.00) node{} +(0.00,0.50) node{} +(1.50,0.50) node{} +(2.00,1.00) node{} +(2.50,1.00) node{} +(2.50,0.00) node{} +(2.00,0.00) node{}
+(1.25,-0.75) node[rectangle, draw=white!0, fill=white!100]{$G_1(0,1)$}
};

\draw(5,15){+(0.00,0.00)--+(0.00,1.00)--+(0.50,0.50)--+(0.00,0.00) +(0.50,0.50)--+(2.00,0.50) +(2.00,0.50)--+(2.50,1.00)--+(3.00,1.00)--+(3.50,0.50)--+(3.00,0.00)--+(2.50,0.00)--+(2.00,0.50)
+(0.00,0.00)node{} +(0.00,1.00)node{} +(0.50,0.50)node{} +(1.00,0.50)node{} +(1.50,0.50)node{} +(2.00,0.50) node{} +(3.50,0.50) node{} +(3.00,1.00) node{} +(3.00,0.00) node{} +(2.50,0.00) node{} +(2.50,1.00) node{}
+(1.25,-0.75) node[rectangle, draw=white!0, fill=white!100]{$G_2(1,0)$}
};

\draw(10,15){+(0.00,0.50)--+(0.50,1.00)--+(1.00,0.50)--+(0.50,0.00)--+(0.00,0.50) +(1.00,0.50)--+(2.00,0.50)  +(2.00,0.50)--+(2.50,1.00)--+(3.00,1.00)--+(3.50,0.50)--+(3.00,0.00)--+(2.50,0.00)--+(2.00,0.50)
+(0.00,0.50) node{} +(1.00,0.50) node{} +(0.50,0.00) node{} +(0.50,1.00) node{} +(1.50,0.50) node{} +(2.00,0.50) node{} +(3.50,0.50) node{} +(3.00,1.00) node{} +(3.00,0.00) node{} +(2.50,0.00) node{} +(2.50,1.00) node{}
+(1.75,-0.75) node[rectangle, draw=white!0, fill=white!100]{$G_2(0,1)$}
};

\end{tikzpicture}
\end{center}
\vskip -0.5 cm \caption{The family $\cG_{b}\subset \cG$ of dumb-bell graphs.} \label{f:Gb}
\end{figure}

We observe that $\cG$ is a subfamily of $\cH$; that is, $\cG \subset \cH$. Let $\cG_b$ be the subfamily of $\cG$ consisting of all dumb-bells. We observe that
\[
\cG_b = \{ G_0(2,0), G_0(1,1), G_0(0,2), G_1(1,0), G_1(0,1), G_2(1,0), G_2(0,1) \}.
\]
The family $\cG_b\subset \cG$ of seven dumb-bell graphs is shown in Figure~\ref{f:Gb}.

Let $\cB=\{B_{1}, B_2, \ldots, B_{11}\}$ be the family of eleven graphs shown in Figure \ref{f:cB}, where the special vertices (Section~\ref{NDTDS}) of each graph are darkened.
\begin{figure}[htb]
\tikzstyle{every node}=[circle, draw, fill=black!0, inner sep=0pt,minimum width=.16cm]
\begin{center}
\begin{tikzpicture}[thick,scale=.7]

\draw(-1,14.5) {+(0.00,0.50)--+(1.00,0.50)--+(0.50,1.00)--+(0.00,0.50)--+(0.50,0.00)--+(1.00,0.50) 
+(0.00,0.50)node[circle, draw, fill=black!100, inner sep=0pt,minimum width=.16cm]{} +(0.50,0.50) node{} +(1.00,0.50)node[circle, draw, fill=black!100, inner sep=0pt,minimum width=.16cm]{} +(0.50,1.00)node{} +(0.50,0.00)node{}
+(0.50,-0.75) node[rectangle, draw=white!0, fill=white!100]{$B_{1}$}
};

\draw(3,14.5) {+(0.00,0.50)--+(0.50,1.00)--+(1.50,1.00)--+(1.00,0.50)--+(1.50,0.00)--+(2.00,0.50)--+(1.50,1.00) +(0.00,0.50)--+(0.50,0.00)--+(1.50,0.00) +(0.00,0.50)node[circle, draw, fill=black!100, inner sep=0pt,minimum width=.16cm]{} +(0.50,1.00) node{} +(1.50,1.00) node[circle, draw, fill=black!100, inner sep=0pt,minimum width=.16cm]{} +(1.00,0.50) node{} +(1.50,0.00) node[circle, draw, fill=black!100, inner sep=0pt,minimum width=.16cm]{} +(2.00,0.50) node{} +(0.50,0.00) node{}
+(1.00,-0.75) node[rectangle, draw=white!0, fill=white!100]{$B_2$}
};

\draw(8,14.5) {+(0.00,0.00)--+(0.00,1.00)--+(2.00,1.00)--+(2.50,0.50)--+(2.00,0.00)--+(0.00,0.00) +(2.00,1.00)--+(1.50,0.50)--+(2.00,0.00) +(0.00,0.00) node[circle, draw, fill=black!100, inner sep=0pt,minimum width=.16cm]{} +(0.00,1.00) node[circle, draw, fill=black!100, inner sep=0pt,minimum width=.16cm]{} +(1.00,0.00) node{} +(1.00,1.00) node{} +(2.00,1.00) node[circle, draw, fill=black!100, inner sep=0pt,minimum width=.16cm]{} +(2.00,0.00) node[circle, draw, fill=black!100, inner sep=0pt,minimum width=.16cm]{} +(2.50,0.50) node{} +(1.50,0.50) node{}
+(1.00,-0.75) node[rectangle, draw=white!0, fill=white!100]{$B_3$}
};

\draw(12.5,14.5){+(0.00,0.00)--+(0.00,1.00)--+(2.00,1.00)--+(2.50,0.50)--+(2.00,0.00)--+(0.00,0.00) +(2.00,1.00)--+(1.50,0.50)--+(2.00,0.00) +(0.00,0.00) node[circle, draw, fill=black!100, inner sep=0pt,minimum width=.16cm]{} +(0.00,1.00) node[circle, draw, fill=black!100, inner sep=0pt,minimum width=.16cm]{} +(1.00,0.00) node{} +(1.00,1.00) node{} +(2.00,1.00) node[circle, draw, fill=black!100, inner sep=0pt,minimum width=.16cm]{} +(2.00,0.00) node[circle, draw, fill=black!100, inner sep=0pt,minimum width=.16cm]{} +(2.50,0.50) node{} +(1.50,0.50) node{} +(2.00,0.50) node{} +(2.00,1.00)--+(2.00,0.00)
+(1.00,-0.75) node[rectangle, draw=white!0, fill=white!100]{$B_{4}$}
};

\draw(1,12) {+(-0.50,0.50)--+(0.00,0.00) +(-0.50,0.50)--+(0.00,1.00)
+(0.00,1.00)--+(2.00,1.00)--+(2.50,0.50)--+(2.00,0.00)--+(0.00,0.00) +(2.00,1.00)--+(1.50,0.50)--+(2.00,0.00) +(0.00,0.00) node[circle, draw, fill=black!100, inner sep=0pt,minimum width=.16cm]{} +(0.00,1.00) node[circle, draw, fill=black!100, inner sep=0pt,minimum width=.16cm]{} +(1.00,0.00) node{} +(1.00,1.00) node{} +(2.00,1.00) node[circle, draw, fill=black!100, inner sep=0pt,minimum width=.16cm]{} +(2.00,0.00) node[circle, draw, fill=black!100, inner sep=0pt,minimum width=.16cm]{} +(2.50,0.50) node{} +(1.50,0.50) node{} +(-0.50,0.50) node[circle, draw, fill=black!100, inner sep=0pt,minimum width=.16cm]{}
+(1.00,-0.75) node[rectangle, draw=white!0, fill=white!100]{$B_{5}$}
};

\draw(5.5,12) {+(0.00,0.50)--+(0.50,1.00)--+(2.50,1.00)--+(3.00,0.50)--+(2.50,0.00)--+(0.50,0.00)--+(0.00,0.50) +(0.00,0.50)--+(3.00,0.50)
+(0.00,0.50) node[circle, draw, fill=black!100, inner sep=0pt,minimum width=.16cm]{} +(0.50,1.00) node{} +(2.50,1.00) node{} +(3.00,0.50) node[circle, draw, fill=black!100, inner sep=0pt,minimum width=.16cm]{} +(2.50,0.00) node{} +(0.50,0.00) node{} +(0.75,0.50) node{} +(1.5,0.50) node[circle, draw, fill=black!100, inner sep=0pt,minimum width=.16cm]{}
+(2.25,0.50) node{}
+(1.50,-0.75) node[rectangle, draw=white!0, fill=white!100]{$B_{6}$}
};

\draw(10.5,12) {+(0.00,0.50)--+(0.50,1.00)--+(2.50,1.00)--+(3.00,0.50)--+(2.50,0.00)--+(0.50,0.00)--+(0.00,0.50) 
+(0.00,0.50) node[circle, draw, fill=black!100, inner sep=0pt,minimum width=.16cm]{} +(0.50,1.00) node{} +(2.50,1.00) node{} +(3.00,0.50) node[circle, draw, fill=black!100, inner sep=0pt,minimum width=.16cm]{} +(2.50,0.00) node{} +(0.50,0.00) node{} +(1.16666,0.00) node[circle, draw, fill=black!100, inner sep=0pt,minimum width=.16cm]{}
+(1.8333,0.00) node[circle, draw, fill=black!100, inner sep=0pt,minimum width=.16cm]{} +(1.50,0.50) node{} +(0.00,0.50)--+(1.50,0.50)--+(3.00,0.50)
+(1.50,-0.75) node[rectangle, draw=white!0, fill=white!100]{$B_{7}$}
};

\draw(-1.0,9.5) {+(0.00,0.50)--+(0.50,0.00)--+(1.50,0.00)--+(2.00,0.50)--+(0.00,0.50)--+(0.5,1.00)--+(1.00,0.5)--+(1.5,1.00)--+(2.00,0.50) +(0.50,0.00) node{} +(1.50,0.00) node{} +(0.00,0.50) node[circle, draw, fill=black!100, inner sep=0pt,minimum width=.16cm]{} +(2.00,0.50) node[circle, draw, fill=black!100, inner sep=0pt,minimum width=.16cm]{} +(0.5,1.00) node{} +(1.00,0.5) node[circle, draw, fill=black!100, inner sep=0pt,minimum width=.16cm]{} +(0.5,0.5) node{} +(1.5,0.50) node{} +(1.5,1.00) node{}
+(1.00,-0.75) node[rectangle, draw=white!0, fill=white!100]{$B_{8}$}
};

\draw(3,9.5) {+(0.50,0.50)--+(0.00,1.00)--+(0.00,0.00)--+(0.50,0.50) +(0.50,0.50)--+(1.50,0.50) +(1.50,0.50)--+(2.00,1.00)--+(2.50,0.50)--+(2.00,0.00)--+(1.50,0.50) +(1.50,0.50)--+(2.50,0.50)
+(2.00,0.50) node{}
+(0.00,0.00) node{} +(0.00,1.00) node{} +(0.50,0.50) node[circle, draw, fill=black!100, inner sep=0pt,minimum width=.16cm]{}
+(1.00,0.50) node[circle, draw, fill=black!100, inner sep=0pt,minimum width=.16cm]{}
+(1.50,0.50) node[circle, draw, fill=black!100, inner sep=0pt,minimum width=.16cm]{} +(2.00,1.00) node{} +(2.50,0.50) node[circle, draw, fill=black!100, inner sep=0pt,minimum width=.16cm]{} +(2.00,0.00) node{}
+(1.25,-0.75) node[rectangle, draw=white!0, fill=white!100]{$B_{9}$}
};

\draw(8,9.5) {+(0.00,0.50)--+(0.50,0.00)--+(1.00,0.50)--+(0.50,1.00)--+(0.00,0.50) +(1.00,0.50)--+(2.50,0.50) +(1.00,0.50)--+(1.50,0.50)--+(2.00,0.00)--+(2.50,0.50)--+(2.00,1.00)--+(1.50,0.50) +(2.00,0.50)node{} +(1.00,0.50)node[circle, draw, fill=black!100, inner sep=0pt,minimum width=.16cm]{} +(0.00,0.50)node[circle, draw, fill=black!100, inner sep=0pt,minimum width=.16cm]{} +(0.50,0.00)node{} +(0.50,1.00)node{} +(1.50,0.50)node[circle, draw, fill=black!100, inner sep=0pt,minimum width=.16cm]{} +(2.50,0.50)node[circle, draw, fill=black!100, inner sep=0pt,minimum width=.16cm]{} +(2.00,0.00)node{} +(2.00,1.00)node{}
+(1.25,-0.75) node[rectangle, draw=white!0, fill=white!100]{$B_{10}$}
};

\draw(12,9.5){+(0.00,0.50)--+(0.50,1.00)--+(2.50,1.00)--+(3.00,0.50)--+(2.50,0.00)--+(0.50,0.00)--+(0.00,0.50) +(2.50,0.00)--+(2.00,0.50)--+(2.50,1.00)
+(0.00,0.50) node[circle, draw, fill=black!100, inner sep=0pt,minimum width=.16cm]{} +(0.50,1.00) node[circle, draw, fill=black!100, inner sep=0pt,minimum width=.16cm]{} +(1.00,1.00) node[circle, draw, fill=black!100, inner sep=0pt,minimum width=.16cm]{} +(1.50,1.00) node[circle, draw, fill=black!100, inner sep=0pt,minimum width=.16cm]{} +(2.00,1.00) node{} +(2.50,1.00) node[circle, draw, fill=black!100, inner sep=0pt,minimum width=.16cm]{} +(3.00,0.50) node{}
+(0.50,0.00) node[circle, draw, fill=black!100, inner sep=0pt,minimum width=.16cm]{} +(1.00,0.00) node[circle, draw, fill=black!100, inner sep=0pt,minimum width=.16cm]{} +(1.50,0.00) node[circle, draw, fill=black!100, inner sep=0pt,minimum width=.16cm]{} +(2.00,0.00) node{} +(2.50,0.00) node[circle, draw, fill=black!100, inner sep=0pt,minimum width=.16cm]{} +(2.00,0.50) node{}
+(1.5,-0.75) node[rectangle, draw=white!0, fill=white!100]{$B_{11}$}
};

\end{tikzpicture}
\end{center}
\vskip -0.5 cm \caption{The family, $\cB$, of graphs.} \label{f:cB}
\end{figure}

\newpage
Let $\cF=\{F_{1}, F_{2}, \ldots, F_{6}\}$ be the family of six graphs shown in Figure~\ref{f:cF}.

\vspace{2mm}
\begin{figure}[htb]
\tikzstyle{every node}=[circle, draw, fill=black!0, inner sep=0pt,minimum width=.16cm]
\begin{center}
\begin{tikzpicture}[thick,scale=.7]
\draw(0,0){+(0.00,0.50)--+(0.50,1.00)--+(1.00,0.50)--+(0.50,0.00)--+(0.00,0.50) +(0.50,1.00)--+(0.50,0.00)
+(1.00,0.50)--+(1.50,0.50)--+(2.00,1.00)--+(2.50,0.50)--+(2.00,0.00)--+(1.50,0.50)
+(0.00,0.50) node{} +(0.50,1.00) node{} +(1.00,0.50) node{} +(0.50,0.00) node{} +(1.00,0.50) node{} +(2.00,1.00) node{} +(2.00,0.00) node{} +(2.50,0.50) node{} +(1.50,0.50) node {}
+(1.25,-0.75) node[rectangle, draw=white!0, fill=white!100]{$F_{1}$}
};
\draw(3.5,0){+(0.00,0.50)--+(0.50,1.00)--+(1.00,0.50)--+(0.50,0.00)--+(0.00,0.50) +(0.50,1.00)--+(0.50,0.00) +(2.00,1.00)--+(2.00,0.00)
+(1.00,0.50)--+(1.50,0.50)--+(2.00,1.00)--+(2.50,0.50)--+(2.00,0.00)--+(1.50,0.50)
+(0.00,0.50) node{} +(0.50,1.00) node{} +(1.00,0.50) node{} +(0.50,0.00) node{} +(1.00,0.50) node{} +(2.00,1.00) node{} +(2.00,0.00) node{} +(2.50,0.50) node{} +(1.50,0.50) node {}
+(1.25,-0.75) node[rectangle, draw=white!0, fill=white!100]{$F_{2}$}
};
\draw(7,0){+(0.50,0.50)--+(0.00,1.00)--+(0.00,0.00)--+(0.50,0.50) +(0.50,0.50)--+(1.50,0.50) +(1.50,0.50)--+(2.00,1.00)--+(2.50,0.50)--+(2.00,0.00)--+(1.50,0.50) +(2.00,1.00)--+(2.00,0.00)
+(0.00,0.00) node{} +(0.00,1.00) node{} +(0.50,0.50) node{}
+(1.00,0.50) node{}
+(1.50,0.50) node{} +(2.00,1.00) node{} +(2.50,0.50) node{} +(2.00,0.00) node{}
+(1.00,-0.75) node[rectangle, draw=white!0, fill=white!100]{$F_{3}$}
};
\draw(10.5,0){+(0.00,0.50)--+(0.50,1.00)--+(2.00,1.00)--+(2.50,0.50)--+(2.00,0.00)--+(0.50,0.00)--+(0.00,0.50) +(0.50,0.00)--+(0.50,1.00)
+(0.00,0.50) node{} +(0.50,1.00) node{} +(2.00,1.00) node{} +(1.25,1.00) node{} +(2.50,0.50) node{} +(2.00,0.00) node{} +(1.25, 0.00) node{} +(0.50, 0.00) node{}
+(1.25,-0.75) node[rectangle, draw=white!0, fill=white!100]{$F_{4}$}
};
\draw(14,0){+(0.00,0.50)--+(0.50,1.00)--+(2.00,1.00)--+(2.50,0.50)--+(2.00,0.00)--+(0.50,0.00)--+(0.00,0.50) +(0.50,0.00)--+(0.50,1.00) +(0.00,0.50)--+(1.25,1.00)
+(0.00,0.50) node{} +(0.50,1.00) node{} +(2.00,1.00) node{} +(1.25,1.00) node{} +(2.50,0.50) node{} +(2.00,0.00) node{} +(1.25, 0.00) node{} +(0.50, 0.00) node{}
+(1.25,-0.75) node[rectangle, draw=white!0, fill=white!100]{$F_{5}$}
};
\draw(17.5,0) {+(0.00,0.00)--+(0.00,1.00)--+(2.00,1.00)--+(2.50,0.50)--+(2.00,0.00)--+(0.00,0.00) +(2.00,1.00)--+(1.50,0.50)--+(2.00,0.00) +(0.00,0.00) node{} +(0.00,1.00) node{} +(1.00,0.00) node{} +(1.00,1.00) node{} +(2.00,1.00) node{} +(2.00,0.00) node{} +(2.50,0.50) node{} +(1.50,0.50) node{} +(1.50,0.50)--+(2.50,0.50)
+(1.00,-0.75) node[rectangle, draw=white!0, fill=white!100]{$F_{6}$}
};
\end{tikzpicture}
\end{center}
\vskip -0.5 cm \caption{The family, $\cF$, of graphs.} \label{f:cF}
\end{figure}

\section{Main Results}

Our aim in this paper is to improve the upper bound of Theorem~\ref{1:thm} on the disjunctive total domination number of a connected graph  when we impose a density condition by restricting the minimum degree to be at least~$2$. In this case, we show that the result of~\cite{He00} given in Table~1 on the total domination number can be improved significantly for the disjunctive total domination number. Several authors (see~\cite{Alfewy04,ChMc,Tu90}) showed that if $G$ is a graph of order~$n$ with $\delta(G) \ge 3$, then $\gt(G) \le n/2$ (see Table~1). We prove that if we relax the minimum degree condition from~$\delta(G) \ge 3$ to~$\delta(G) \ge 2$, then $n/2$ is an upper bound on the disjunctive total domination number $\dtd(G)$, provided $n \ge 8$. More precisely, we prove the following result. A proof of Theorem~\ref{n2:thm} is given in Section~\ref{n2:subsec}.

\begin{thm}
\label{n2:thm}
Let $G$ be a connected graph of order $n \ge 8$ with $\delta(G)\ge 2$. Then, $\dtd(G)\le n/2$ with equality if and only if $G \in \{C_{8}, C_{12}, B_3, D_{b}(4,4), D_{b}(3,4,1), D_{b}(3,3,2)\} \cup \cF$.
\end{thm}

The connected graphs with minimum degree at least~$2$ and order at least~$18$ that have maximum possible disjunctive total domination number are characterized in the following result, a proof of which is given in Section~\ref{nminus:subsec}.

\begin{thm}
\label{nminus:thm}
Let $G$ be a connected graph of order $n \ge 18$ with $\delta(G)\ge 2$. Then, $\dtd(G)\le (n-1)/2$ with equality if and only if $G \in\cH$.
\end{thm}

Since the graphs that achieve equality in the upper bound of Theorem~\ref{n2:thm} all have order at most~$12$, we remark that Theorem~\ref{t:main1} is an immediate consequence of Theorem~\ref{n2:thm} and Theorem~\ref{nminus:thm}. Further, we remark that the upper bound of Theorem~\ref{t:main1} is sharp even for the class of bipartite graphs, for example, for $k\ge 0$ each graph $G_{0}(0,k)\in \cH$ is bipartite.

\section{Preliminary Results}
\label{pre:subsec}

Before presenting a proof of our main results, we first establish some preliminary results. We omit the proofs of these preliminary results which are straightforward, albeit tedious, to check.\footnote{Proof of several of the preliminary results can be found in the Appendix}

\begin{ob}
\label{sd:ob}
Let $G$ be a graph with no isolated vertex and let $F$ be a $4$-subdivision of $G$. Then, $\dtd(F) \le \dtd(G)+2$.
\end{ob}

The disjunctive total domination number of a cycle $C_n$ on $n$ vertices is established in~\cite{HeN13}.

\begin{prop}{\rm (\cite{HeN13})}
\label{cycle:ob}
For $n\ge 3$, $\dtd(C_{n})=2n/5$ if $n \equiv 0 \, (\mod \, 5)$ and $\dtd(C_{n})=\lceil 2(n+1)/5\rceil$ otherwise.
\end{prop}

The daisies and dumb-bells with large disjunctive total domination number are characterized in Proposition~\ref{d:ob} and Proposition~\ref{db:ob}, respectively.

\begin{prop}
\label{d:ob}
If $G$ is a daisy of order $n$, then $\dtd(G)\le (n-1)/2$. Furthermore, $\dtd(G)=(n-1)/2$ if and only if $G \in  \cD$.
\end{prop}

\begin{prop}
\label{db:ob}
If $G$ is a dumb-bell of order $n$, then $\dtd(G)\le n/2$. Furthermore, $\dtd(G)\ge (n-1)/2$ if and only if $G \in \cD_{b} \cup \cG_{b}$.
\end{prop}

The disjunctive total domination number of graphs in the family $\cB \cup \cC \cup \cD \cup \cD_{b} \cup \cG$ is given by Observation~\ref{f1:ob}.

\begin{ob}
\label{f1:ob}
Let $G \in \cB \cup \cC \cup \cD \cup \cD_{b} \cup \cG$ have order $n$. Then, $G$ is a connected graph with $\delta(G)=2$. Further,
$\dtd(G)=(n+1)/2$ if $G\in\{C_{3}, C_{7}\}$, $\dtd(G)=n/2$ if $G \in \{C_{4}, C_{6}, C_{8}, C_{12}, B_3, D_{b}(4,4), D_{b}(3,4,1), D_{b}(3,3,2)\}$ and $\dtd(G)=(n-1)/2$ otherwise.
\end{ob}

\subsection{$\frac{1}{2}$-Minimal Graphs}

In order to prove our two main results, we study so-called $\frac{1}{2}$-minimal graphs. A graph $G$ is $\frac{1}{2}$-\emph{minimal graph} if $G$ is edge-minimal with respect to the following three conditions: (i)~$\delta(G)\ge 2$, (ii) $G$ is connected, and (iii) $\dtd(G)\ge (n-1)/2$,
where $n$ is the order of $G$. If $G$ is edge-minimal with respect to conditions (i) and (ii) but does not necessarily satisfy condition (iii) we refer to $G$ as an \emph{edge-minimal graph}.
Thus if $G$ is an edge-minimal graph and $e \in E(G)$, then $\delta(G - e) = 1$ or $G - e$ is disconnected. It is evident that, an edge-minimal graph which satisfies condition (iii) is a $\frac{1}{2}$-minimal graph.

As an immediate consequence of Proposition~\ref{cycle:ob},~\ref{d:ob} and~\ref{db:ob}, we obtain a characterization of the $\frac{1}{2}$-minimal graphs that are cycles, daisies and dumb-bells.

\begin{cor}
\label{cycle:cor}
Let $G$ be a $\frac{1}{2}$-minimal graph. Then the following holds. \\
\indent {\rm (a)} $G$ is a cycle if and only if $G \in \cC$. \\
\indent {\rm (b)} $G$ is a daisy if and only if $G \in \cD$. \\
\indent {\rm (c)} $G$ is a dumb-bell if and only if $G \in \cD_{b} \cup \cG_{b}$.
\end{cor}

We next establish properties of graphs in the family $\cB \cup \cC \cup \cD \cup \cD_{b} \cup \cG$.

\begin{ob}
\label{f2:ob}
Each graph in $\cB \cup \cC \cup \cD \cup \cD_{b} \cup \cG$ is a $\frac{1}{2}$-minimal graph.
\end{ob}

Recall that an edge $e$ in a graph $G$ is a good-edge of $G$ if there is a $\dtd(G)$-set which contains both ends of $e$; otherwise, it is a bad-edge.

\begin{ob}
\label{f3a:ob}
Let $G \in \cB \cup \cC \cup \cD \cup \cD_{b} \cup \cG$. Then every edge of $G$ is a good-edge of $G$ unless $G \in \{D(4,4), D_{b}(3,4), D_{b}(3,3,1)\} \cup \cG\setminus\cG_1$. Furthermore, if $e$ is a bad-edge of $G$, then the following holds. \\
\indent {\rm (a)} If $G=D(4,4)$, then $e$ is not incident with the vertex of degree~$4$ in $G$. \\
\indent {\rm (b)} If $G=D_{b}(3,4)$, then $e$ is not incident with a vertex of degree~$3$ in $G$. \\
\indent {\rm (c)} If $G=D_{b}(3,3,1)$, then $e$ is a cycle edge of $G$. \\
\indent {\rm (d)} If $G \in \cG\setminus\cG_1$, then $e$ is incident with the identified vertex of $G$.
\end{ob}

Recall that a vertex $v$ in a graph $G$ is a  good-vertex of $G$ if it belongs to some $\dtd(G)$-set; otherwise, it is a bad-vertex of $G$.

\begin{ob}
\label{f3b:ob}
Let $G \in \cB \cup \cC \cup \cD \cup \cD_{b} \cup \cG$. Then every vertex of $G$ is a good-vertex of $G$ unless $G \in \{D(4,4), D_{b}(3,4), D_{b}(3,3,1)\} \cup \cG_0$. Furthermore, if $v$ is a bad-vertex of $G$, then the following holds. \\
\indent {\rm (a)} If $G \in \cG_{0}$, then $v$ is the identified vertex of $G$. \\
\indent {\rm (b)} If $G \in \{D(4,4), D_{b}(3,4), D_{b}(3,3,1)\}$, then $v$ is a vertex at distance~$2$ from the central \\ \hspace*{0.8cm} vertex of $G$.
\end{ob}

The following observation establishes a property of good-graphs which are not dumb-bells.

\begin{ob}
\label{f4:ob}
If $G \in \cG$ and $v$ is the identified vertex of $G$, then there is a $\dtd(G)$-set, $S$, such that $v \notin S$ and $N(v) \subset S$. Further if $G \notin \cG_1$, then $S$ can be chosen to totally dominate $N(v)$.
\end{ob}

The following observation characterizes the $\frac{1}{2}$-minimal graphs of small order.

\begin{ob}
\label{s:ob}
If $G$ is a $\frac{1}{2}$-minimal graph of order $n$, $3 \le n \le 7$, then $G \in \{B_{1}, B_2, C_3$, $C_4, C_5, C_6, D(3,3), D(4,4), D_b(3,4), D_{b}(3,3,1)\} \subset \cB \cup \cC \cup \cD \cup  \cD_{b}$.
\end{ob}

\subsection{Near Disjunctive Total Dominating Sets}
\label{NDTDS}

Let $G$ be a graph and let $v$ be a vertex of $G$. We denote the graph obtained from $G$ by adding a new vertex $v'$ and adding the pendant edge $vv'$ by $G^v$. We define a \emph{near-disjunctive total dominating set}, abbreviated NDTD-set, of $G^v$ to be a set $S$ of vertices in $G^v$ such that $v' \in S$ and every vertex in $V(G)$ is DT-dominated by the set $S$ in $G^v$. We observe that if $S$ is a NDTD-set of $G^v$, then possibly the vertex $v'$ may not be DT-dominated by $S$. The \emph{near-disjunctive total domination number}, $\ndtd(G^v)$, of $G^v$ is the minimum cardinality of a NDTD-set in $G^v$. A NDTD-set of cardinality $\ndtd(G^v)$ is called a $\ndtd(G^v)$-set.

We say that a vertex $v$ of an edge-minimal graph $G$ is \emph{special} if for every edge $e = uv$ incident with $v$ we have that $G - e$ is disconnected or $d_G(u) = 2$. In particular, we note that if $v$ is a special-vertex of $G$, then $\delta(G - v) = 1$ or $G - v$ is disconnected. The special vertices of each graph in the family $\cB$ are indicated in Figure~\ref{f:cB}.

We present next two observations that establish useful properties of graphs that belong to the family $\cB \cup \cC \cup \cD \cup \cD_{b} \cup \cG$.

\begin{ob}
\label{special:ob}
Let $G \in \cB \cup \cC \cup \cD \cup \cD_{b} \cup \cG$ and let $v$ be a special vertex of $G$. Then, $\ndtd(G^v) \le \dtd(G)$, unless one of the following four conditions hold. \\
\indent {\rm (a)} $G \in \{C_4, C_5\}$.  \\
\indent {\rm (b)} $G = B_1$ and $v$ has degree~$3$ in $G$.  \\
\indent {\rm (c)} $G \in \cG$ and $v$ is the identified vertex of $G$. \\
\indent {\rm (d)} $G \in \cG_0$ and $v$ is a neighbor of the identified vertex of $G$.
\end{ob}


\begin{ob}
\label{ndtd1:ob}
Let $G \in \cB \cup \cC \cup \cD \cup \cD_{b} \cup \cG$ and let $v$ be a special vertex of $G$. Let $v'$ be the vertex added to $G$ when constructing $G^v$. If $\ndtd(G^v) \le \dtd(G)$, then there exists a $\ndtd(G^v)$-set that DT-dominates $v'$ unless one of the following two conditions hold.
 \\
\indent {\rm (a)} $G = D(4,4)$ and $v$ is at distance~2 from the vertex of degree~$4$ in $G$.\\
\indent {\rm (b)} $G \in\{D_{b}(3,4), D_{b}(3,3,1)\}$, $d_G(v) = 3$, and $v$ belongs to a $3$-cycle in $G$.
\end{ob}

\section{Proof of Main Results}

In this section, we present a proof of our main results, namely Theorem~\ref{n2:thm} and Theorem~\ref{nminus:thm}. We begin with a characterization of $\frac{1}{2}$-minimal graphs.

\subsection{A Characterization of $\frac{1}{2}$-Minimal Graphs}

A key result to enable us to prove our main results is the following characterization of $\frac{1}{2}$-minimal graphs.

\begin{thm}
\label{em:thm}
A graph $G$ is $\frac{1}{2}$-minimal if and only if $G \in \cB \cup \cC \cup \cD \cup \cD_{b} \cup \cG$.
\end{thm}
\noindent{\bf Proof of Theorem 19.} The sufficiency follows from Observation \ref{f2:ob}. To prove the necessity, we proceed by induction on the order $n \ge 3$ of a $\frac{1}{2}$-minimal graph. By Observation~\ref{s:ob}, the result is true for $n \le 7$. Suppose $n\ge 8$, and assume that the result is true for all $\frac{1}{2}$-minimal graphs $G'$ of order $n'$, where $3 \le n' < n$. Let $G=(V,E)$ be a $\frac{1}{2}$-minimal graph of order $n$. We first present two useful observations. If $e$ is an edge of $G$, then $\dtd(G-e) \ge \dtd(G)$. Hence, by the minimality of $G$, we have the following observation.

\begin{ob}
\label{bridge:ob}
If $e\in E$, then either $e$ is a bridge of $G$ or $\delta(G-e)=1$.
\end{ob}

Since the disjunctive total domination number cannot decrease if edges are removed, the next observation follows as a consequence of the inductive hypothesis.

\begin{ob}
\label{sg:ob}
If $G'$ is a connected subgraph of $G$ of order $n'<n$ with $\delta(G')\ge 2$, then either $G'\in \cB\cup\cC\cup\cD\cup\cD_{b}\cup\cG$ or $\dtd(G')\le (n-2)/2$.
\end{ob}

We now return to the proof of Theorem~\ref{em:thm}. Suppose $G=C_{n}$ (and still $n\ge 8$). Then, by Corollary~\ref{cycle:cor}, $G \in \cC$. Hence we may assume $G$ is not a cycle. Let $\cL$ be the set of all large vertices of $G$ and let $\cS$ be the set of small vertices in $G$, i.e., $\cL=\{v\in V \mid d_{G}(v)\ge 3\}$ and $\cS = \{v\in V \mid d_{G}(v) = 2\}$. Since $G$ is not a cycle, $|\cL|\ge 1$. If $|\cL|=1$, then $G$ is a daisy, and by Corollary~\ref{cycle:cor}(b), $G\in \cD$. Hence, we may assume $|\cL|\ge 2$. Further, if $|\cL|=2$ and $G$ is a dumb-bell, then  by Corollary~\ref{cycle:cor}(c),  $G\in\cD_{b}\cup\cG_{b}$. Hence, we may assume that if $|\cL|=2$, then $G$ is not a dumb-bell.

Let $C$ be any component of $G - \cL$; it is a path. If $C$ has only one vertex, or has at least two vertices but the ends of $C$ are adjacent in $G$ to different large vertices, then we say that $C$ is a $2$-\emph{path}. Otherwise we say that $C$ is a $2$-\emph{handle}.

\begin{lem}
\label{i:lem}
If $\cL$ is not an independent set, then $G\in\{B_{10}, D_{b}(3,4), D_{b}(4,4)\} \cup \cG$.
\end{lem}
\textbf{Proof of Lemma~\ref{i:lem}.} Suppose that $\cL$ is not an independent set. Let $u$ and $v$ be two adjacent vertices in $\cL$ and let $e = uv$. By Observation \ref{bridge:ob}, $e$ is a bridge. Let $G_1=(V_{1}, E_{1})$ and $G_2=(V_{2}, E_{2})$ be the two components of $G-e$, where $u \in V_{1}$. For $i=1,2$, let $|V_i| = n_{i}$, and so $n=n_{1}+n_{2}$. Since $u,v \in \cL$ in $G$, we note that $\delta(G_1)\ge 2$ and $\delta(G_2)\ge 2$. Hence, by Observation~\ref{sg:ob}, for $i=1,2$, either $\dtd(G_{i})\le (n_{i}-2)/2$ or $G_{i}\in \cB \cup \cC \cup \cD \cup \cD_{b} \cup \cG$. If $\dtd(G_{i})\le (n_{i}-1)/2$ for $i=1,2$, then, $\dtd(G)\le (n_{1}-1)/2+(n_{2}-1)/2=(n-2)/2$, a contradiction. Hence we may assume that $\dtd(G_1)\ge n_{1}/2$. By Observation \ref{f1:ob}, if $\dtd(G_1)>n_{1}/2$, then $G_1\in \{C_{3}, C_{7}\}$ and if $\dtd(G)=n_{1}/2$, then $G_1\in\{C_{4}, C_{6}, C_{8}, C_{12}\} \cup \{B_3, D_{b}(4,4), D_{b}(3,4,1), D_{b}(3,3,2)\}$.

Let $N_v$ denote the set of neighbors of $v$ in $G_2$, and so $N_v = N_G(v) \setminus \{u\}$. If $V(G_2) = N_v \cup \{v\}$, then $\{u,v\}$ is a DTD-set of $G$, and so recalling that $n \ge 8$, $\dtd(G) = 2 < (n-2)/2$, a contradiction. Hence, $V(G_2) \ne N_v \cup \{v\}$. We proceed further with the following series of claims.

\begin{unnumbered}{Claim A}
If $G_1=C_{3}$, then $G \in \cG_{3} \cup \cG_{5}$.
\end{unnumbered}
\textbf{Proof of Claim~A.} Suppose that $G_1=C_{3}$. Since $G$ is not a dumb-bell, the graph $G_2$ is not a cycle. The vertex $v$ has at least two neighbors in $G_2$, and so $|N_v| \ge 2$. Further, every vertex in $V(G_2) \setminus N[v]$ has degree at least~$2$ in $G_2-v$.


\begin{unnumbered}{Claim A.1}
The set $N_v$ is an independent set in $G$.
\end{unnumbered}
\textbf{Proof of Claim A.1} Suppose, to the contrary, that $x_1$ and $x_2$ are two adjacent vertices in $N_v$. If $x_1$ and $x_2$ are both large vertices, let $f = x_1x_2$, while if $x_1$ is a large vertex and $x_2$ is a small vertex, let $f = vx_1$. In both cases, $G - f$ is connected and $\delta(G - f) \ge 2$, contradicting the edge-minimality of $G$. Hence, both $x_1$ and $x_2$ are small vertices. Thus, $G$ contains a $2$-handle $C$ with $|C| = 2$ and with both ends of $C$ adjacent to $v$. Since $V(G_2) \ne N_v \cup \{v\}$, we observe that $|N_v|\ge 3$ and $d_{G}(v)\ge 4$. Let $H=G-\{x_{1}, x_{2}\}$. The degree of each vertex in $H$ different from $v$ remains unchanged from its degree in $G$, and so $H$ is a connected graph with $\delta(H)\ge 2$. Further, $H$ is edge-minimal. Let $H$ have order $n_{H} = n-2$. By Observation~\ref{sg:ob}, $H\in\cB\cup\cD\cup\cD_{b}\cup\cG$ or $\dtd(H) \le (n_{H}-2)/2$.

Suppose $\dtd(H)\le (n_{H}-2)/2$. Every $\dtd(H)$-set can be extended to a DTD-set of $G$ by adding to it the vertex $v$, implying that $\dtd(G)\le \dtd(H)+1 \le (n_{H}-2)/2+1 = (n-2)/2$, a contradiction. Hence, $H\in\cB\cup\cD\cup\cD_{b}\cup\cG$. Since $G_1$ is a subgraph of $H$, we note that $H$ is not a cycle. Therefore by Observation~\ref{f1:ob}, $\dtd(H) \le n_{H}/2$.
If $v$ is a good-vertex of $H$, then a $\dtd(H)$-set that contains~$v$ is a DTD-set of $G$, implying that $\dtd(G) \le \dtd(H) \le n_{H}/2 = (n-2)/2$, a contradiction. Therefore, $v$ is a bad-vertex of $H$. By Observation~\ref{f1:ob}, either $H \in \{D(4,4), D_{b}(3,4), D_{b}(3,3,1)\}$ and $v$ is a vertex at distance~$2$ from the central vertex of $H$ or $H \in \cG_{0}$ and $v$ is the identified vertex of $H$. In both cases, $G_1$ is not a subgraph of $H$, a contradiction.~\smallqed


\begin{unnumbered}{Claim A.2}
If $N_v$ contains two small vertices in $G$, then these two vertices have only $v$ as their common neighbor.
\end{unnumbered}
\textbf{Proof of Claim A.2} Suppose, to the contrary, that $w_1$ and $w_2$ are two small vertices (of degree~$2$) in $N_v$ and that $N(w_1) = N(w_2) = \{v,x\}$. Suppose that $d_G(x) = 2$. In this case we consider the connected graph $H = G - \{w_1,w_2,x\}$. Let $H$ have order $n_{H} = n-3$. Since $G$ is not a dumb-bell, $d_H(v) \ge 2$ and $H$ is edge-minimal. Further since $H$ is not a cycle, $\dtd(H) \le n_H/2$. Since $G_1$ is a subgraph of $H$, we can choose a $\dtd(H)$-set to contain $u$ and $v$ or to contain $u$ and at least two vertices in $N_v$. In both cases, such a $\dtd(H)$-set can be extended to a DTD-set of $G$ by adding to it the vertex $w_1$, implying that $\dtd(G) \le \dtd(H) + 1$. If $\dtd(H) \le (n_H-1)/2$, then $\dtd(G) \le (n-2)/2$, a contradiction. Hence, $\dtd(H) = n_H/2$, implying by Observation~\ref{f1:ob} and our earlier observations, that $H = D_{b}(3,4,1)$ or $H = D_{b}(3,3,2)\}$. In both cases, $\dtd(G) \le 4 = (n-2)/2$, a contradiction. Hence, $d_G(x) \ge 3$.

We now consider the graph $H = G-w_1$. The degree of each vertex in $H$ different from $v$ and $x$ remains unchanged from its degree in $G$. Further, $d_H(v) \ge 2$ and $d_H(x) \ge 2$, and so $H$ is a connected graph with $\delta(H)\ge 2$, implying that $H$ is edge-minimal. Let $H$ have order $n_{H} = n-1$. Since $G_1$ is a subgraph of $H$, we can choose a $\dtd(H)$-set to contain $u$ and $v$ or to contain $u$ and at least two vertices in $N_v$. Such a $\dtd(H)$-set is a DTD-set of $G$, implying that $\dtd(G) \le \dtd(H)$. Since $H$ is not a cycle, $\dtd(H) \le n_H/2$. If $\dtd(H) \le (n_H-1)/2$, then $\dtd(G) \le (n-2)/2$, a contradiction. Hence, $\dtd(H) = n_H/2$, implying by Observation~\ref{f1:ob} and our earlier observations, that $H = D_{b}(3,3,2)$. But then $\{u,w_1,w_2\}$ is a DTD-set in $G$, and so $\dtd(G) \le 3 = (n-3)/2$, a contradiction.~\smallqed

\medskip
Among all vertices in $N_v$, let $w$ be one of minimum degree. We note that $d_{G_2 - v}(w) \ge 1$. Let $G_2^*$ be the graph obtained from $G_2 - v$ by adding as few edges as possible joining $w$ to vertices in $N_v$ so that the resulting graph is connected and has minimum degree at least~$2$. We note that by Claim~A.1 and Claim~A.2, and by the fact that $G$ is edge-minimal, the graph $G_2^*$ is edge-minimal. By the inductive hypothesis, $G_2^* \in \cB \cup \cC \cup \cD \cup \cD_{b} \cup \cG$ or $\dtd(G_2^*)\le (n_{2}^{*}-2)/2$. Let $S^*$ be a $\dtd(G_2^*)$-set.

\begin{unnumbered}{Claim~A.3}
$G_2^* \in \cB \cup \cC \cup \cD \cup \cD_{b} \cup \cG$.
\end{unnumbered}
\textbf{Proof of Claim~A.3} Suppose $\dtd(G_2^*) \le (n_{2}^*-2)/2 = (n-6)/2$. If $w \in S^*$ or if $S^* \cap N_v  = \emptyset$, let $S = S^* \cup \{u,v\}$. If $w \notin S^*$ and $|S^* \cap N_v| \ge 1$, let $S = S^* \cup \{u,w\}$. We show that $S$ is a DTD-set in $G$. Suppose that there is a vertex $x$ in $G$ that is not DT-dominated by the set $S$. Then, $x$ has no neighbor in $S$ and is at distance~$2$ from at most one vertex of $S$ in $G$. If $x = u$, then $v \notin S$, implying that $S = S^* \cup \{u,w\}$. In this case, $|S \cap N_v| \ge 2$ and therefore $x$ is at distance~$2$ from at least two vertices of $S$, a contradiction. Hence, $x \ne u$. Since $u \in S$, we have that $x \notin N(u)$, and so $x \in V(G_2^*)$. Suppose $x \in N_v$. Then, $v \notin S$, implying that $w \in S$, $w \notin S^*$ and $|S^* \cap N_v| \ge 1$. Let $w^* \in S^* \cap N_v$. If $x \ne w$, then $x$ is at distance~$2$ from both $u$ and $w$. If $x = w$, then $x$ is at distance~$2$ from both $u$ and $w^*$. In both cases, $x$ is at distance~$2$ from at least two vertices of $S$, a contradiction. Hence, $x \notin N_v$. Thus, the neighbors of $x$ in $G$ and $G_2^*$ are the same. Since $S^* \subset S$ and $x$ has no neighbor in $S$, the vertex $x$ has no neighbor in $S^*$. However, $S^*$ is a DTD-set of $G_2^*$, implying that $x$ is at distance~$2$ from at least two vertices of $S^*$ in $G_2^*$. Since $x$ is at distance~$2$ from at most one vertex of $S$ in $G$ and since $S^* \subset S$, there is a vertex $x^* \in S^*$ at distance~$2$ from $x$ in $G_2^*$ but at distance greater than~$2$ from $x$ in $G$.  This is only possible if $w$ is the only common neighbor of  $x$ and $x^*$ in $G_2^*$ and if $wx^*$ was an edge added to $G_2 - v$ when forming $G_2^*$. Therefore, $x^* \in S^* \cap N_v$ and $x^* \ne w$. In particular, we note that $w \notin S^*$ and $|S^* \cap N_v| \ge 1$. But then $w \in S$, implying that $x$ has a neighbor in $S$, a contradiction. Therefore, $S$ is a DTD-set in $G$. Thus, $\dtd(G) \le |S| = |S^*| + 2 = \dtd(G_2^*) + 2 \le (n-6)/2 + 2 = (n-2)/2$, a contradiction. Hence, $G_2^* \in\cB \cup \cC \cup \cD \cup \cD_{b} \cup \cG$.~\smallqed

\begin{unnumbered}{Claim~A.4}
At least one edge was added to $G_2 - v$ when forming $G_2^*$.
\end{unnumbered}
\textbf{Proof of Claim~A.4}
On the one hand, suppose that $d_G(w) = 2$. Then, $d_{G_2 - v}(w) = 1$ and therefore at least one new edge incident with $w$ was added to $G_2 - v$ to guarantee that $\delta(G_2^*) \ge 2$. On the other hand, suppose that $d_G(w) \ge 3$. By our choice of $w$, every neighbor of $v$ in $G_2$ has degree at least~$3$ in $G$ in this case. If $G_2 - v$ is connected, then removing from $G$ an arbitrary edge joining $v$ with a vertex in $N_v$ produces a connected graph with minimum degree~$2$, contradicting the edge-minimality of $G$. Hence, $G_2 - v$ is disconnected. Therefore at least one new edge incident with $w$ was added to $G_2 - v$ to guarantee that $G_2^*$ is connected.~\smallqed

\medskip
By Claim~A.4, at least one edge was added to $G_2 - v$ when forming $G_2^*$. Let $f = wx$ be such an added edge. We note that $\{w,x\} \subseteq N_v$.

\begin{unnumbered}{Claim~A.5}
$G_2^* \in \cB \cup \cD \cup \cD_{b} \cup \cG$.
\end{unnumbered}
\textbf{Proof of Claim~A.5} By Claim~A.3, $G_2^* \in \cB \cup \cC \cup \cD \cup \cD_{b} \cup \cG$. Suppose that $G_2^* \in \cC$. Then the edge $f$ was the only edge added to $G_2 - v$ when forming $G_2^*$, and so $G_2 - v$ is a path whose ends, namely $w$ and $x$, are both adjacent to $v$ in $G_2$. By the edge-minimality of $G$, the vertex $v$ is adjacent to no other vertex on this path, implying that $G_2$ is a cycle and therefore $G$ is a dumb-bell, a contradiction.~\smallqed

\begin{unnumbered}{Claim A.4}
The edge $f$ is a bad-edge in $G_2^*$.
\end{unnumbered}
\textbf{Proof of Claim~A.6}  Suppose to the contrary that $f$ is a good-edge in $G_2^*$. Then the $\dtd(G_2^*)$-set, $S^*$, can be chosen to contain both $w$ and $x$. With this choice of $S^*$, let $S = S^* \cup \{u\}$. Then, $S$ is a DTD-set of $G$, and so $\dtd(G) \le |S| + 1 = \dtd(G_2^*)+1$. By Observation~\ref{f1:ob} and Claim~A.5, we have that $\dtd(G_2^*) \le n_2^*/2 = (n-4)/2$, implying that $\dtd(G) \le (n-2)/2$, a contradiction.~\smallqed

\begin{unnumbered}{Claim~A.7}
$G_2^* \in \cG \setminus \cG_1$.
\end{unnumbered}
\textbf{Proof of Claim~A.7} By Claim~A.6, the edge $f$ is a bad-edge in $G_2^*$. By Observation~\ref{f3a:ob} and Claim A.3, we have that $G_2^* \in  \{D(4,4), D_{b}(3,4), D_{b}(3,3,1)\} \cup \cG \setminus \cG_1$.
Suppose to the contrary that $G_2^* \in  \{D(4,4), D_{b}(3,4), D_{b}(3,3,1)\}$. Then, by Observation~\ref{f3a:ob} and by the edge-minimality of $G$ and the construction of $G_2^*$, the edge $f$ was the only edge added to $G_2 - v$ when forming $G_2^*$. The graph $G$ is therefore determined and has order~$n = 11$, and it can be readily checked that $\dtd(G) = 4 < (n-1)/2$, a contradiction.~\smallqed

\medskip
By Claim~A.6, the edge $f$ is a bad-edge in $G_2^*$. By Claim~A.7, we have that $G_2^* \in \cG \setminus \cG_1$. Let $v^*$ be the identified vertex of $G_2^*$. By Observation~\ref{f3a:ob}(d), the edge $f$ is incident with the vertex $v^*$.

\begin{unnumbered}{Claim~A.8}
$G_2^* \in \cG_{0}$.
\end{unnumbered}
\textbf{Proof of Claim~A.8} Suppose to the contrary that $G_2^* \notin \cG_{0}$. Then, $G_2^* \in \cG_{i}$ for some $i$, $2 \le i \le 6$. Let $x^*$ be the neighbor of $v^*$ that  belongs to the $X_i$-unit in $G_2^*$. By Observation~\ref{f3b:ob}, the vertex $v^*$ is a good-vertex of $G_2^*$. Thus the $\dtd(G_2^*)$-set, $S^*$, can be chosen to contain $v^*$. Further, in each unit of type-1 or type-2 in $G_2^*$, we can choose $S^*$ to contain the neighbor of $v^*$ in that unit as well as a vertex at distance~$2$ from $v^*$ is that unit. If $v$ is not adjacent to $x^*$ in $G$, then the set $S^* \cup \{u\}$ is a DTD-set of $G$, implying that $\dtd(G) \le |S^*| + 1 = \dtd(G_2^*) + 1 = (n_2^* - 1)/2 + 1 = (n - 3)/2$, a contradiction. Hence, $vx^*$ is an edge of $G$. In this case, the set $(S^* \cup \{u,v\}) \setminus \{v^*\}$ is a DTD-set of $G$, implying that $\dtd(G) \le |S^*| + 1 = (n - 3)/2$, a contradiction.~\smallqed

\medskip
By Claim~A.8, $G_2^* \in \cG_{0}$. The $\dtd(G_2^*)$-set, $S^*$, can be chosen to contain the neighbor of $v^*$ in each unit of type-1 or type-2 in $G_2^*$ as well as a vertex at distance~$2$ from $v^*$ in each unit. By the edge-minimality of $G$, and the way in which $G_2^*$ is constructed, if $z$ is a neighbor of $v$ in $G_2$, then either $z = v^*$ or $z$ is a neighbor of $v^*$ in $G_2^*$. Thus if $d_{G_2}(v) \ge 3$, then the set $S^* \cup \{u\}$ is a DTD-set of $G$, implying that $\dtd(G) \le |S^*| + 1 = \dtd(G_2^*) + 1 = (n_2^* - 1)/2 + 1 = (n - 3)/2$, a contradiction. Hence, $d_{G_2}(v) = 2$, implying that the ends of the edge $f$, namely $w$ and $x$, are the only neighbors of $v$ in $G_2$. As observed earlier, the edge $f$ is incident with the vertex $v^*$. Let $\{u^*,v^*\} = \{w,x\}$. If $u^*$ belongs to a type-1 unit in $G_2^*$, then $G \in \cG_3$. If $u^*$ belongs to a type-2 unit in $G_2^*$, then $G \in \cG_5$. This completes the proof of Claim~A.~\smallqed

\medskip
By Claim~A, we may assume that $G_1 \ne C_3$, for otherwise the desired result follows. Analogously, we may assume that $G_2 \ne C_3$.

\begin{unnumbered}{Claim B}
$G_1 \ne C_7$.
\end{unnumbered}
\textbf{Proof of Claim~B.} Suppose to the contrary that $G_1 = C_7$. Then,  $n = n_2 + 7$. Let $G_1$ be the cycle $uu_{1}u_{2}\ldots u_{6}u$. By Observation~\ref{sg:ob}, $\dtd(G_2)\le (n_2-2)/2$ or $G_2 \in \cB \cup \cC \cup \cD \cup \cD_{b} \cup \cG$. Let $S$ be a $\dtd(G_2)$-set. Suppose that $\dtd(G_2) \le (n_2-1)/2 = n/2 - 4$. Then the set $S \cup \{u_{3}, u_{4},v\}$ is a DTD-set of $G$, implying that $\dtd(G) \le \dtd(G_2) + 3 = n/2 - 1$, a contradiction. Hence, $\dtd(G_2) \ge n_2/2$ and $G_2 \in \cB \cup \cC \cup \cD \cup \cD_{b} \cup \cG$. Since $G$ is not a dumb-bell, the graph $G_2$ is not a cycle. Thus by Observation~\ref{f1:ob}, $\dtd(G_2) = n_2/2$ and $G_2 \in \{B_3, D_{b}(4,4), D_{b}(3,4,1), D_{b}(3,3,2)\}$. By Observation~\ref{f3b:ob}, every vertex of $G_2$ is a good-vertex. In particular, the vertex $v$ is a good-vertex of $G_2$. Choosing the $\dtd(G_2)$-set, $S$, to contain the vertex $v$, we have that the set $S \cup \{u_{3}, u_{4}\}$ is a DTD-set of $G$, implying that $\dtd(G) \le |S| + 2 = \dtd(G_2) + 2 = n_2/2 + 2 = (n-3)/2$, a contradiction.~\smallqed

\medskip
By Claim~B, we have $G_1 \ne C_7$. Analogously, we have $G_2 \ne C_7$.  Since $G_1 \notin \{C_3,C_7\}$, we have that $\dtd(G) = n_{1}/2$ and $G_1 \in \{C_{4}, C_{6}, C_{8}, C_{12}\} \cup \{B_3, D_{b}(4,4), D_{b}(3,4,1), D_{b}(3,3,2)\}$. If $\dtd(G_2) \le (n_2-2)/2$, then $\dtd(G) \le \dtd(G_1) + \dtd(G_2) \le n_1/2 + (n_2-2)/2 = (n-2)/2$, a contradiction. Hence, $\dtd(G_2) \ge (n_{2}-1)/2$. By Observation~\ref{sg:ob}, $G_2 \in \cB \cup \cC \cup \cD \cup \cD_{b} \cup \cG$. By our earlier assumptions and observations, $G_2 \notin \{C_3,C_7\}$, implying that $\dtd(G_2) \le n_{2}/2$.

\begin{unnumbered}{Claim C} If $G_1 = C_4$, then $G = B_{10}$ or $G \in \cG$.
\end{unnumbered}
\textbf{Proof of Claim~C.} Suppose that $G_1 = C_4$. Then, $n=n_{2}+4$. As observed earlier, $G_2 \notin \{C_3,C_7\}$ and $G_2 \in \cB \cup \cC \cup \cD \cup \cD_{b} \cup \cG$. Since $G$ is not a dumb-bell, the graph $G_2$ is not a cycle. Thus, $G_2 \notin \cC$, and so $G_2 \in \cB \cup \cD \cup \cD_{b} \cup \cG$. By the edge-minimality of $G$, the vertex $v$ is a special vertex in $G_2$. Recall that the graph $G_2^v$ is obtained from $G_2$ by adding a new vertex $v'$ and adding the pendant edge $vv'$.

Suppose $\ndtd(G_2^v) \le \dtd(G_2)$. Let $S_2$ be a $\ndtd(G_2^v)$-set. Let $u_1$ be a neighbor of $u$ in $G_1$.  Then the set $(S_2 \setminus \{v'\}) \cup \{u,u_1\}$ is a DTD-set of $G$, and so $\dtd(G) \le |S_2| + 1 = \ndtd(G_2^v) + 1 \le \dtd(G_2) + 1 \le n_2/2 + 1 = (n-2)/2$, a contradiction. Hence, $\ndtd(G_2^v) > \dtd(G_2)$. Applying Observation~\ref{special:ob} to the graph $G_2$ and the special vertex $v$ of $G_2$, one of the three conditions (b), (c) or (d) in the statement of the observation hold. We consider each of the three conditions in turn.

If $G_2 = B_1$ and $v$ has degree~$3$ in $G_2$, then $G = B_{10}$.
If $G_2 \in \cG$ and $v$ is the identified vertex of $G_2$, then $G \in \cG$ and $G$ has one additional type-$2$ unit than does $G_2$.
Suppose $G_2 \in \cG_0$ and $v$ is a neighbor of the identified vertex of $G_2$. Then, $G_2 = G_0(i,j)$ for some $i,j$ where $i + j \ge 2$ and either $v$ belongs to a type-$1$ unit or a type-$2$ unit of $G_2$. If $v$ belongs to a type-$1$ unit of $G_2$, then $G \in G_5(i-1,j) \in \cG_5$. If $v$ belongs to a type-$2$ unit of $G_2$, then $G \in G_6(i,j-1) \in \cG_6$. Hence if $G_2 \in \cG_0$, then $G \in \cG$.~\smallqed

\medskip
By Claim~C, we may assume that $G_1 \ne C_4$, for otherwise the desired result follows. Analogously, we may assume that $G_2 \ne C_4$.
By the edge-minimality of $G$, the vertex $u$ is a special vertex in $G_1$ and the vertex $v$ is a special vertex in $G_2$. Recall that the graph $G_1^u$ is obtained from $G_1$ by adding a new vertex $u'$ and adding the pendant edge $uu'$.
By our earlier observations, $G_1 \in \{C_{6}, C_{8}, C_{12}\} \cup \{B_3, D_{b}(4,4), D_{b}(3,4,1), D_{b}(3,3,2)\}$. Further, $G_2 \notin \{C_3,C_4,C_7\}$ and $G_2 \in \cB \cup \cC \cup \cD \cup \cD_{b} \cup \cG$.

\begin{unnumbered}{Claim D}
$G_2 \in \cG_0$ and that $v$ is the identified vertex of $G_2$.
\end{unnumbered}
\textbf{Proof of Claim~D.} Suppose that $v$ is a good-vertex of $G_2$. Let $S_2$ be a $\dtd(G_2)$-set that contains the vertex~$v$. By Observation~\ref{special:ob}, we have that $\ndtd(G_1^u) \le \dtd(G_1)$. Let $S_1$ be a $\ndtd(G_1^u)$-set and note that $u' \in S_1$. Then the set $(S_1 \setminus \{u'\}) \cup S_2$ is a DTD-set of $G$, and so $\dtd(G) \le |S_1| + |S_2| - 1 = \ndtd(G_1^u) + \dtd(G_2) - 1 \le \dtd(G_1) + \dtd(G_2) - 1 \le n_1/2 + n_2/2 - 1 = (n-2)/2$, a contradiction. Hence, $v$ is a bad-vertex of $G_2$.
Applying Observation~\ref{f3b:ob} to the graph $G_2$, we have that $G_2 \in \{D(4,4), D_{b}(3,4), D_{b}(3,3,1)\}$ or $G_2 \in \cG_{0}$. Further, if $G_2 \in \cG_0$, then $v$ is the identified vertex of $G_2$.

By Observation~\ref{f3b:ob}, every vertex of $G_1$ is a good-vertex of $G_1$. Let $D_1$ be a  $\dtd(G_1)$-set that contains the vertex~$u$. Suppose that $\ndtd(G_2^v) \le \dtd(G_2)$. Let $D_2$ be a $\ndtd(G_2^v)$-set. Then the set $D_1 \cup (D_2 \setminus \{v'\})$ is a DTD-set of $G$, and so $\dtd(G) \le |D_1| + |D_2| - 1 \le \dtd(G_1) + \ndtd(G_2^v) - 1 \le \dtd(G_1) + \dtd(G_2) -1 \le n_1/2 + n_2/2 - 1 = (n-2)/2$, a contradiction. Hence, $\ndtd(G_2^v) > \dtd(G_2)$. By Observation~\ref{special:ob}, we therefore have that $G_2 \notin \{D(4,4), D_{b}(3,4), D_{b}(3,3,1)\}$, implying that $G_2 \in \cG_0$ and that $v$ is the identified vertex of $G_2$.~\smallqed

\medskip
We now return to the proof of Lemma~\ref{i:lem} one last time. By Claim~D, $G_2 \in \cG_0$ and $v$ is the identified vertex of $G_2$. Let $S_1$ be a $\ndtd(G_1^u)$-set and note that $u' \in S_1$. Let $S_2$ be a $\dtd(G_2)$-set that contains all neighbors of $v$ in $G_2$ and a vertex at distance~$2$ from $v$ in each unit of $G_2$. Then, $\dtd(G_2) = |S_2| = (n_2 - 1)/2$.
As observed earlier, we have that $G_1 \in \{C_{6}, C_{8}, C_{12}, B_3, D_{b}(4,4)$, $D_{b}(3,4,1), D_{b}(3,3,2)\}$. We consider the seven possibilities in turn.

If $G = C_6$, then $G \in \cG_2$.

If $G_1 = C_8$, let $G_1$ be given by $uu_1u_2 \ldots u_7u$. Then the set $S_2 \cup \{u,u_3,u_4\}$ is a DTD-set of $G$, and so $\dtd(G) \le |S_2| + 3 = \dtd(G_1) + \dtd(G_2) - 1 = n_1/2 + (n_2-1)/2 - 1 = (n-3)/2$, a contradiction. Hence, $G_1 \ne C_8$.

Suppose $G_1 \in \{C_{12},B_3\}$. Then, $\ndtd(G_1^u) = \dtd(G_1) - 1$. The set $(S_1 \setminus \{u'\}) \cup (S_2 \cup \{v\})$ is a DTD-set of $G$, and so $\dtd(G) \le |S_1| + |S_2| = \dtd(G_1) + \dtd(G_2) - 1 = n_1/2 + (n_2-1)/2 - 1 = (n-3)/2$, a contradiction. Hence, $G_1 \notin \{C_{12},B_3\}$.

Suppose $G_1 = D_{b}(4,4)$. If $u$ has degree~$2$ in $G_1$, then $\ndtd(G_1^u) = \dtd(G_1) - 1$, and as before we obtain a contradiction. Hence, $d_{G_1}(u) = 3$. But then $G \in \cG_6$.

Suppose $G_1 = D_{b}(3,4,1)$. Then, $G_1$ can be obtained from a path $u_1u_2 \ldots u_8$ by adding the edges $u_1u_3$ and $u_5u_8$. Since $u$ is a special vertex of $G_1$, we note that $u \in \{u_3,u_4,u_5,u_7\}$. If $u = u_3$, let $S = S_2 \cup \{u_3,u_5,u_6\}$. If $u = u_7$, let $S = S_2 \cup \{u_3,u_4,u_5\}$. In both cases, $S$ is a DTD-set of $G$, and so $\dtd(G) \le |S| = |S_2| + 3 = \dtd(G_1) + \dtd(G_2) - 1 = (n-3)/2$, a contradiction. Hence, $u \in \{u_4,u_5\}$. If $u = u_4$, then $G \in \cG_5$. If $u = u_5$, then $G \in \cG_4$.

Suppose $G_1 = D_{b}(3,3,2)$. Then, $G_1$ can be obtained from a path $u_1u_2 \ldots u_8$ by adding the edges $u_1u_3$ and $u_6u_8$. Since $u$ is a special vertex of $G_1$, we note that $u \in \{u_3,u_4,u_5,u_6\}$. By symmetry, we may assume that $u \in \{u_3,u_4\}$. If $u = u_3$, then $S_2 \cup \{u_3,u_5,u_6\}$ is a DTD-set of $G$, and so $\dtd(G) \le |S_2| + 3 = \dtd(G_1) + \dtd(G_2) - 1 = (n-3)/2$, a contradiction. Hence, $u = u_4$. But then $G \in \cG_3$.

Hence we have shown that $G_1 \in \{C_{6}, D_{b}(4,4)$, $D_{b}(3,4,1), D_{b}(3,3,2)\}$ and that $G \in \cG$. This completes the proof of Lemma~\ref{i:lem}.~\qed

\medskip
By Lemma~\ref{i:lem}, we may assume that $\cL$ is not an independent set, for otherwise the desired result follows.

\begin{lem}
\label{path6:lem}
If $G$ contains a path on six vertices each internal vertex of which has degree~$2$  in $G$ and whose end vertices are not adjacent, then $G\in \{B_3, B_{4}, B_{5}, B_{7}, B_{11}\}$.
\end{lem}
\textbf{Proof of Lemma~\ref{path6:lem}.} Let $u$ and $v$ be non-adjacent vertices in $G$ joined by a path $uw_{1}w_{2}w_{3}w_{4}v$ every internal vertex of which has degree~$2$ in $G$. Let $G'$ be the graph obtained from $G$ by removing the vertices $w_{1}, w_{2}, w_{3}$ and $w_{4}$, and adding the edge $uv$.
Then, $G'$ is a connected graph of order $n'=n-4$ with $\delta(G')\ge 2$.
By Observation~\ref{sd:ob}, $\dtd(G)\le \dtd(G')+2$. It follows that if $\dtd(G')\le (n'-2)/2$, then $\dtd(G)\le (n-6)/2+2=(n-2)/2$, a contradiction. Hence, $\dtd(G')\ge (n'-1)/2$. Let $F=G'-e$. Then, $F$ has order~$n'$ and either $F$ is disconnected or $\delta(F) = 1$ or $F$ is connected and $\delta(F) \ge 2$.

\begin{unnumbered}{Claim~E} If $F$ is disconnected or $\delta(F) = 1$, then $G \in \{B_{5}, B_{11}\} \subset \cB$.
\end{unnumbered}
\textbf{Proof of Claim~E.} Assume that $F$ is disconnected or $\delta(F) = 1$. Then, $G'$ is  an edge-minimal graph. As observed earlier, $\dtd(G')\ge (n'-1)/2$. Thus, $G'$ is a $\frac{1}{2}$-minimal graph. Applying the induction hypothesis to the graph $G'$, we see that $G' \in \cB \cup \cC \cup \cD \cup \cD_{b} \cup \cG$. If $G' \in \cC \cup \cD \cup \cD_{b} \cup \cG_{b}$, then $G \in \cC \cup \cD \cup \cD_{b} \cup \cG_{b}$, contradicting our previous assumptions. Hence, $G' \in \cB \cup (\cG \setminus \cG_{b})$.

\begin{unnumbered}{Claim~E.1} $G' \in \cB$.
\end{unnumbered}
\textbf{Proof of Claim~E.} Suppose to the contrary that $G' \notin \cB$. Then, $G' \in \cG \setminus \cG_{b}$. Let $x$ be the identified vertex of $G$.
Suppose that $G' \in \cG_0$, and so $G' = G_0(n_1,n_2)$ for some $n_1 \ge 0$ and $n_2 \ge 0$ where $n_1 + n_2 \ge 2$. Since $G' \notin \cG_{b}$, we note that $n_1 + n_2 \ge 3$. By Observation~\ref{f4:ob}, there is a $\dtd(G')$-set, $S'$, such that $v \notin S'$, $N(v) \subset S'$ and $S'$ totally dominates $N(v)$ in $G'$. In particular, we note that in each unit in $G'$ there is exactly one vertex at distance~$2$ from $x$ in $G'$ that belongs to $S'$. If $e$ is a bridge of $G'$, then we may assume, renaming vertices if necessary, that $x$ is incident with $e$ and that $x = u$. In this case, we let $S = S' \cup \{w_1\}$. If $e$ belongs to a $3$-cycle in $G'$, then renaming vertices, if necessary, we may assume that this $3$-cycle is given by $zuvz$ and that $xyz$ is a path in $G'$. In this case, we note that $\{y,z\} \subset S'$ and we let $S = (S \setminus \{z\}) \cup \{w_2,w_3\}$. If $e$ belongs to a $4$-cycle in $G'$, then renaming vertices, if necessary, we may assume that this $4$-cycle is given by $yuvzy$ and that $xy$ is a path in $G'$. In this case, we may assume that $\{y,z\} \subset S'$ (since if $u \in S'$, we simply replace $u$ in $S'$ by $z$) and let $S = (S \setminus \{z\}) \cup \{w_2,w_3\}$. In all three cases, the set $S$ is a DTD-set of $G$, and so $\dtd(G) \le |S| + 1 = (n'-1)/2 + 1 = (n-3)/2$, a contradiction. Hence, $G' \notin \cG_0$.

Suppose that $G' \in \cG_1$, and so $G' = G_1(n_1,n_2)$ for some $n_1 \ge 0$ and $n_2 \ge 0$ where $n_1 + n_2 \ge 1$. Since $G' \notin \cG_{b}$, we note that $n_1 + n_2 \ge 2$. By Observation~\ref{f4:ob}, there is a $\dtd(G')$-set, $S'$, such that $v \notin S'$ and $N(v) \subset S'$. Further we can choose such a set $S'$ so that each unit in $G'$ contains a vertex at distance~$2$ from $x$ in $G'$ that belongs to $S'$.
If $n_2 \ge 1$, then since $\cL$ is an independent set, we have that $n_2 = 1$ and that $e$ is the edge that joins $x$ to the vertex of the $4$-cycle in $G'$. Renaming $u$ and $v$, if necessary, we may assume that $u = x$. In this case, we let $S = S' \cup \{w_1\}$. Then, $S$ is a DTD-set of $G$, and so $\dtd(G) \le |S| + 1 = (n'-1)/2 + 1 = (n-3)/2$, a contradiction. Hence, $G' = G_1(n_1,0)$ for some $n_1 \ge 2$.
If $e$ belongs to the $5$-cycle in $G'$, then we may assume, renaming vertices if necessary, that $x$ is incident with $e$ and that $x = u$. In this case, we note that $v \in S'$ and we let $S = (S' \setminus \{v\}) \cup \{w_3,w_4\}$.
If $e$ is a bridge of $G$, then we may assume, renaming vertices if necessary, that $x$ is incident with $e$ and that $x = u$. In this case, we let $S = S' \cup \{w_1\}$.
If $e$ belongs to a $3$-cycle in $G'$, then renaming vertices, if necessary, we may assume that this $3$-cycle is given by $zuvz$ and that $xyz$ is a path in $G'$. In this case, we note that $\{y,z\} \subset S'$ and we let $S = (S \setminus \{z\}) \cup \{w_2,w_3\}$.
In all three cases, the set $S$ is a DTD-set of $G$, and so $\dtd(G) \le |S| + 1 = (n'-1)/2 + 1 = (n-3)/2$, a contradiction. Hence, $G' \notin \cG_1$.

Suppose that $G' \in \cG_2$, and so $G' = G_2(n_1,n_2)$ for some $n_1 \ge 0$ and $n_2 \ge 0$ where $n_1 + n_2 \ge 1$. Since $G' \notin \cG_{b}$, we note that $n_1 + n_2 \ge 2$. By Observation~\ref{f4:ob}, there is a $\dtd(G')$-set, $S'$, such that $v \notin S'$, $N(v) \subset S'$ and $S'$ totally dominates $N(v)$ in $G'$. Since $\cL$ is an independent set, $e$ is the edge that joins $x$ to the vertex of the $6$-cycle in $G'$. Renaming $u$ and $v$, if necessary, we may assume that $u = x$. Let $C \colon vv_1v_2v_3v_4v_5v$ be the $6$-cycle in $G'$. We may choose $S'$ so that $S \cap V(C) = \{v,v_1,v_5\}$. With this choice of the set $S'$, the set $S' \cup \{w_4\}$ is a DTD-set of $G$, and so $\dtd(G) \le |S'| + 1 = (n'-1)/2 + 1 = (n-3)/2$, a contradiction. Hence, $G' \notin \cG_2$.

Suppose that $G' \in \cG_3$. Since $\cL$ is an independent set in $G$, either $G' = G_3(1,0)$ or $G' = G_3(0,1)$. In both cases, $e$ is the edge joining the two vertices of degree~$3$ in $G'$. The graph $G$ is therefore determined and has order~$n = 17$ and $\dtd(G) \le 7 = (n-3)/2$, a contradiction. Hence, $G' \notin \cG_3$.

Suppose that $G' \in \cG_4$. Since $\cL$ is an independent set in $G$, either $G' = G_4(1,0)$ or $G' = G_4(0,1)$. In both cases, $G'$ has order $n' = 13$ and, up to isomorphism, there are five different choice for the edge $e$. This gives rise to a total of ten possible (non-isomorphic) constructions for the graph $G$. However in all ten cases, we have $n = 17$ and $\dtd(G) \le 7 = (n-3)/2$, a contradiction. Hence, $G' \notin \cG_4$.

If $G' \in \cG_5$, then for every possible choice of the edge $e$, we will always produce two adjacent large vertices in $G$, contradicting our assumption that $\cL$ is an independent set. Hence, $G' \notin \cG_5$.

Suppose that $G' \in \cG_6$. Since $\cL$ is an independent set in $G$, either $G' = G_6(1,0)$ or $G' = G_6(0,1)$. In both cases, $e$ is the edge joining the two large vertices in $G'$. The graph $G$ is therefore determined and has order~$n = 17$ and $\dtd(G) \le 7 = (n-3)/2$, a contradiction. Hence, $G' \notin \cG_6$.~\smallqed

\begin{unnumbered}{Claim~E.2} $G' \in \{B_1,B_5\}$.
\end{unnumbered}
\textbf{Proof of Claim~E.} By Claim~E.1, $G' \in \cB$. We wish to show that $G' \in \{B_1,B_5\}$. Suppose to the contrary that $G' \notin \{B_1,B_5\}$.

Suppose that $G'=B_2$. Up to isomorphism, there are two different choices for the edge~$e$. This gives rise to two possible (non-isomorphic) constructions for the graph $G$ from the graph $G'$. However in both cases, we have $n = 11$ and $\dtd(G) \le 4 = (n-3)/2$, a contradiction.

Suppose that $G'=B_3$. Up to isomorphism, there are two different choices for the edge~$e$. This gives rise to two possible (non-isomorphic) constructions for the graph $G$ from the graph $G'$. However in both cases, we have $n = 12$ and $\dtd(G) \le 5 = (n-2)/2$, a contradiction.

Suppose that $G' \in \{B_4,B_6,B_7,B_8,B_9,B_{10} \}$. Then, $n' = 9$, and so $n = 13$. If $G' \in \{B_4,B_6,B_8\}$, then up to isomorphism, there are two different choices for the edge~$e$. If $G' \in \{B_7,B_9,B_{10}\}$, then up to isomorphism, there are three different choices for the edge~$e$. In all cases, the resulting graph $G$ constructed from $G'$ satisfies $\dtd(G) \le 5 = (n-3)/2$, a contradiction.

Suppose that $G'=B_{11}$. Up to isomorphism, there are two different choices for the edge~$e$. This gives rise to two possible (non-isomorphic) constructions for the graph $G$ from the graph $G'$. However in both cases, we have $n = 17$ and $\dtd(G) \le 7 = (n-3)/2$, a contradiction.
Since all the above cases produce a contradiction, the desired result of the claim follows.~\smallqed

\medskip
We now return to the proof of Claim~E. By Claim~E.2, $G' \in \{B_1,B_5\}$. If $G' = B_1$, then $G = B_5$. Suppose $G' = B_5$. Then, $n' = 9$ and $n = 13$. Up to isomorphism, there are two different choices for the edge~$e$. If $e$ is on one of the two paths of length~$2$ that join the two vertices of degree~$3$ in $G'$, then $\dtd(G) \le 5 = (n-3)/2$, a contradiction. Hence, $e$ is on the path of length~$6$ that join the two vertices of degree~$3$ in $G'$. In this case, $G = B_{11}$. This completes the proof of Claim~E.~\qed

\medskip
We next consider the case when $F$ is connected and $\delta(F) \ge 2$.

\begin{unnumbered}{Claim~F} If $F$ is connected and $\delta(F) \ge 2$, then $G \in \{B_3,B_4,B_7\}$.
\end{unnumbered}
\textbf{Proof of Claim~F.} Assume that $F$ is connected and $\delta(F) \ge 2$. By the edge-minimality of $G$ and since $\dtd(F) \ge \dtd(G') \ge (n'-1)/2$, the graph $F$ is a $\frac{1}{2}$-minimal graph. By the inductive hypothesis, $F \in \cB \cup \cC \cup \cD \cup \cD_{b} \cup \cG$.

\begin{unnumbered}{Claim F.1} Both $u$ and $v$ are special vertices in $F$.
\end{unnumbered}
\textbf{Proof of Claim~F.1} Suppose to the contrary that $u$ is not a special vertex in $F$. Then, there is an edge $f$ in $F$ which may be removed without disconnecting $F$ and with no vertex, except possibly for the vertex $u$, having degree~$1$ in $F$. But then $\delta(G - f) \ge 2$ and $G$ is connected, contradicting the edge-minimality of $G$. Hence, $u$ is a special vertex in $F$. Analogously, $v$ is a special vertex in $F$.~\smallqed

\medskip
Recall from Section~\ref{NDTDS}, that for a graph $H$ and vertex $x$ in $H$, the graph obtained from $H$ by adding a new vertex $x'$ and adding the pendant edge $xx'$ is denoted by $H^x$.

\begin{unnumbered}{Claim F.2} Let $x \in \{u,v\}$. If $\ndtd(F^x) \le \dtd(F)$, then there is no $\ndtd(F^x)$-set which DT-dominates $x'$.
\end{unnumbered}
\textbf{Proof of Claim~F.2}  Let $x \in \{u,v\}$ and suppose, to the contrary, that $\ndtd(F^x) \le \dtd(F) = \dtd(G' -e)$ but there exists a $\ndtd(F^x)$-set, $S$ say, which DT-dominates $x'$. Renaming $u$ and $v$ if necessary, we may assume that $x = v$. Then the set $(S \setminus \{v'\}) \cup \{u,w_4\}$ is a DTD-set of $G$, and so $\dtd(G) \le |S| + 1 = \dtd(F) + 1$. If $\dtd(F) \le n'/2$, then $\dtd(G) \le n'/2 + 1 = (n-4)/2 + 2 = (n-2)/2$, a contradiction. Hence, $\dtd(F)=(n'+1)/2$. By Observation~\ref{f1:ob}, $F \in \{C_3,C_7\}$. Since $u$ and $v$ are not adjacent in $F$, the case $F = C_3$ cannot occur. Therefore, $F = C_7$ and $u$ and $v$ are at distance either~$2$ or~$3$ apart in $F$. This gives rise to two possible constructions for the graph $G$ from the graph $F$. However in both cases, we have $n = 11$ and $\dtd(G) \le 4 = (n-3)/2$, a contradiction.~\smallqed

\begin{unnumbered}{Claim~F.3} If $x \in \{u,v\}$, then $\ndtd(F^x) > \dtd(F)$.
\end{unnumbered}
\textbf{Proof of Claim~F.3} Let $x \in \{u,v\}$ and suppose, to the contrary, that $\ndtd(F^x) \le \dtd(F)$. Renaming $u$ and $v$ if necessary, we may assume that $x = v$. By Claim~F.2, there is no $\ndtd(F^v)$-set which DT-dominates $v'$. By Observation~\ref{ndtd1:ob}, either $F = D(4,4)$ and $v$ is at distance~2 from the vertex of degree~$4$ in $F$ or $F \in \{D_{b}(3,4), D_{b}(3,3,1)\}$, $d_F(v) = 3$, and $v$ belongs to a $3$-cycle in $F$.

Suppose that $F = D(4,4)$ and $v$ is at distance~2 from the vertex of degree~$4$ in $F$. Since $\cL$ is an independent set in $G$, and since $u$ and $v$ are not adjacent in $F$, either $u$ is the vertex of degree~$4$ in $F$ or $u$ is the vertex at distance~$2$ from $v$ in $F$.

Suppose that $F = D_{b}(3,4)$. Then, $v$ is the vertex of degree~$3$ in $F$ that belongs to a $3$-cycle. Since $\cL$ is an independent set in $G$, and since $u$ and $v$ are not adjacent in $F$, the vertex $u$ is the vertex in $F$ at distance~$3$ from $v$ in $F$.

Suppose that $F = D_{b}(3,3,1)$. Then, $v$ is one of the two vertices of degree~$3$ in $F$ that belongs to a $3$-cycle. Since $\cL$ is an independent set in $G$, and since $u$ and $v$ are not adjacent in $F$, the vertex $u$ is the other vertex in $F$ of degree~$3$ in $F$ that belongs to a $3$-cycle.

In all the above cases, the graph $G$ is determined. In particular, $G$ has order $n = 11$ and $\dtd(G) \le 4 \le (n-3)/2$, a contradiction.~\smallqed

\medskip
By Claim~F.3, $\ndtd(F^u) > \dtd(F)$, and so  Observation~\ref{special:ob} applies to the graph $F$ and the special vertex $u$. Analogously, by Claim~F.3, $\ndtd(F^v) > \dtd(F)$, and so Observation~\ref{special:ob} applies to the graph $F$ and the special vertex $v$. Since there is only one identified vertex in a good-graph (that belongs to the family $\cG$), we therefore have that either $F \in \{C_4, C_5\}$ or $F = B_1$ and $u$ and $v$ are the two vertices of degree~$3$ in $F$. If $F = C_4$, then $u$ and $v$ are at distance~$2$ apart in $F$, implying that $G = B_3$. If $F = C_5$, then $u$ and $v$ are at distance~$2$ apart in $F$, implying that $G = B_7$. If $F = B_1$, then $u$ and $v$ are the two vertices of degree~$3$ in $F$, implying that $G = B_4$. Thus we have shown that $G \in \{B_3,B_4,B_7\}$. This completes the proof of Claim~F.~\qed

\medskip
The proof of Lemma~\ref{path6:lem} follows immediately from Lemma~E and Lemma~F.~\qed

\medskip
By Lemma~\ref{path6:lem}, we may assume that $G$ does not contain a path on six vertices each internal vertex of which has degree~$2$ in $G$ and whose end vertices are not adjacent, for otherwise the desired result follows.

\begin{lem}
\label{degree3:lem}
If $G$ contains a vertex of degree~$3$ adjacent to ends of a $2$-handle, then $G = B_9$ or $G \in \cG$.
\end{lem}
\textbf{Proof of Lemma~\ref{degree3:lem}.}
Assume that $G$ contains a vertex of degree~$3$ adjacent to ends of a $2$-handle, $C$. Let $P$ be the $2$-path which has an end adjacent to $u$, and let $v$ be the other large vertex adjacent with an end of $P$. Let $w$ be a neighbor of $u$ on $C$, let $z$ be the neighbor of $u$ on $P$, and let $y$ be the end of $P$ different from $z$ (possibly, $y = z$). Let $C$ contain $r-1$ vertices and $P$ contain $s$ vertices. By Lemma~\ref{path6:lem}, $3\le r \le 6$, and $1\le s \le 3$.
Let $G_1 = G[V(C)\cup \{u\}\cup V(P)]$ and let $G_2 = G-V(G_1)$. Then, $G_1$ is a key $L_{r,s}$. Further, $G_2$ is a connected graph with $\delta(G_2)\ge 2$ and $v\in V(G_2)$. For $i=1,2$, let $G_{i}$ have order $n_{i}$, and so $n=n_{1}+n_{2}$. We note that $4 \le n_1 \le 9$. The following result is straightforward to verify.

\begin{unnumbered}{Claim~G}
$(n_1 - 2)/2 \le \dtd(G_1) \le (n_1 + 2)/2$. More precisely, the following holds. \\
\indent $\bullet$ $\dtd(G_1)=(n_{1}+1)/2$ if and only if $G_1=L_{4,3}$. \\
\indent $\bullet$ $\dtd(G_1)=n_{1}/2$ if and only if $G_1\in\{L_{3,1}, L_{3,3}, L_{4,2}, L_{5,3}\}$. \\
\indent $\bullet$  $\dtd(G_1)=(n_{1}-1)/2$ if and only if $G_1\in\{L_{3,2}, L_{4,1}, L_{5,2}, L_{6,1}, L_{6,3}\}$.
\\
\indent $\bullet$ $\dtd(G_1)=(n_{1}-2)/2$ if and only if $G_1\in\{L_{5,1}, L_{6,2}\}$.
\end{unnumbered}

If $G_2$ is a cycle, then $G$ is a dumb-bell, a contradiction. Hence, $G_2$ is not a cycle. In particular, $G_2 \notin \cC$. Since $G$ is an edge-minimal graph, so too is $G_2$.  By Observation~\ref{sg:ob}, $G_2 \in \cB \cup\cD \cup\cD_{b} \cup \cG$ or $\dtd(G_2) \le (n_{2}-2)/2$.
%

\begin{unnumbered}{Claim~H}
The following holds. \\
\indent {\rm (a)} $(n_1 - 1)/2 \le \dtd(G_1) \le n_1/2$. \\
\indent {\rm (b)} $(n_2 - 1)/2 \le \dtd(G_2) \le n_2/2$. \\
\indent {\rm (c)} $\dtd(G_1) = n_1/2$ or $\dtd(G_2) = n_2/2$.
\end{unnumbered}
\textbf{Proof of Claim~H} (a) Since $\dtd(G) \le \dtd(G_1) + \dtd(G_2)$, if $\dtd(G_1) \le (n_{1}-2)/2$, then $\dtd(G)\le (n_{1}-2)/2+n_{2}/2=(n-2)/2$, a contradiction. Hence, $\dtd(G_1) \ge (n_{1}-1)/2$.  It remains for us to show that $\dtd(G_1) \le n_1/2$. Suppose to the contrary that $\dtd(G_1) > n_1/2$. Then, by Claim~G, $G_1=L_{4,3}$ and $\dtd(G_1) = 3$. Further, $G$ has order~$n = n_7 + 7$. If $\dtd(G_2) \le (n_{2}-1)/2$, then $\dtd(G) \le \dtd(G_1) + \dtd(G_2) \le 3 + (n-8)/2 =(n-2)/2$, a contradiction. Therefore, $\dtd(G_2) \ge n_{2}/2$. As observed earlier, $\dtd(G_2) \le n_{2}/2$. Consequently, $\dtd(G_2) = n_{2}/2$. By Observation~\ref{f1:ob}, $G_2 \in \{B_3, D_{b}(4,4), D_{b}(3,4,1), D_{b}(3,3,2)\}$. This implies by Observation~\ref{f3b:ob} that every vertex of $G_2$ is a good-vertex. Let $S$ be a $\dtd(G_2)$-set that contains the vertex $v$. Let $w$ be a neighbor of $u$ on $C$. Then, $S \cup \{u,w\}$ is a DTD-set of $G$, and so $\dtd(G) \le |S| + 2 = \dtd(G_2)+2 = n_{2}/2+2 = (n-3)/2$, a contradiction. This establishes Part~(a).

(b) If $\dtd(G_2) \le (n_2 - 2)/2$, then $\dtd(G) \le \dtd(G_1) + \dtd(G_2) \le n_1/2 + (n_2 - 2)/2 =(n-2)/2$, a contradiction. Hence, $G_2 \in \cB \cup\cD \cup\cD_{b} \cup \cG$. Thus, by Observation~\ref{f1:ob}, $\dtd(G_2) \le n_{2}/2$.

(c) Part~(c) follows from Parts~(a) and~(b) above and the observation that $(n-1)/2 = (n_1 + n_2 - 1)/2 \le \dtd(G) \le \dtd(G_1) + \dtd(G_2)$.~\smallqed

\begin{unnumbered}{Claim~I} If $\dtd(G_1) = n_1/2$, then $G = B_9$ or $G \in \cG$.
\end{unnumbered}
\textbf{Proof of Claim~I} Suppose that $\dtd(G_1) = n_1/2$. By Claim~G, $G_1 \in \{L_{3,1}, L_{3,3}, L_{4,2}, L_{5,3}\}$.

\begin{unnumbered}{Claim~I.1}  $G_1 = L_{3,1}$.
\end{unnumbered}
\textbf{Proof of Claim~I.1} Suppose the contrary that $G_1 \in \{L_{3,3}, L_{4,2}, L_{5,3}\}$. By Claim~H, $\dtd(G_2) \le n_{2}/2$. If $G_1=L_{5,3}$, then $n_1 = 8$ and every $\dtd(G_2)$-set can be extended to a DTD-set of $G$ by adding to it the vertices $u$, $v$ and $w$, and so $\dtd(G) \le \dtd(G_2) + 3 \le n_{2}/2 + 3 = (n-2)/2$, a contradiction. Hence, $G_1 \in \{L_{3,3}, L_{4,2}\}$ and $n_1 = 6$.

We show firstly that $v$ is a bad-vertex in $G_2$. Suppose to the contrary that $v$ is a good-vertex in $G_2$ and let $S$ be a $\dtd(G_2)$-set that contains~$v$. Then, $S \cup \{u,w\}$ is a DTD-set of $G$, and so $\dtd(G) \le |S| + 2 = \dtd(G_2) + 2 = n_{2}/2+2 \le (n-2)/2$, a contradiction. Hence, $v$ is a bad-vertex in $G_2$. Therefore, by Observation~\ref{f3b:ob}, $G_2 \in \{D(4,4), D_{b}(3,4), D_{b}(3,3,1)\} \cup \cG_0$. Further, if $G_2 \in \cG_{0}$, then $v$ is the identified vertex of $G_2$, while if $G_2 \in \{D(4,4), D_{b}(3,4), D_{b}(3,3,1)\}$, then $v$ is a vertex at distance~$2$ from the central vertex of $G_2$.

Suppose $G_2 \in \{D(4,4), D_{b}(3,4), D_{b}(3,3,1)\}$. Then, $n_2 = 7$ and $v$ is a vertex at distance~$2$ from the central vertex of $G_2$. In this case, $n = 13$ and in all cases we have $\dtd(G) \le 5 = (n-2)/2$, irrespective of whether $G_1 = L_{3,3}$ or $G_1 =  L_{4,2}$. This produces a contradiction. Therefore, $G_2 \in \cG_{0}$ and $v$ is the identified vertex of $G_2$. But then by Observation~\ref{f4:ob}, there is a $\dtd(G_2)$-set, $S$, such that all neighbors of $v$ in $G_2$ belong to the set $S$. Such a set $S$ can be extended to a DTD-set in $G$ by adding to it at most two vertices of $G_1$. Hence, $\dtd(G) \le |S| + 2 \le (n_2 - 1)/2 + 2 = (n - 3)/2$, a contradiction.~\smallqed

\medskip
We now return to the proof of Claim~I. By Claim~I.1, $G_1 = L_{3,1}$. Thus, $n_1 = 4$ and $n = n_2 + 4$. By Claim~H(b) and our earlier observations, $G_2 \in \cB \cup\cD \cup\cD_{b} \cup \cG$.

\begin{unnumbered}{Claim~I.2}
The vertex $v$ is a special vertex of $G_2$ and $\ndtd(G^v_2) > \dtd(G_2)$.
\end{unnumbered}
\textbf{Proof of Claim~I.2} Suppose to the contrary that $v$ is not a special vertex in $G_2$. Then, there is an edge $f$ in $F$ which may be removed without disconnecting $F$ and with no vertex, except possibly for the vertex $v$, having degree~$1$ in $G_2$. But then $\delta(G - f) \ge 2$ and $G$ is connected, contradicting the edge-minimality of $G$. Hence, $v$ is a special vertex of $G_2$. We show next that $\ndtd(G^v_2) > \dtd(G_2)$. Suppose to the contrary that $\ndtd(G^v_2) \le \dtd(G_2)$. Let $S$ be a $\ndtd(G^v_2)$-set. In particular, $v' \in S$. The set $(S \cup \{u,z\}) \setminus \{v'\}$  is a DTD-set in $G$, and so
$\dtd(G) \le |S|+1 \le n_{2}/2 + 1 = (n-2)/2$, a contradiction.~\smallqed

\medskip
By Claim~I.2, the vertex $v$ is a special vertex of $G_2$ and $\ndtd(G^v_2) > \dtd(G_2)$. The graph $G_2$ is therefore one of the graphs listed in (b), (c) or (d) in the statement of Observation~\ref{special:ob}. We consider the three possibilities in turn.

If $G_2 = B_1$ with $v$ a vertex of degree~$3$ in $G_2$, then $G = B_9$.
If $G_2 \in \cG$ and $v$ is the identified vertex of $G_2$, then $G_2 = G_k(i,j)$ for some non-negative integers $i,j,k$, where $k \le 6$, $i + j \ge 2$ if $k = 0$ and $i + j \ge 1$ if $k \ge 1$.  In this case, $G = G_k(i+1,j) \in \cG$.
Finally, if $G_2 \in \cG_0$ and $v$ is a neighbor of the identified vertex of $G_2$, then $G_2 = G_0(i,j)$ for some non-negative integers $i,j$ where $i + j \ge 2$. If $v$ belongs to a type-1 unit in $G_2$, then $G = G_3(i-1,j) \in \cG_3$. If $v$ belongs to a type-2 unit in $G_2$, then $G = G_4(i,j-1) \in \cG_4$. In both cases, $G \in \cG$. This completes the proof of Claim~I.~\qed

\begin{unnumbered}{Claim~J} If $\dtd(G_1) = (n_1-1)/2$, then $G \in \cG_3 \cup \cG_4 \cup \cG_5$.
\end{unnumbered}
\textbf{Proof of Claim~J} Suppose that $\dtd(G_1) = (n_1-1)/2$. By Claim~G, we see that $G_1 \in \{L_{3,2}, L_{4,1}, L_{5,2}, L_{6,1}, L_{6,3}\}$. Recall that $w$ is a neighbor of $u$ on $C$, $z$ is the neighbor of $u$ on $P$, and $y$ is the end of $P$ different from $z$ (possibly, $y = z$). Let $S_1$ be a $\dtd(G_1)$-set chosen so that the neighbor of $y$ on $P$ belongs to $S_1$ (this is possible due to the structure of $G_1$).

Since $\dtd(G_1) = (n_1-1)/2$, by Claim~H we have $\dtd(G_2) = n_2/2$. Thus by Observation~\ref{f1:ob}, $G_2 \in \{B_3, D_{b}(4,4), D_{b}(3,4,1), D_{b}(3,3,2)\}$. Since the set $\cL$ is independent in $G$, we note that $G_2 \ne D_{b}(4,4)$ since the two adjacent vertices of degree~$3$ in $G_2$ will also be adjacent in $G$. Further, the vertex $v$ is a special vertex of $G_2$.

Suppose that $G_2 = B_3$. Since $G$ is edge-minimal, there are only two choices for the vertex $v$, namely, $v$ is one of the two large vertex (of degree~$3$) in $G_2$ or $v$ is a small vertex (of degree~$2$) at distance~$2$ from a large vertex in $G_2$. In both cases, the set $S_1$ can be extended to a DTD-set of $G$ by adding to it three vertices from $G_2$, and so $\dtd(G) \le \dtd(G_1) + \dtd(G_2) - 1 = (n_1 + n_2 - 3)/2 = (n-3)/2$, a contradiction. Hence, $G_2 \ne B_3$.

\begin{unnumbered}{Claim~J.1} If $G_2 = D_{b}(3,4,1)$, then $G \in \cG_4 \cup \cG_5$.
\end{unnumbered}
\textbf{Proof of Claim~J.1} Suppose that $G_2 = D_{b}(3,4,1)$. Then, $n = n_1 + 8$. The graph $G_2$ can be obtained from a path $u_{1}u_{2}\ldots u_{8}$ by joining the vertices $u_{1}$ and $u_{4}$, and the vertices $u_{6}$ and $u_{8}$. Since $v$ is a special vertex of $G_2$, we note that $v \in \{u_{2}, u_{4}, u_{5}, u_{6}\}$. If $v = u_2$, let $S = S_1 \cup \{u_4,u_5,u_6\}$, while if $v = u_6$, let $S = S_1 \cup \{u_1,u_4,u_6\}$. In both cases, $S$ is a DTD-set of $G$, and so $\dtd(G) \le |S| + 3 = \dtd(G_1) + \dtd(G_2) - 1 = (n-3)/2$, a contradiction. Hence, $v = u_4$ or $u = u_5$.
If $G_1 \in \{L_{5,2}, L_{6,1}\}$, then $n = 15$ and the set $\{u_1,u_5,u_6\}$ can be extended to a DTD-set of $G$ by adding to it three vertices from $G_1$, implying that $\dtd(G) \le 6 = (n-3)/2$, a contradiction. If $G_1 = L_{6,3}$, then $n = 17$ the set $\{u_1,u_5,u_6\}$ can be extended to a DTD-set of $G$ by adding to it four vertices from $G_1$, implying that $\dtd(G) \le 7 = (n-3)/2$, a contradiction. Hence, $G_1 = L_{3,2}$ or $G_1 = L_{4,1}$. If $v = u_4$, then $G \in \cG_4$. If $v = u_5$, then $G \in \cG_5$.~\smallqed

\begin{unnumbered}{Claim~J.2} If $G_2 = D_{b}(3,3,2)$, then $G \in \cG_3$.
\end{unnumbered}
\textbf{Proof of Claim~J.2} Suppose that $G_2 = D_{b}(3,3,2)$. Then, $n = n_1 + 8$. The graph $G_2$ can be obtained from a path $u_{1}u_{2}\ldots u_{8}$ by joining the vertices $u_{1}$ and $u_{3}$, and the vertices $u_{6}$ and $u_{8}$. Since $v$ is a special vertex of $G_2$, we note that $v \in \{u_{3}, u_{4}, u_{5}, u_{6}\}$. Renaming vertices if necessary, we may assume that $v = u_3$ or $v = u_4$. If $v = u_3$, then the set $S_1 \cup \{u_3,u_5,u_6\}$ is a DTD-set of $G$, implying that $\dtd(G) \le \dtd(G_1) + 3 = \dtd(G_1) + \dtd(G_2) - 1 = (n-3)/2$, a contradiction. Hence, $v = u_4$.
If $G_1 \in \{L_{5,2}, L_{6,1}\}$, then $n = 15$ and the set $\{u_3,u_5,u_6\}$ can be extended to a DTD-set of $G$ by adding to it three vertices from $G_1$, implying that $\dtd(G) \le 6 = (n-3)/2$, a contradiction. If $G_1 = L_{6,3}$, then $n = 17$ the set $\{u_3,u_5,u_6\}$ can be extended to a DTD-set of $G$ by adding to it four vertices from $G_1$, implying that $\dtd(G) \le 7 = (n-3)/2$, a contradiction. Hence, $G_1 = L_{3,2}$ or $G_1 = L_{4,1}$. In both cases, $G \in \cG_3$.~\smallqed

\medskip
Claim~J now follows from Claim~J.1 and Claim~J.2.~\qed

\medskip
We now return to the proof of Lemma~\ref{degree3:lem}. By Claim~H, $\dtd(G_1) = n_1/2$ or $\dtd(G_1) = (n_1-1)/2$. If $\dtd(G_1) = n_1/2$, then by Claim~I, $G = B_9$ or $G \in \cG$. If $\dtd(G_1) = (n_1-1)/2$, then by Claim~J, $G \in \cG_4 \cup \cG_5$. This completes the proof of Lemma~\ref{degree3:lem}.~\qed

\medskip
By Lemma~\ref{degree3:lem}, we may assume that if $G$ contains a vertex adjacent to ends of a $2$-handle, then such a vertex has degree at least~$4$ in $G$, for otherwise the desired result follows. By our earlier observations, every $2$-handle in $G$ has order at most~$5$.

\begin{lem}\label{34h:lem}
Every $2$-handle in $G$ has order~$2$ or~$5$.
\end{lem}
\textbf{Proof of Lemma~\ref{34h:lem}.} Suppose that $G$ contains a 2-handle $C$ with $|V(C)| \in \{3,4\}$. Let $C$ be the path $v_{1}v_{2} \ldots v_{t}$, where $t\in\{3,4\}$ and let $v$ be the vertex in $G$ adjacent to both ends of $C$ (the vertices $v_{1}$ and $v_{t}$). By our assumptions to date, $d_{G}(v)\ge 4$. Let $N_{v}=N(v)\setminus\{v_{1}, v_{t}\}$. By Lemma~\ref{i:lem}, $N_{v}$ comprises only small vertices.  We note that $|N_{v}|\ge 2$.

Let $G'=G-V(C)$. The graph $G'$ is a connected subgraph of $G$ with $\delta(G')\ge 2$. Let $G'$ have order $n'=n-t$. The degree of every large vertex different from $v$ is unchanged in $G$ and $G'$, and so $G'$ contains at least one large vertex. Hence, $G'$ is not a cycle and, by Observation~\ref{sg:ob}, $G'\in\cB\cup\cD\cup\cD_{b}\cup\cG$ or $\dtd(G')\le (n'-2)/2$. We proceed with the following claim.

\begin{unnumbered}{Claim K}
$G'\in\cB\cup\cD\cup\cD_{b}\cup\cG$.
\end{unnumbered}
\textbf{Proof of Claim K.} Suppose to the contrary that $\dtd(G')\le (n'-2)/2$. Let $S$ be a $\dtd(G')$-set. If $t=4$, then $S\cup\{v,v_{1}\}$ is a DTD-set of $G$, and so $\dtd(G) \le |S| + 2 = \dtd(G')+2 = (n'-2)/2+2=(n-2)/2$, a contradiction. Hence, $t=3$.

\begin{unnumbered}{Claim K.1}
$N_{v}$ is an independent set in $G$.
\end{unnumbered}
\textbf{Proof of Claim K.1} Suppose there are two adjacent vertices, $w_{1}$ and $w_{2}$, in $N_{v}$. Let $H = G-\{w_{1}, w_{2}\}$ and let $H$ have order $n_{H}=n-2$. Then, $H$ is a connected graph with minimum degree at least~$2$. By Observation~\ref{sg:ob}, either $\dtd(H)\le (n_{H}-2)/2$ or $H\in\cB\cup\cC\cup\cD_{b}\cup\cD\cup\cG$. We note that $d_H(v) \ge 3$ and that the degree of each large vertex in $G$ other than $v$ is unchanged in $H$. Since $G$ contains at least two large vertices, the graph $H$ is therefore not a cycle. Hence, by Observation~\ref{f1:ob}, $\dtd(H)\le n_{H}/2$. By considering the cycle $vv_1v_2v_3v$ in $H$, we can choose a $\dtd(H)$-set to contain $v$ or to contain at least two vertices in $N_{H}(v)$. In both cases this set is a DTD-set of $G$. Hence, $\dtd(G)\le \dtd(H)\le n_{H}/2=(n-2)/2$, a contradiction.~\smallqed

\begin{unnumbered}{Claim K.2.}
The only common neighbor of two vertices in $N_{v}$ is the vertex~$v$.
\end{unnumbered}
\textbf{Proof of Claim K.2} Suppose to the contrary that two vertices $w_{1}$ and $w_{2}$ in $N_{v}$ have a common neighbor, $x$, different from~$v$. We show first that $d_{G}(x) \ge 3$. Suppose to the contrary that $d_{G}(x)=2$. In this case we consider the connected graph $H=G-\{w_{1}, w_{2}, x\}$. Let $H$ have order $n_{H}=n-3$. Since $|\cL|\ge 2$ in $G$, $H$ has at least two large vertices, and so $d_{H}(v)\ge 3$. Thus, $H$ is not a cycle. Applying the inductive hypothesis and Observation~\ref{f1:ob} to the graph $H$, we have that $\dtd(H)\le n_{H}/2$. By considering the cycle $vv_1v_2v_3v$ in $H$, we can choose a $\dtd(H)$-set to contain $v$ or to contain at least two vertices in $N_{H}(v)$. In both cases adding the vertex $w_{1}$ to such a set produces a DTD-set of $G$, and so $\dtd(G) \le \dtd(H) + 1$. If $\dtd(H)\le (n_{H}-1)/2$, then $\dtd(G) \le (n_{H}-1)/2 + 1=(n-2)/2$, a contradiction. Hence, $\dtd(H)=n_{H}/2$. Since $vv_1v_2v_3v$ is a cycle in $H$ and $d_H(v) \ge 3$, by Observation~\ref{f1:ob} we deduce that $H \in\{D_{b}(4,4), D_{b}(3,4,1)\}$. If $H=D_{b}(4,4)$, then $\cL$ is not an independent set in $G$, a contradiction. If $H=D_{b}(3,4,1)$, then the graph $G$ is determined and $\dtd(G) = 4 = (n-3)/2$, a contradiction. Hence, $d_{G}(x) \ge 3$.

We now consider the graph $H=G-w_{1}$. The degree of each vertex in $H$ different from $v$ and $x$ remains unchanged from its degree in $G$. Further, $d_H(v) \ge 3$ and $d_H(x) \ge 2$, and so $H$ is a connected graph with $\delta(H)\ge 2$. Since $H$ is a subgraph of $G$, $H$ is edge-minimal. Let $H$ have order $n_{H} = n-1$. By considering the cycle $vv_1v_2v_3v$ in $H$, we can choose a $\dtd(H)$-set to contain $v$ or to contain at least two vertices in $N_{H}(v)$. Such a $\dtd(H)$-set is a DTD-set of $G$ implying that $\dtd(G) \le \dtd(H)$. Applying the inductive hypothesis to the edge-minimal graph $H$ (which we recall is not a cycle), $\dtd(H) \le n_{H}/2$. If $\dtd(H)\le (n_{H}-1)/2$, then $\dtd(G)\le (n-2)/2$, a contradiction. Hence, $\dtd(H)=n_{H}/2$. Since $vv_1v_2v_3v$ is a cycle in $H$ and $d_H(v) \ge 3$, by Observation~\ref{f1:ob} we deduce that $H = D_{b}(3,4,1)$. But then, $G$ is determined and $\dtd(G) = 3 = (n-3)/2$, a contradiction.~\smallqed

\medskip
We now return to the proof of Claim~K. Let $G^*$ be the graph obtained from $G'-v$ by selecting a vertex $w \in N_{v}$ (of degree~$1$ in $G' - v$) and joining it to every other vertex in $N_v$. Let $G^{*}$ have order $n^{*}=n'-1=n-4$. We note that by Claim~K.1 and Claim~K.2, and by the fact that $G'$ is edge-minimal, the graph $G^{*}$ is edge-minimal. By construction $G^*$ has $|\cL| - 1 \ge 1$ large vertices, implying that $G^*$ is not a cycle. Hence, by the inductive hypothesis $G^{*}\in\cB\cup\cD\cup\cD_{b}\cup\cG$ or $\dtd(G^{*})\le (n^{*}-2)/2$.

Suppose $\dtd(G^{*})\le (n^{*}-2)/2$. Let $S^*$ be a $\dtd(G^{*})$-set. If $|S^* \cap N_v| = 1$, let $S = S^* \cup \{v,v_1\}$. If $|S^* \cap N_v| \ge 2$, let $S = S^* \cup \{v_1,w\}$. In both cases, $S$ is a NTD-set of $G$ and $|S| \le |S^*| + 2$, implying that $\dtd(G)\le \dtd(G^{*})+2 = (n^{*}-2)/2+2=(n-2)/2$, a contradiction. Hence, $G^{*}\in \cB\cup\cD\cup\cD_{b}\cup\cG$. Since $G^{*}$ is not a cycle, Observation~\ref{f1:ob} implies that $\dtd(G^{*}) \le n^{*}/2$.

Let $x \in N_v \setminus \{w\}$ and consider the edge $e = wx$ that was added to $G' - v$ when constructing $G^*$. If $e$ is a good-edge in $G^*$, then every $\dtd(G^{*})$-set containing $w$ and $x$ may be extended to a DTD-set of $G$ by adding to it the vertex $v_{1}$, implying that $\dtd(G) \le \dtd(G^{*}) + 1 \le n^*/2 + 1 = (n-2)/2$, a contradiction. Therefore, the edge $e$ is a bad-edge of $G^*$. By Observation~\ref{f3a:ob}, $G^{*}\in\{D(4,4), D_{b}(3,4), D_{b}(3,3,1)\} \cup \cG \setminus \cG_{1}$.

Suppose $G^{*}\in\{D(4,4), D_{b}(3,4), D_{b}(3,3,1)\}$. If $G^{*}=D(4,4)$, then the edge $e$ satisfies Observation~\ref{f3a:ob}(a) and the graph $G$ is determined. In this case, $n=11$ and $\dtd(G) \le 4 = (n-3)/2$, a contradiction.
If $G^{*}=D_{b}(3,4)$, then the edge $e$ satisfies Observation~\ref{f3a:ob}(b), implying that $\cL$ is not an independent set in $G$, a contradiction.
If $G^{*}=D_{b}(3,3,1)$, then the edge $e$ satisfies Observation~\ref{f3a:ob}(c). In this case, we contradict Claim~K.1 or Claim~K.2.
Since all three cases produce a contradiction, $G^* \in \cG \setminus \cG_{1}$.

Let $z$ be the identified vertex of $G^{*}$. By Observation~\ref{f3a:ob}(d), the edge $e = wx$ is incident with the identified vertex $z$. Hence either $w = z$ or $x = z$. In both cases, when we reconstruct the original graph $G$ from $G^*$ the component of $G - vz$ that contains the vertex~$z$ contains at least one $2$-handle whose ends are adjacent to a vertex of degree~$3$ in $G$, contradicting our earlier assumption that a vertex of $G$ adjacent to the ends of a $2$-handle has degree at least~$4$ in $G$. This completes the proof of Claim~K.~\qed

\medskip
We now return to the proof of Lemma~\ref{34h:lem}. By Claim~K, $G'\in\cB\cup\cD\cup\cD_{b}\cup\cG$.  Since $G'$ is not a cycle, Observation~\ref{f1:ob} implies that $\dtd(G')\le n'/2$. If $v$ is a good-vertex in $G'$, then every $\dtd(G')$-set can be extended to DTD-set of $G$ by adding to it the vertex $v_1$, implying that in this case $\dtd(G) \le \dtd(G') + 1 \le n'/2+1$.

\newpage
\begin{unnumbered}{Claim L}
$t=3$.
\end{unnumbered}
\textbf{Proof of Claim L.} Suppose to the contrary that $t=4$. Then, $n'=n-4$. If $v$ is a good-vertex in $G'$, then $\dtd(G) \le n'/2+1 = (n-2)/2$, a contradiction. Hence, $v$ is a bad-vertex in $G$, and so, by  Observation~\ref{f3b:ob}, $G\in\{D(4,4), D_{b}(3,4), D_{b}(3,3,1)\} \cup \cG_{0}$.
Suppose $G'\in \cG_{0}$. Then, $G'=G_{0}(n_{1},n_{2}) \in \cG_{0}$, where $n_{1}+n_{2}\ge 2$. By Observation~\ref{f3b:ob}(a), the vertex $v$ is the identified vertex of $G$. We note that $G'$ contains $n_1 + n_2$ $2$-handles whose ends are adjacent to a vertex of degree~$3$ in $G'$. Reconstructing the original graph $G$ from $G'$ we therefore note that $G$ contains $n_1 + n_2$ $2$-handles whose ends are adjacent to a vertex of degree~$3$ in $G$, a contradiction. Hence, $G' \in \{D(4,4), D_{b}(3,4), D_{b}(3,3,1)\}$. By Observation~\ref{f3b:ob}(b), $v$ is a vertex at distance~$2$ from the central vertex in $G'$. If $G'=D_{b}(3,4)$, then $\cL$ is not a independent set in $G$, a contradiction. Hence, $G' \in \{D(4,4), D_{b}(3,3,1)\}$. In both cases, the graph $G$ is determined and $\dtd(G) = 4 = (n-3)/2$, a contradiction.~\smallqed

\medskip
By Claim~L, $t = 3$. Thus, $n'=n-3$. Suppose $\dtd(G')=n'/2$. By Observation~\ref{f1:ob} and by our assumption that $\cL$ is an independent set in $G$, $G' \in \{B_3, D_{b}(3,4,1), D_{b}(3,3,2)\}$. Since $v \in \cL$ and $\cL$ is an independent set in $G$, the vertex $v$ is not adjacent to a large vertex in $G'$. If $G'=B_{3}$, then, up to isomorphism, there are two possible graphs $G$. In both cases, $\dtd(G) = 4 = (n-3)/2$, a contradiction. Hence $G' \in \{D_{b}(3,4,1), D_{b}(3,3,2)\}$. In both cases, however we reconstruct the original graph $G$, we produce a $2$-handle whose ends are adjacent to a vertex of degree~$3$ in $G$, a contradiction. Hence, $\dtd(G) \le (n'- 1)/2$.

If $v$ is a good-vertex in $G'$, then $\dtd(G) \le (n'-1)/2+1=(n-2)/2$, a contradiction. Hence, $v$ is a bad-vertex in $G'$. By Observation~\ref{f3b:ob}, $G' \in\{D(4,4), D_{b}(3,4), D_{b}(3,3,1)\}\cup\cG_{0}$. If $G'\in\cG_{0}$, then analogously as in the proof of Claim~L the graph $G$ would contain a $2$-handle whose ends are adjacent to a vertex of degree~$3$ in $G$, a contradiction. If $G'=D_{b}(3,4)$, then $\cL$ is not a independent set in $G$, a contradiction. Hence, $G' \in \{D(4,4), D_{b}(3,3,1)\}$. By Observation~\ref{f3b:ob}(b), $v$ is a vertex at distance~$2$ from the central vertex in $G'$. In both cases, the graph $G$ is determined and $\dtd(G) = 4 = (n-2)/2$, a contradiction. This completes the proof of Lemma~\ref{34h:lem}.~\qed

\medskip
By Lemma~\ref{34h:lem}, every $2$-handle in $G$ has order~$2$ or~$5$. We show next that in fact $G$ contains no $2$-handle.

\begin{lem}
There is no $2$-handle in $G$.
\label{no2h:lem}
\end{lem}
\textbf{Proof of Lemma~\ref{no2h:lem}.} Suppose to the contrary that $G$ contains a $2$-handle, $C$. Then, $|V(C)| \in \{2,5\}$. Let $C$ be the path $v_{1}v_{2} \ldots v_{t}$, where $t \in \{2,5\}$ and let $v$ be the vertex in $G$ adjacent to both ends of $C$. By our assumptions to date, $d_{G}(v) \ge 4$ and every neighbor of $v$ has degree~$2$.
Let $F=G-V(C)$ and let $F$ have order $n_{F}=n-t$. By construction, $F$ is a connected subgraph of $G$ with $\delta(F)\ge 2$. By Observation~\ref{sg:ob}, $F \in \cB\cup\cD\cup\cD_{b} \cup \cG$ or $\dtd(F)\le (n_{F}-2)/2$.

Suppose $\dtd(F)\le (n_F-2)/2$. If $t=2$, then $n_F = n - 2$ and every $\dtd(F)$-set can be extended to a DTD-set of $G$ by adding to it the vertex $v$, implying that $\dtd(G) \le \dtd(F) + 1 \le (n_{F}-2)/2+1 = (n-2)/2$, a contradiction. If $t=5$, then $n_F = n - 5$ and every $\dtd(F)$-set can be extended to a DTD-set of $G$ by adding to it the vertices $v$ and $v_{1}$, implying that $\dtd(G) \le \dtd(F) + 2 \le (n_{F}-2)/2 + 2 \le (n-3)/2$, a contradiction. Hence, $F\in\cB\cup\cD\cup\cD_{b}\cup\cG$. Since $|\cL| \ge 2$, the graph $F$ has at least one large vertex different from~$v$, implying that $F$ is not a cycle. Hence by Observation~\ref{f1:ob}, $\dtd(F) \le n_{F}/2$.

We show that $v$ is a bad-vertex in $F$. Suppose to the contrary that $v$ is a good-vertex in~$F$. If $t = 2$, then every $\dtd(F)$-set which contains $v$ is a DTD-set of $G$, implying in this case that $\dtd(G) \le \dtd(F) \le n_{F}/2 = (n-2)/2$, a contradiction. If $t = 5$, then every $\dtd(F)$-set which contains $v$ can be extended to a DTD-set of $G$ by adding to it the vertex $v_1$, implying in this case that $\dtd(G) \le \dtd(F) + 1 \le n_{F}/2 - 1 = (n-3)/2$, a contradiction. Therefore, $v$ is a bad-vertex in $F$.

By Observation~\ref{f3b:ob}, either $F \in \{D(4,4), D_{b}(3,4), D_{b}(3,3,1)\}$ and $v$ is the vertex at distance~$2$ from the central vertex in $F$, or $F \in\cG_{0}$ and $v$ is the identified vertex of $F$. In all cases, the graph $G$ is determined. If $F \in \{D_{b}(3,4), D_{b}(3,3,1)\}$, then the set $\cL$ is not independent, a contradiction. If $F \in \cG_{0}$, then $G$ contains a vertex of degree~$3$ adjacent to both ends of a $2$-handle, a contradiction. If $F = D(4,4)$, then either $t = 2$, in which case $\dtd(G) = 3 = (n-3)/2$, or $t = 5$, in which case $\dtd(G) \le 5 = (n-2)/2$. Both cases produce a contradiction. This completes the proof of Lemma~\ref{no2h:lem}.~\qed

\medskip
By Lemma~\ref{no2h:lem}, there is no $2$-handle in $G$. Recall that by our earlier assumptions, $n \ge 8$ and the set $\cL$ is an independent set. Further, $G$ does not contain a path on six vertices each internal vertex of which has degree~$2$ in $G$ and whose end vertices are not adjacent. In particular, every $2$-path in $G$ has order at most~$3$.

\begin{lem}
\label{noc4:lem}
If $G$ contains a $4$-cycle, then $G = B_{8}$.
\end{lem}
\textbf{Proof of Lemma~\ref{noc4:lem}.} Suppose that $G$ contains a $4$-cycle $uvwxu$. Renaming vertices if necessary, we may assume that $u$ and $w$ are large vertices of $G$ (and so $v$ and $x$ are small vertices). We now consider the connected subgraph $G'=G-v$ of $G$ that satisfies $\delta(G')\ge 2$. Let $G'$ have order $n'$, and so $n' = n-1 \ge 7$. By Observation~\ref{sg:ob}, $G'\in \cB\cup\cC\cup\cD\cup\cD_{b}\cup\cG$ or $\dtd(G')\le (n'-2)/2$. Every DTD-set in $G'$ is also a DTD-set in $G$, implying that $\dtd(G) \le \dtd(G')$.  If $\dtd(G')\le (n'-1)/2$, then $\dtd(G)\le (n-2)/2$, a contradiction. Hence, $\dtd(G') \ge n'/2$. If $G'$ is a cycle, then $G$ would contain a $2$-path of order at least~$4$, a contradiction. Hence by Observation~\ref{f1:ob}, $G' \in \{B_{3}, D_{b}(4,4), D_{b}(3,4,1), D_{b}(3,3,2)\}$. If $G'\in\{D_{b}(4,4), D_{b}(3,4,1), D_{b}(3,3,2)\}$, then $G$ contains a $2$-handle, a contradiction. Hence, $G' = B_3$, implying that $G = B_8$.~\qed

\medskip
By Lemma~\ref{noc4:lem}, we may assume that $G$ contains no $4$-cycle, for otherwise the desired result follows.

\begin{lem}
\label{noc5:lem}
There is no $5$-cycle in $G$.
\end{lem}
\textbf{Proof of Lemma~\ref{noc5:lem}.} Suppose to the contrary that $G$ contains a $5$-cycle $uvwxyu$. Renaming vertices if necessary, we may assume that $u$ and $w$ are large vertices of $G$ (and so $v$, $x$ and $y$ are small vertices). We now consider the connected subgraph $G'=G-v$ of $G$ that satisfies $\delta(G') \ge 2$. Let $G'$ have order $n'$, and so $n' = n-1 \ge 7$. By Observation~\ref{sg:ob}, $G'\in \cB\cup\cC\cup\cD\cup\cD_{b}\cup\cG$ or $\dtd(G')\le (n'-2)/2$. Every DTD-set in $G'$ is also a DTD-set in $G$, implying that $\dtd(G) \le \dtd(G')$.  If $\dtd(G')\le (n'-1)/2$, then $\dtd(G)\le (n-2)/2$, a contradiction. Hence, $\dtd(G') \ge n'/2$. By Observation~\ref{f1:ob} and since every $2$-path in $G$ has order at most~$3$, $G' \in \{C_7,B_{3}, D_{b}(4,4), D_{b}(3,4,1), D_{b}(3,3,2)\}$. We note that $uyxw$ is an induced path in $G'$ where $x$ and $y$ have degree~$2$ in $G'$. This implies that $G' \notin \{D_{b}(4,4), D_{b}(3,4,1)\}$. If $G' = C_7$, then $\dtd(G) \le |\{u,v,w\}| = (n-2)/2$, a contradiction. If $G' = D_{b}(3,3,2)$, then $u$ and $w$ are necessarily the two large vertices in $G'$. But then $G$ contains a $2$-handle, a contradiction. If $G' = B_3$, then $G$ contains a $4$-cycle, a contradiction.~\qed

\medskip
By Lemma~\ref{noc5:lem} and our assumptions to date, we may assume that a shortest cycle in $G$ has length at least~$6$; that is, $G$ has girth at least~$6$. Recall that every $2$-path in $G$ has order at most~$3$.

\begin{lem}
\label{no2path3:lem}
If $G$ has a $2$-path of order~$3$, then $G = B_6$.
\end{lem}
\textbf{Proof of Lemma~\ref{no2path3:lem}.} Let $P \colon v_1v_2v_3$ be a $2$-path in $G$ and let $u$ and $v$ be the large vertices adjacent to $v_1$ and $v_3$, respectively.

\begin{unnumbered}{Claim~M}
If $u$ and $v$ do not have a common neighbor, then $G = B_6$.
\end{unnumbered}
\textbf{Proof of Claim~M.}
Suppose that $u$ and $v$ do not have a common neighbor. We show that $G = B_6$. Let $N_u^1$ and $N_u^2$ be the set of vertices at distance~$1$ and~$2$, respectively, from $u$ in $G - V(P)$. Let $N_v^1$ and $N_v^2$ be defined analogously. By assumption, $N_u^1 \cap N_v^1 = \emptyset$. Every neighbor of a large vertex is a small vertex. In particular, every vertex in $N_u \cup N_v$ has degree~$2$.  Let $G'$ be the graph of order $n'=n-4$ obtained from $G - (V(P)\cup\{u\})$ by joining the vertices in $N_u^1$ to the vertex $v$. Since $G$ is an edge-minimal graph, so too is the graph $G'$. Since $|N_u^1|\ge 2$ and $|N_v^1| \ge 2$, we note that $d_{G'}(v) = |N_u^1| + |N_v^1| \ge 4$. Hence, $G'$ is not a cycle and $G'$ contains a vertex of degree at least~$4$. Therefore by the induction hypothesis, either $G'\in \{B_{4}, B_{8}, B_{9}, B_{10}\} \cup \cD \cup \cG$ or $\dtd(G')\le (n'-2)/2$.

\begin{unnumbered}{Claim~M.1}
$\dtd(G') \ge (n'-1)/2$.
\end{unnumbered}
\textbf{Proof of Claim~M.1} Suppose to the contrary that $\dtd(G') \le (n'-2)/2$. Let $S'$ be a $\dtd(G')$-set. By construction, $N_{G'}(v) = N_u^1 \cup N_v^1$. Suppose that $S'$ totally dominates the vertex $v$. Thus, $S'$ contains at least one vertex in $N_u^1 \cup N_v^1$.
If $S'$ contains a vertex in $N_u^1$ and a vertex in $N_v^1$, let $S = S' \cup \{u,v\}$.
If $S'$ contains a vertex of $N_v^1$ and no vertex of $N_u^1$, let $S = S' \cup \{u,v_1\}$.
If $S'$ contains no vertex of $N_v^1$ and $v \in S'$, let $S = S' \cup \{u,v_3\}$.
If $S'$ contains no vertex of $N_v^1$ and $v \notin S'$, let $S = S' \cup \{v,v_3\}$.
In all four cases, the set $S$ is a DTD-set of $G$  and $|S|= |S'|+2$, implying that $\dtd(G) \le |S| = |S'|+2 = (n'-2)/2 + 2 = (n-2)/2$, a contradiction. Hence, $S'$ does not totally dominates the vertex $v$.

Since $S'$ contains no vertex in $N_u^1 \cup N_v^1$, the set $S'$ disjunctively dominates $v$ and therefore contains at least two vertices in $N_u^2 \cup N_v^2$.
If $v \in S'$ and $S'$ contains no vertex in $N_u^2$, let $S = S' \cup \{u,v_1\}$.
If $v \in S'$ and $S'$ contains exactly one vertex in $N_u^2$, let $S = S' \cup \{u,v_2\}$.
If $v \in S'$ and $S'$ contains at least two vertices in $N_u^2$, let $S = S' \cup \{u,v_3\}$.
If $v \notin S'$ and $S'$ contains no vertex in $N_u^2$, let $S = S' \cup \{v_1,v_2\}$.
If $v \notin S'$ and $S'$ contains exactly one vertex in $N_u^2$, let $S = S' \cup \{u,v_2\}$.
If $v \notin S'$ and $S'$ contains at least two vertices in $N_u^2$, let $S = S' \cup \{v_2,v_3\}$.
In all six cases, the set $S$ is a DTD-set of $G$  and $|S|= |S'|+2$, implying that $\dtd(G) \le |S| = |S'|+2 = (n'-2)/2 + 2 = (n-2)/2$, a contradiction.~\smallqed

\medskip
By Claim~M, $G'\in \{B_{4}, B_{8}, B_{9}, B_{10}\} \cup \cD \cup \cG$. Since $\cL$ is an independent set in $G$, by construction the set of large vertices in $G'$, namely the set $\cL \setminus \{u\}$, form an independent set in $G'$. If $G' \in \cG$, then since $v$ has degree at least~$4$ in $G'$ and since the large vertices in $G'$ are independent, $G' \in \cG_0 \cup \cG_1$ (and $G'$ consists only of type-1 units). But then the graph $G$ would contain a $2$-handle, a contradiction. Hence, $G' \notin \cG$.

If $G' = B_{4}$, then $v$ is one of the two vertices of degree~$4$ in $G'$, implying that the graph $G$ contains a $2$-path on at least four vertices as well as a $4$-cycle, a contradiction. If $G'\in \{B_{8}, B_{9}, B_{10}\}$, then $v$ is the vertex of degree~$4$ in $G'$. If $G'\in \{B_{9}, B_{10}\}$, then $G$ would contain a $2$-handle, a contradiction. If $G' = B_8$, then since $G$ contains no $4$-cycles, there is only one way to reconstruct the graph $G$ from $G'$. In this case, $n = 13$ and the four large vertices in $G$ form a DTD-set of $G$, implying that $\dtd(G) \le 4 = (n-5)/2$, a contradiction. Hence, $G' \notin \{B_{4}, B_{8}, B_{9}, B_{10}\}$. Therefore, $G \in \cD$. But then $\cL = \{u,v\}$ and $\{u,v_1,v_3,v\}$ is a DTD-set of $G$, implying that $\dtd(G) \le 4$. If $n \ge 10$, then $\dtd(G) \le (n-2)/2$, a contradiction. Hence, $n = 9$, implying that $G' = D(3,3)$ and $G = B_6$. This completes the proof of Claim~M.1\smallqed

\medskip
By Claim~M.1, we may assume that $u$ and $v$ have a common neighbor, $w$ say, for otherwise $G = B_6$, and the desired result follows. Since $G$ has no $4$-cycle, the vertex $w$ is the only common neighbor of $u$ and $v$. Let $G'$ be the graph obtained from $G - (V(P) \cup \{w\})$ by adding the edge $e = uv$. Then, $G'$ is a connected graph with $\delta(G') \ge 2$. Let $G'$ have order~$n'$, and so $n' =  n- 4$.

\begin{unnumbered}{Claim~M.2}
The graph $G' - e$ is an edge-minimal graph.
\end{unnumbered}
\textbf{Proof of Claim~M.2} Suppose to the contrary that $G' - e = G - (V(P) \cup \{w\})$ is not an edge-minimal graph. Then $G' - e$ is disconnected or at least one of $u$ or $v$ has degree~$1$ in $G' - e$. This implies that the graph $G'$ is an edge-minimal graph. Then, $\dtd(G') \le (n'-2)/2$ or $G' \in \cB \cup \cC \cup \cD \cup \cD_{b} \cup \cG$. Let $S'$ be a $\dtd(G')$-set. Suppose $\dtd(G') \le (n'-2)/2$. If $u$ or $v$ belong to $S'$, let $S = S' \cup \{u,v,w\}$. Suppose that neither $u$ nor $v$ belong to $S'$. If both $u$ and $v$ have a neighbor in $S'$, let $S = S' \cup \{u,v\}$. If $u$ has a neighbor in $S'$, let $S = S' \cup \{v_1,w\}$. If $v$ has a neighbor in $S'$, let $S = S' \cup \{v_3,w\}$. If neither $u$ nor $v$ have a neighbor in $S'$, let $S = S' \cup \{u,v\}$. In all the above cases, the set $S$ is a DTD-set of $G$ and $|S| \le |S'| - 2$, implying that $\dtd(G) \le |S| \le (n'-2)/2 = (n-2)/2$, a contradiction. Hence, $G' \in \cB \cup \cC \cup \cD \cup \cD_{b} \cup \cG$.

If $G'$ is a cycle, then since $G$ has girth at least~$6$ and every $2$-path in $G$ has order at most~$3$, we deduce that $G' = C_5$. But then the graph $G$ is determined. In this case, $n = 9$ and $\{u,v,w\}$ is a DTD-set of $G$, and so $\dtd(G) = 3 = (n-3)/2$, a contradiction. Hence, $G'$ is not a cycle, implying that $\dtd(G') \le n'/2$.
If $e$ is a good-edge of $G'$, then choosing $S'$ to be a $\dtd(G')$-set that contains both $u$ and $v$, the set $S' \cup \{w\}$ is a DTD-set of $G$, implying that $\dtd(G) \le |S'| + 1 \le n'/2 + 1 = (n-2)/2$, a contradiction. Hence, $e$ is a bad-edge of $G'$. However applying Observation~\ref{f3a:ob} to the graph $G'$ and the bad-edge $e$ of $G'$, the graph $G$ necessarily contains a $3$-cycle or a $4$-cycle, a contradiction.~\smallqed

\medskip
By Claim~M.2, the graph $G' - e = G - (V(P) \cup \{w\})$ is an edge-minimal graph. Thus, $\dtd(G' - e) \le (n'-2)/2$ or $G'-e \in \cB \cup \cC \cup \cD \cup \cD_{b} \cup \cG$. Let $S'$ be a $\dtd(G'-e)$-set. If $\dtd(G'-e) \le (n'-2)/2$, then $S' \cup \{u,v\}$ is a DTD-set of $G$, implying that $\dtd(G) \le |S| \le (n'-2)/2 = (n-2)/2$, a contradiction. Hence, $G'-e \in \cB \cup \cC \cup \cD \cup \cD_{b} \cup \cG$. Since both $u$ and $v$ have degree at least~$2$ in $G' - e$, we note that $G' \notin \cD$. If $G'-e \in \cD_{b} \cup \cG$, then $G$ would contain a $3$-cycle or a $4$-cycle, a contradiction. If $G'-e$ is a cycle, then since $G$ has girth at least~$6$ and every $2$-path in $G$ has order at most~$3$, we note that $G' \in \{C_6,C_7\}$ and $u$ and $v$ are at distance~$3$ apart on the cycle. But then the graph $G$ is determined. In this case, $n \in \{10,11\}$ and $\{u,v,w\}$ is a DTD-set of $G$, and so $\dtd(G) = 3 \le (n-4)/2$, a contradiction. If $G'-e \in \cB$, then since $G$ has girth at least~$6$, $G'-e = B_6$ and $u$ and $v$ are the two vertices of degree~$3$ in $G' - e$. But then $n = 13$ and  $\{u,v,w\}$ is a DTD-set of $G$, and so $\dtd(G) = 3 = (n-7)/2$, a contradiction. This completes the proof of Lemma~\ref{no2path3:lem}.~\qed

\medskip
By Lemma~\ref{no2path3:lem}, if $G$ has a $2$-path of order~$3$, then $G = B_6$. Hence we may assume that every $2$-path in $G$ has order~$1$ or~$2$. Thus every small vertex (of degree~$2$) has either two large neighbors or one large neighbor and one small neighbor depending on whether it belongs to a $2$-path of order~$1$ or a $2$-path of order~$2$, respectively.

Let $\cS = (\cS_1,\cS_2)$ be a weak partition of $\cS$ (a partition where some of the sets may be empty), where $\cS_1$ is the set of small vertices with two large neighbors and $\cS_2$ is the set of small vertices with exactly one large neighbor. We note that $G[\cS_2]$ consists of the disjoint union of paths of order~$2$.

Let $\cL = (\cL_0,\cL_1,\cL_2)$ be a weak partition of the large vertices $\cL$, where $\cL_0$, $\cL_1$ and $\cL_2$ are the set of large vertices adjacent to zero, one and at least two vertices in $\cS_1$, respectively, respectively. Let $\cS_{1,1}$ be the set of vertices in $\cS_1$ with both neighbors in $\cL_1$ and let $\cS_{1,2}$ be the set of vertices in $\cS_1$ with exactly one neighbor in $\cL_1$ (and the other neighbor in $\cL_2$). Further, for $i \in \{1,2\}$ let $\cL_{1,i}$ be the set of vertices in $\cL_1$ adjacent to a vertex in $\cS_{1,i}$. Thus, $(\cL_{1,1},\cL_{1,2})$ is a partition of $\cL_1$.

Let $|\cS| = s$, $|\cS_1| = s_1$ and $|\cS_2| = 2s_2$, and so $s = s_1 + 2s_2$. Let $|\cS_{1,1}| = s_{1,1}$ and $|\cS_{1,2}| = s_{1,2}$, and so $s_1 \ge s_{1,1} + s_{1,2}$. Let $|\cL| = \ell$,  $|\cL_0| = \ell_0$, $|\cL_1| = \ell_1$, and $|\cL_2| = \ell_2$, and so $\ell = \ell_0 + \ell_1 + \ell_2$. Let $|\cL_{1,1}| = \ell_{1,1}$ and $|\cL_{1,2}| = \ell_{1,2}$, and so $\ell_1 = \ell_{1,1} + \ell_{1,2}$ and
\begin{equation}
\ell = \ell_0 + \ell_{1,1} + \ell_{1,2} + \ell_2.
 \label{Eq1}
\end{equation}

The subgraph $G[\cL_{1,1} \cup \cS_{1,1}]$ induced by $\cL_{1,1} \cup \cS_{1,1}$ consists of a disjoint union of paths $P_3$ on three vertices (where the internal vertices of these paths form the set $\cS_{1,1}$), while $G[\cL_{1,2} \cup \cS_{1,2}]$ consists of a disjoint union of paths $P_2$ on two vertices (where each path contains one vertex of $\cL_{1,2}$ and one vertex of $\cS_{1,2}$.) Therefore,
\[
\ell_{1,1} = 2s_{1,1} \hspace*{1cm} \mbox{and} \hspace*{1cm} \ell_{1,2} = s_{1,2},
\]
implying that
\begin{equation}
s_1 \ge \frac{1}{2} \ell_{1,1} + \ell_{1,2}.
 \label{Eq2}
\end{equation}

Counting the edge joining the large vertices to the ends of $2$-paths of order~$2$, this sum is exactly $2s_2$ and at least $3\ell_0 + 2\ell_1 + \ell_2$. Thus,
\begin{equation}
2s_2 \ge 3\ell_0 + 2\ell_{1,1} + 2\ell_{1,2} + \ell_2.
 \label{Eq3}
\end{equation}
By (\ref{Eq1}), (\ref{Eq2}) and (\ref{Eq3}) we therefore have that
\[
\begin{array}{lcl}
n &  = & \ell + s_1 + 2s_2  \1 \\
    & \ge & 4\ell_0 + \frac{7}{2}\ell_{1,1} + 4\ell_{1,2} + 2\ell_2 \2 \\
    & \ge & 4\ell_0 + 3\ell_{1,1} + 4\ell_{1,2} + 2\ell_2 \2 \\
    & = & 2\ell + 2\ell_0 + \ell_{1,1} + 2\ell_{1,2},
\end{array}
\]
or, equivalently,
\begin{equation}
\frac{n}{2} \ge \ell + \ell_0 + \frac{1}{2} \ell_{1,1} + \ell_{1,2}.
 \label{Eq4}
\end{equation}

For each vertex $v \in \cL_0$, let $v'$ be an arbitrary neighbor of $v$. Let $D' = \bigcup \{v'\}$, where the union is taken over all vertices $v \in \cL_0$. Then, $|D'| = \ell_0$. We now consider the set
\[
D = \cL \cup D' \cup \cS_{1,1} \cup \cS_{1,2}.
\]

Every vertex in $G$ is totally dominated by the set $D$, except possibly for vertices in $\cL_2$ which are disjunctively dominated by $D$ (since each vertex in $\cL_2$ is at distance~$2$ from at least two vertices in $\cL$). The set $D$ is therefore a DTD-set of $G$, implying by (\ref{Eq4}) that
\[
\begin{array}{lcl}
\dtd(G) & \le & |D| \1 \\
& = & \ell + \ell_0 + s_{1,1} + s_{1,2} \2 \\
& = & \ell + \ell_0 + \frac{1}{2} \ell_{1,1} + \ell_{1,2} \1 \\
& \le & \frac{n}{2}.
\end{array}
\]

\begin{lem}
\label{cLcL0:lem}
Every $2$-path has order~$2$; that is, $\cL = \cL_0$.
\end{lem}
\textbf{Proof of Lemma~\ref{cLcL0:lem}.}
If $\cL_{1,2} \ne \emptyset$, then removing from $D$ an arbitrary vertex that belongs to the set $\cL_{1,2}$ produces a DTD-set of cardinality~$|D| - 1$, implying that $\dtd(G) \le (n-2)/2$, a contradiction. Hence, $\cL_{1,2} = \emptyset$. This in turn implies that $\cS_{1,2} = \emptyset$. Thus, $\ell_{1,2} = s_{1,2} = 0$ and Inequality~(\ref{Eq2}) simplifies to $s_1 \ge \frac{1}{2}\ell_{1,1}$.

If $\cL_2 \ne \emptyset$, then there are at least two vertices in $\cS_1$ that do not belong to $\cS_{1,1}$, and so $s_1 \ge 2 + \frac{1}{2}\ell_{1,1}$. But then Inequality~(\ref{Eq4}) can be strengthened to $\frac{n}{2} \ge \ell + \ell_0 + \frac{1}{2}\ell_{1,1} + 1$, implying that $\dtd(G) \le |D| \le n/2 - 1$, a contradiction. Hence, $\cL_2 = \emptyset$. Thus, $\ell_2 = 0$ and $s_1 = \frac{1}{2}\ell_{1,1}$.

Suppose that $\cL_1 \ne \emptyset$. Then, by our earlier observations, $\cL_1 = \cL_{1,1}$. As observed earlier, $\ell_{1,1}$ is even and the subgraph $G[\cL_{1,1} \cup \cS_{1,1}]$ induced by $\cL_{1,1} \cup \cS_{1,1}$ consists of a disjoint union of paths $P_3$ with the internal vertices of these paths forming the set $\cS_{1,1}$. Removing from $D$ an arbitrary vertex that belongs to the set $\cL_{1,1}$ produces a DTD-set of cardinality~$|D| - 1$, implying that $\dtd(G) \le (n-2)/2$, a contradiction.  Hence, $\cL_1 = \emptyset$, implying that $\cL = \cL_0$. This completes the proof of Lemma~\ref{cLcL0:lem}.~\qed

\medskip
By Lemma~\ref{cLcL0:lem}, every $2$-path in $G$ has order~$2$, implying that $\cL = \cL_0$ and $D = \cL_0 \cup D'$.

\begin{lem}
\label{no6cycle:lem}
There is no $6$-cycle in $G$.
\end{lem}
\textbf{Proof of Lemma~\ref{no6cycle:lem}.}  Suppose to the contrary that $G$ contains a $6$-cycle $C \colon v_1v_2v_3v_4v_5v_6v_1$. Renaming vertices if necessary, we may assume that $v_1$ and $v_4$ are large vertices. Recall that for each vertex $v \in \cL_0$, we let $v'$ be an arbitrary neighbor of $v$ and the set $D' = \bigcup \{v'\}$, where the union is taken over all vertices $v \in \cL_0$. We now remove from $D$ the vertex $v_1$ and the two vertices $v_1'$ and $v_4'$ (that belong to $D'$) and we replace them with the two neighbors of $v_1$ on $C$ (namely, the vertices $v_2$ and $v_6$). The resulting set is a DTD-set of cardinality~$|D| - 1$, implying that $\dtd(G) \le (n-2)/2$, a contradiction.~\smallqed

\medskip
By Lemma~\ref{no6cycle:lem}, $G$ contains no $6$-cycle, implying that the girth of $G$ is at least~$7$. This in turn implies that $|\cL| \ge 4$. We now choose an arbitrary vertex $v \in \cL$. Let $v_1,v_2, \ldots, v_q$ be the vertices in $\cL$ at distance~$3$ from $v$. Since $G$ contains no $6$-cycle, we note that $q = d_G(v)$. For $i = 1,2,\ldots,q$, let $v_ia_ib_iv$ be a path in $G$ (and so, $a_ib_i$ is a $2$-path in $G$ of order~$2$ whose one end $a_i$ is adjacent to $v_i$ and whose other end $b_i$ is adjacent to $v$). If we now choose the vertex $v_i'$ to be the vertex $a_i$ for each $i = 1,2,\ldots,q$, then the vertex $v$ can be removed from $D$ to produce a DTD-set of cardinality~$|D| - 1$, implying that $\dtd(G) \le (n-2)/2$, a contradiction. This completes the proof of Theorem~\ref{em:thm}.~\qed

\subsection{Proof of Theorem~\ref{n2:thm}}
\label{n2:subsec}

Let $G$ be a connected graph of order $n\ge 8$ with $\delta(G)\ge 2$. Since $\dtd(G)$ cannot increase if edges are added, it follows from Theorem~\ref{em:thm} and Observation~\ref{f1:ob} that $\dtd(G)\le n/2$. Further, suppose $\dtd(G)=n/2$. We produce a $\frac{1}{2}$-minimal graph $G'$ from $G$ by removing edges if necessary so that $G'$ satisfies $\dtd(G')=n/2$. By Theorem~\ref{em:thm} and Observation~\ref{f1:ob}, $G'\in\{B_{3}, D_{b}(4,4), D_{b}(3,4,1), D_{b}(3,3,2), C_{8}, C_{12}\}$. In all cases it can be readily checked that $G=G'$ or $G'\in\cF$ where $\cF$ is the family of graphs shown in Fig.~\ref{f:cF}.

\subsection{Proof of Theorem~\ref{nminus:thm}}
\label{nminus:subsec}

Let $G$ be a connected graph of order $n\ge 18$ with $\delta(G)\ge 2$. Since $\dtd(G)$ cannot increase if edges are added, it follows from Theorem~\ref{em:thm} and Observation~\ref{f1:ob} that $\dtd(G)\le (n-1)/2$. Further, suppose $n\ge 18$ and $\dtd(G)=(n-1)/2$. We produce a $\frac{1}{2}$-minimal graph $G'$ from $G$ by removing edges if necessary so that $G'$ satisfies $\dtd(G')=(n-1)/2$. Since $n\ge 18$, $G'\in \cG$ by Observation~\ref{f1:ob}. It can readily be checked that $G=G'$ or $G\in\cH$ where $\cH$ is the family of graphs constructed in Section~\ref{s3:sec}.

\newpage
\begin{center}
\underline{\textbf{APPENDIX}}
\end{center}

This appendix contains proofs of some of selected preliminary results from Section~\ref{pre:subsec}. We introduce the following notation for the distance between sets. Let $G=(V,E)$ be a graph. For vertex sets $X,Y\subseteq V$, the distance between the sets $X$ and $Y$, denoted $d_G(X,Y)$ or simply by $d(X,Y)$ if $G$ is clear from context, is the minimum distance $d(x,y)$ taken over all possible pairs of vertices $x\in X$, and $y\in Y$. We begin by establishing the value of $\dtd(G)$ when $G$ is a path.

\begin{prop}
\label{path:prop:ap}
If $G=P_{n}$ and $n\ge 2$ then, $\dtd(G)=\lceil 2(n+1)/5\rceil+1$ if $n\equiv 1 (\mod \, 5)$, and $\dtd(G)=\lceil 2(n+1)/5\rceil$ otherwise.
\end{prop}
\textbf{Proof of Proposition~\ref{path:prop:ap}.} We proceed by induction on $n \ge 2$. The result is easily verified for $n \le 9$. Suppose that $n \ge 10$ and the result is true for all paths of order less than~$n$. Let $G=P_{n}$ be a path given by $v_1v_{2} \ldots v_n$. We first establish upper bounds on $\dtd(G)$. Let
\[
S = \bigcup_{i = 0}^{ \lfloor n/5 \rfloor - 1} \{v_{5i+2},v_{5i+3}\}.
\]

If $n \equiv 0 \, (\mod \, 5)$, let $D = S\cup\{v_{n-1}\}$.
If $n \equiv i \, (\mod \, 5)$, where $i \in \{1,2,3,4\}$, let $D = S \cup \{v_{n-2},v_{n-1}\}$.
In all cases, the set $D$ is a DTD-set of $G$. Further if $n\equiv 1\,(\mod \, 5)$, then $|D| = \lceil 2(n+1)/5\rceil+1 $, while if $n \not\equiv 1 \,(\mod \, 5)$, then $|D| = \lceil 2(n+1)/5 \rceil$. Hence, $\dtd(G) \le \lceil 2(n+1)/5\rceil+1 $ if $n \equiv 1 \, (\mod \, 5)$ and $\dtd(G) \le \lceil 2(n+1)/5 \rceil$ if $n \not\equiv 1 \,(\mod \, 5)$.

To prove the reverse inequality, let $T$ be a $\dtd(G)$-set. Since $v_{1}$ is a leaf vertex in $G$, we may suppose $v_{2}\in T$. Further, since $T$ cannot disjunctively dominate $v_{2}$ and $n\ge 10$ we may suppose $v_{3}\in T$. Hence, $\{v_{2}, v_{3}\}\subset T$. We show next that we can choose $T$ so that $T \cap \{v_4,v_5,v_6\} = \emptyset$.

Suppose that $v_4 \in T$. If $v_5 \in T$, let $i$ be the largest integer such that $i \ge 5$ and $v_i \in T$, and replace $v_4$ in $T$ with the vertex $v_{i+1}$. Suppose that $v_5 \notin T$. If $v_6 \in T$, replace $v_4$ in $T$ with $v_5$. If $v_6 \notin T$ and $v_7 \in T$, replace $v_4$ in $T$ with $v_6$. If $v_6 \notin T$ and $v_7 \notin T$, then $v_8 \in T$ and replace $v_4$ in $T$ with $v_7$. In all the above cases, we can choose $T$ so that $v_4 \notin T$.

Suppose that $v_5 \in T$. If $v_6 \notin T$, then $v_7 \in T$ and we can replace $v_5$ in $T$ with the vertex $v_6$. If $v_6 \in T$, let $i$ be the largest integer such that $i \ge 6$ and $v_i \in T$, and replace $v_5$ in $T$ with the vertex $v_{i+1}$. In both cases we can choose $T$ so that $v_5 \notin T$.

Suppose that $v_6 \in T$. Then, $v_7 \in T$. Let $i$ be the largest integer such that $i \ge 7$ and $v_i \in T$, and replace $v_6$ in $T$ with the vertex $v_{i+1}$. Hence we can choose $T$ so that $v_6 \notin T$. Therefore, $T \cap \{v_4,v_5,v_6\} = \emptyset$, implying that $\{v_7,v_8\} \subset T$. Let $T'  = T \setminus \{v_1,v_2\}$, and note that $|T'| = |T| - 2$.

We now let $G'$ be obtained from $G$ by deleting the vertices $v_i$, $1 \le i \le 5$. Then, $G' = C_{n'}$, where $n' = n - 5 \ge 5$. Since $T$ is a DTD-set of $G$, the set $T'$ is a DTD-set of $G'$. Hence, $\dtd(G') \le |T'| = |T| - 2$. Applying the inductive hypothesis to $G'$, we have that $\dtd(G') = \lceil 2(n'+1)/5 \rceil+1=\lceil 2(n+1)/5\rceil-1$ if $n \equiv 1 \, (\mod \, 5)$ and $\dtd(G') = \lceil 2(n'+1)/5 \rceil = \lceil 2(n+1)/5 \rceil - 2$ if $n \not\equiv 1 \, (\mod \, 5)$. This implies that $\dtd(G) = |T| \ge
\lceil 2(n+1)/5\rceil+1$ if $n \equiv 1 \, (\mod \, 5)$ and $\dtd(G) = |T| \ge \lceil 2(n+1)/5 \rceil$ if $n \not\equiv 1 \, (\mod \, 5)$. The desired bounds now follow as a consequence of the upper bounds on $\dtd(G)$ established earlier.~\qed

In order to establish upper bounds on the value of $\dtd(G)$ for a daisy we require the following proposition.

\begin{prop}
\label{daisy:prop:ap}
Let $G$ be a daisy with $k$~petals. If $v$ is the vertex of degree~$2k$ in $G$ and no $\dtd(G)$-set contains $v$, then $v$ is totally dominated in $G$.
\end{prop}
\textbf{Proof of Proposition~\ref{daisy:prop:ap}.} Let $G$ be a daisy of order~$n$ with $k$ petals. Let $v$ be the vertex of degree~$2k$ in $G$ and suppose that $v$ is not contained in any $\dtd(G)$-set. We proceed by induction on the number of petals $k\ge 2$. We establish the base case with the following claim.

\begin{unnumbered}
{Claim N} $k\ge 3$.
\end{unnumbered}
\textbf{Proof of Claim~N.} Suppose $k=2$. Then, $G=D(n_{1}+1, n_{2}+1)$ is a daisy with two petals and $v$ is the vertex of degree~$4$ in $G$. Let $F_{1}$ and $F_{2}$ denote the two cycles passing through $v$, where $F_{i}\cong C_{n_{i}+1}$ for $i=1,2$. Let $F_{1}=vu_{1}\ldots u_{n_{1}}v$ and let $F_{2}=vv_{1}\ldots v_{n_{2}}v$. We note that $n=n_{1}+n_{2}+1\ge 5$. Assume $v$ is not contained in any $\dtd(G)$-set and that no $\dtd(G)$-set totally dominates $v$ in $G$.

Let $S$ be a $\dtd(G)$-set. Then, $v\notin S$ and $S$ disjunctively dominates $v$ in $G$. Since $v\notin S$ $|S\cap V(F_{i})|\le \dtd(F_{i})-1$ for $i=1,2$, otherwise restricting $S$ to one of the cycles $F_{i}$ gives a DTD-set in that cycle whereby (rearranging vertices) we may choose $v\in S$, a contradiction.
Further, $N[v]\cap S=\emptyset$ implying that $n_{1},n_{2}\ge 4$, otherwise $S$ is either not a $\dtd(G)$-set or $N(v)\cap S\ne \emptyset$ both of which are contradictions. Hence, $n\ge 9$. Observe that each vertex in $N(v)$ has exactly one vertex at distance 2 from it in $V(G)\setminus N(v)$. Hence, since $N(v)\cap S=\emptyset$ every vertex in $N(v)$ is totally dominated by $S$. Furthermore, by the same argument, if $x\in N(v)$, then the neighbor of $x$ in $N(x)\setminus\{v\}$ is contained in $S$. Suppose $n_{1}=5$. Then, $\{u_2,u_3\}\subset S$ and $|S\cap V(F_1)|=\dtd(F_{1})$, a contradiction. Suppose $n_{1}=6$. Then, $\{u_{2},u_{3},u_{4}\}\subset S$ and $|S\cap V(F_1)|=\dtd(F_{1})$, a contradiction. By symmetry, $n_{2}\notin\{5,6\}$. Hence, for $i=1,2$, suppose $n_{i}\ge 7$. Thus, $\{v_{2},v_{3},v_{n-2},v_{n-3}\}\subset S$ and $\{u_{2},u_{3},u_{n-2},u_{n-3}\}\subset S$. It follows that $|S\cap V(F_{i})|\ge 4$. Then, by Proposition~\ref{cycle:ob} and the fact that $|S\cap V(F_{i})|\le \dtd(F_{i})-1$, $n_{i}\ge 10$ for $i=1,2$.

We now show that there is a $\dtd(G)$-set, $S^{*}$, which contains the vertices $u_{2+5j}$ and $u_{3+5j}$ for $j\in\{0,\ldots, \lfloor (n_{1}+1)/5 \rfloor-1\}$. We proceed by induction on $j\ge 0$. Suppose that for $j^*\le j'=j-1<j$, the set $S$ contains the required vertices, and so  $\{u_{2+5j'}, u_{3+5j'}\}\subset S$ but the vertices $u_{2+5j}$ and $u_{3+5j}$ are not necessarily contained in $S$. If $j=0$, then the result is true, since $\{u_{2},u_{3}\}\subset S$. This establishes the base case. Note that by our assumptions $\{u_{n_{1}-1}, u_{n_{1}-2}\}\subset S$, regardless. Let $T=\{u_{4+5j'}, u_{5+5j'}, u_{6+5j'}\}$.

Suppose $n_{1}\equiv 0\,(\mod\,5)$. Then, $u_{2+5j}=u_{n_{1}-3}$ and $u_{3+5j}=u_{n_{1}-2}$. Thus, $u_{3+5j}\in S$. Suppose $S\cap T=\emptyset$. Then, since $S$ is a DTD-set of $G$, $u_{2+5j}=u_{n_{1}-3}\in S$, and $S^{*}=S$ contains the required vertices. Hence, $S\cap T\ne \emptyset$. Suppose $|S\cap T|\ge 2$, then $S^{*}=(S\setminus T)\cup\{u_{2+5j}\}$ is a DTD-set of $G$ of smaller cardinality than $S$, a contradiction. Hence $|S\cap T|=1$. In that case, $S^{*}=(S\setminus T)\cup \{u_{2+5j}\}$ is a $\dtd(G)$-set which contains the required vertices.

Suppose $n_{1}\equiv 1\,(\mod\,5)$. Then, $u_{2+5j}=u_{n_{1}-4}$ and $u_{3+5j}=u_{n_{1}-3}$. Suppose $S\cap T=\emptyset$. Then, since $S$ is a DTD-set of $G$, $\{u_{2+5j}, u_{3+5j}\}\subset S$, and $S^{*}=S$ contains the required vertices. Hence, $S\cap T\ne \emptyset$. Suppose $|S\cap T|=3$. Then, $S^{*}=(S\setminus T)\cup\{u_{2+5j}, u_{3+5j}\}$, is a DTD-set of $G$ of smaller cardinality than $S$, a contradiction. Suppose $|S\cap T|=2$. Then, $S^{*}=(S\setminus T)\cup\{u_{2+5j}, u_{3+5j}\}$, is a $\dtd(G)$-set which contains the required vertices. Hence, $|S\cap T|=1$. Assume $u_{4+5j'}\in S$. Then, $u_{3+5j}\in S$, otherwise $S$ is not a DTD-set of $G$, a contradiction. Then, let $S^*=(S\setminus\{u_{4+5j'}\})\cup\{u_{2+5j}\}$. Assume $u_{5+5j'}\in S$. Then, $u_{2+5j}\in S$, otherwise $S$ is not a DTD-set of $G$, a contradiction. Then, let $S^*=(S\setminus\{u_{5+5j'}\})\cup\{u_{3+5j}\}$. Assume $u_{6+5j'}\in S$. Then, $u_{2+5j}\in S$, otherwise $S$ is not a DTD-set of $G$, a contradiction. Then, let $S^*=(S\setminus\{u_{5+5j'}\})\cup\{u_{3+5j}\}$. In each case $S^*$ is a $\dtd(G)$-set which contains the required vertices.

Suppose $n_{1}\equiv 2\,(\mod\,5)$. Then, $u_{2+5j}=u_{n_{1}-5}$ and $u_{3+5j}=u_{n_{1}-4}$. Suppose $S\cap T=\emptyset$. Then, since $S$ is a DTD-set of $G$, $\{u_{2+5j}, u_{3+5j}\}\subset S$, and $S^{*}=S$ contains the required vertices. Hence, $S\cap T\ne \emptyset$. Suppose $|S\cap T|=3$, then $S^{*}=(S\setminus T)\cup\{u_{2+5j}, u_{3+5j}\}$, is a DTD-set of $G$ of smaller cardinality than $S$, a contradiction. Suppose $|S\cap T|=2$, then $S^{*}=(S\setminus T)\cup\{u_{2+5j}, u_{3+5j}\}$, is a $\dtd(G)$-set which contains the required vertices. Hence, $|S\cap T|=1$. Let $R_{1}=\{u_{2+5j}, u_{3+5j}, u_{4+5j}\}$. We must have $|S\cap R_{1}|= 1$, otherwise $S$ is not a $\dtd(G)$-set, a contradiction. In each possible case $S^{*}=(S\setminus(T\cup R_{1}))\cup\{u_{2+5j}, u_{3+5j}\}$ is a $\dtd(G)$-set which contains the required vertices.

Suppose $n_{1}\equiv 3\,(\mod\,5)$. Then, $u_{2+5j}=u_{n_{1}-6}$ and $u_{3+5j}=u_{n_{1}-5}$. Suppose $S\cap T=\emptyset$. Then, since $S$ is a DTD-set of $G$, $\{u_{2+5j}, u_{3+5j}\}\subset S$, and $S^{*}=S$ contains the required vertices. Hence, $S\cap T\ne \emptyset$. Suppose $|S\cap T|=3$, then $S^{*}=(S\setminus T)\cup\{u_{2+5j}, u_{3+5j}\}$, is a DTD-set of $G$ of smaller cardinality than $S$, a contradiction. Suppose $|S\cap T|=2$, then $S^{*}=(S\setminus T)\cup\{u_{2+5j}, u_{3+5j}\}$, is a $\dtd(G)$-set which contains the required vertices. Hence, $|S\cap T|=1$. Let $R_{2}=\{u_{2+5j}, u_{3+5j}, u_{4+5j}\}=\{u_{n_{1}-6}, u_{n_{1}-5}, u_{n_{1}-4}\}$. We must have $|S\cap R_{2}|= 1$, otherwise $S$ is not a $\dtd(G)$-set, a contradiction. In each possible case $S^{*}=(S\setminus(T\cup R_{2}))\cup\{u_{2+5j}, u_{3+5j}\}$ is a $\dtd(G)$-set which contains the required vertices.

Suppose $n_{1}\equiv 4\,(\mod\,5)$. Then, $u_{2+5j}=u_{n_{1}-2}$ and $u_{3+5j}=u_{n_{1}-1}$. Thus, $S^{*}=S$ by our previous assumptions and the inductive hypothesis, and so we are done.

We now use this result to obtain a contradiction. Let $X=S^{*}\cap V(F_{1})$ and consider $|X|$. If $n_{1}\equiv 0\,(\mod\,5)$, then $|X|=2\lfloor (n_{1}+1)/5 \rfloor+1$. If $n_{1}\equiv 1$, $2$ or $3\,(\mod\,5)$, then $|X|=2\lfloor (n_{1}+1)/5 \rfloor+2$. If $n_{1}\equiv 4\,(\mod\,5)$, then $|X|=2\lfloor (n_{1}+1)/5 \rfloor=2(n_{1}+1)/5$. We may rewrite $|X|$ as : $|X|=2(n_{1}+1)/5$ if $(n_{1}+1)\equiv 0\,(\mod\,5)$ and $|X|=\lceil 2(n_{1}+2)/5\rceil$ otherwise. Then, by Proposition~\ref{cycle:ob}, in each case, $|X|= \dtd(F_{1})> \dtd(F_{1})-1$, a contradiction.~\qed

\newpage
We now return to the proof of Proposition~\ref{daisy:prop:ap}. Suppose $v$ is not contained in any $\dtd(G)$ and no $\dtd(G)$-set totally dominates $v$. Let $S$ be a $\dtd(G)$-set. Then, $S$ disjunctively dominates $v$ and $v\notin S$. Assume that the result is true for all daisies with $k'\le k-1<k$ petals. By Claim N, $k\ge 3$. Let $G'$ be a daisy obtained from $G$ by removing every small vertex (vertex of degree~$2$) of a single petal so that $G'$ has $2\le k'<k$ petals. By the inductive hypothesis, there is a $\dtd(G')$-set, $S'$, which totally dominates $v$. Let $X=S\cap V(G')$ and let $Y=S\cap (V(G) \setminus V(G'))$. Since $v$ is disjunctively dominated by $S$, $d(X,Y)\ge 4$, and so $X$ and $Y$ DT-dominate $G'$ and $G-V(G')$ respectively. Furthermore, we must have $|S'|=|X|$. Then, $(S\setminus X)\cup S'$ is a $\dtd(G')$-set which totally dominates $v$.~\qed

\begin{unnumbered}
{Proposition~\ref{d:ob}}
If $G$ is a daisy of order $n$, then $\dtd(G)\le (n-1)/2$. Furthermore, $\dtd(G)=(n-1)/2$ if and only if $G \in  \cD$.
\end{unnumbered}
\textbf{Proof of Proposition~\ref{d:ob}.} Let $G$ be a daisy of order $n$ with $k$ petals. We first establish the following claim.

\begin{unnumbered}
{Claim O} $k\ge 3$.
\end{unnumbered}
\textbf{Proof of Claim O.} We show that the desired result is true if $k=2$.

Suppose $k=2$. The sufficiency is straightforward to check. In order to prove the necessity, we proceed by induction on $n$. If $n=5$, then $G=D(3,3)\in \cD$ and $\dtd(G)=2=(n-1)/2$. If $n=6$, then $G=D(3,4)$ and $\gamma_{t}^{d}(G)=2 < (n-1)/2$. If $n=7$, then $G=D(3,5)$ and $\dtd(G)=2<(n-1)/2$  or $G=D(4,4)\in \cD$ and $\dtd(G)=3=(n-1)/2$. Hence, $n\ge 8$ otherwise we are done.

Let $G=D(n_{1}+1, n_{2}+1)$ so $n=n_{1}+n_{2}+1$. Let $v$ denote the vertex of degree~4 in $G$. Let $F_{1}$ and $F_{2}$ denote the two cycles passing through $v$ where $F_{i}\cong C_{n_{i}+1}$ for $i=1,2$. Let $F_{1}$ be the cycle $vv_{1} \ldots v_{n_{1}}v$ and let $F_{2}$ be the cycle $vu_{1}\ldots u_{n_{2}}v$. Let $S_{1}=\{v_{i} \mid i\equiv 0$ or $1\,(\mod\, 5)\}$ and let $S_{2}$ be a $\dtd(F_{2})$-set that contains $v$ and $v_{1}$. By Proposition~\ref{cycle:ob}, $|S_{2}|=2(n_{2}+1)/5$ if $n_{2}\equiv 4\,(\mod\,5)$ and $|S_{2}|=\lceil 2(n_{2}+2)/5\rceil$ otherwise. Let $D_{1}$ be a $\dtd(F_{1})$-set which contains $v$ and $u_{1}$ and let $D_{2}=\{v_{i} \mid i\equiv 0$ or $4\,(\mod\,5)\}$. By Proposition~\ref{cycle:ob}, $|D_{1}|=2(n_{1}+1)/5$ if $n_{1}\equiv 4\,(\mod\,5)$ and $|D_{1}|=\lceil 2(n_{1}+2)/5\rceil$ otherwise.

Suppose $n_{1}\equiv 4\,(\mod\,5)$. Then, $|S_{1}|=(2n_{1}-3)/5$. If $n_{2}\equiv 0$, $2$, $3$ or $4\,(\mod\,5)$, then $\dtd(G)\le |S_{1}\cup S_{2}|\le (2n+1)/5<(n-1)/2$. On the other hand, if $n_{2}\equiv 1(\,\mod\,5)$, then $|D_{1}|=2(n_{1}+1)/5$ and $|D_{2}|=2(n_{2}-1)/5$ implying that $\dtd(G)\le |D_{1}\cup D_{2}|\le (2n-2)/5<(n-1)/2$. Hence, $n_{i}\not\equiv 4\,(\mod\,5)$ for $i=1,2$.

Suppose $n_{1}\equiv 0\,(\mod\,5)$. Then, $|S_{1}|=2n_{1}/5$. If $n_{2}\equiv 0$ or $3 \, (\mod\,5)$, then $T=S_{1}\cup (S_{2}\setminus\{v\})$ is a DTD-set of $G$ and $\dtd(G)\le |T|\le 2(n-1)/5<(n-1)/2$. If $n_{2}\equiv 1\,\textnormal{or}\,2\,(\mod\,5)$, then $|D_{1}|=(2n_{1}+5)/5$ and $|D_{2}|\le 2(n_{2}-1)/{5}$. Since $D_{1}\cup D_{2}$ is a DTD-set of $G$ it follows that $\gamma_{t}^{d}(G)\le (2n+1)/5 <(n-1)/2$. Hence, $n_{i}\not\equiv 0\,(\mod\,5)$ for $i=1,2$ and, further, $n\ge 9$.

Suppose $n_{1}\equiv 3\,(\mod\,5)$. Then, $|D_{1}|=(2n_{1}+2)/5$. If $n_{2}\equiv 1,2$ or $3\,(\mod\,5)$, then $|D_{2}|\le 2(n_{2}-1)/5$ and $\dtd(G)\le |D_{1}\cup D_{2}|\le (2n-1)/5< (n-1)/2$. Hence, $n_{i}\not\equiv 3\,(\mod\,5)$ for $i=1,2$.

Suppose $n_{1}\equiv 2\,(\mod\, 5)$. Then, $|D_{1}|=(2n_{1}+6)/{5}$. If $n_{2}\equiv 2\,(\mod\, 5)$, then $|D_{2}|=(2n_{2}-4)/5$ and $\dtd(G)\le |D_{1}\cup D_{2}|\le 2n/5<(n-1)/2$. If $n_{2}\equiv 1\,(\mod\, 5)$, then $|D_{2}|=(2n_{2}-2)/5$ and $\dtd(G)\le |D_{1}\cup D_{2}|\le (2n+2)/5\le (n-1)/2$ with equality if and only if $G=D(3,7)\in \cD$. Hence, $n_{i}\not\equiv 2\,(\mod\, 5)$ for $i=1,2$.

Finally, suppose $n_{i}\equiv 1\,(\mod\, 5)$ for $i=1,2$. Then, $n\ge 13$ since $D(7,7)$ is the graph of smallest order with $n_{i}\equiv 1\,(\mod\, 5)$ for $i=1,2$. Let  $R_{1}=\{u_{i} \mid i\equiv 2$ or $3\,(\mod\, 5)\}$ and $R_{2}=\{v_{i} \mid i\equiv 2$ or $3\,(\mod\, 5)\}$. Then $|R_{1}|=2(n_{1}-1)/5$ and $|R_{2}|=2(n_{2}-1)/5$. The set $|R_{1}\cup R_{2}|\cup\{v\}$ is DTD-set of $G$. Hence, $\dtd(G)\le 2(n_{1}-1)/5+2(n_{2}-1)/5+1=(2n-1)/5<(n-1)/2$. Thus, $k\ge 3$.~\qed

We now return to the proof of Proposition~\ref{d:ob}. We proceed by induction on $n$. Assume that for all daisies of order less than $n$ the result is true. By Claim~O, $k\ge 3$, and so $n\ge 7$. If $n=7$, then $G=D(3,3,3)$ with $\dtd(G)=2<(n-1)/2$, as required. Hence, $n\ge 8$.

Let $v$ denote the vertex of degree~$2k$ in $G$, and let $F_{1}$, $F_{2},\ldots, F_{k}$ denote the $k$ cycles passing through $v$, where $F_{i}\cong C_{n_{i}+1}$ for $i=1,2\ldots k$. Thus, $n=1+\sum_{i=1}^{k}n_{i}$.   Let $F_{1}$ be the cycle $vv_{1}\ldots v_{n_{1}}v$.

Let $G'=D(n_{2},\ldots, n_{k}).$ Then, $G'$ is a daisy of order $n'=n-n_{1}$ with $k-1\ge 2$ petals. Applying the inductive hypothesis to $G'$, we have $\dtd(G')\leq (n'-1)/2$ with equality if and only if $G'\in\cD$. Let $S'$ be a $\dtd(G')$-set.

Suppose $v\in S'$. If $n_{1}=2$, then $S'$ is a DTD-set of $G$, and so $\dtd(G)\le |S'|<(n-1)/2$. Hence, $n_{1}\ge 3$. We now extend $S'$ to a DTD-set of $G$ as follows. If $n_{1}\equiv 0$ or $4\,(\mod\,5)$, then let $S_{1}=\{v_i \mid i\equiv 0$ or $1\,(\mod\,5)\}$. If $n_{1}\equiv 1$ or $2\,(\mod\,5)$, then let $S_{1}=\{v_{i} \mid i\equiv 3$ or $4\,(\mod\,5)\}$. If $n_{1}\equiv 3\,(\mod\,5)$, then let $S_{1}=\{v_{i} \mid i\equiv 0$ or $4\,(\mod\,5)\}$. Then, $|S_{1}|\le 2n_{1}/5$ and $\dtd(G)\le |S'\cup S_{1}|\le (n'-1)/2+2n_{1}/5\le (n-1)/2-n_{1}/10\le (n-1)/2$. Further, if $\dtd(G)=(n-1)/2$ we must have equality throughout this inequality chain. In particular, $\dtd(G)=(n'-1)/2$ and $n_{1}=0$. However, this contradicts Claim~O, implying $\dtd(G)<(n-1)/2$ for all daisies with $k\ge 3$ and $v\in S'$.

Suppose $v\notin S'$. We may suppose that no $\dtd(G')$-set contains $v$, for otherwise, we are done. For $2\le i\le k$, in each cycle $F_{i}$ we require that $|S'\cap V(F_{i})|\le \dtd(F_{i})-1$, otherwise the restriction of $S'$ to the vertices of some cycle in $G'$ is a DTD-set of that cycle and we may then choose $v\in S'$, a contradiction. Further, since each graph $H\in \cD$ has a $\dtd(H)$-set which contains $v$, if $v\notin S'$ then $G'\notin \cD$. Hence, $G'\notin \cD$ and so $\dtd(G')=|S'|\le (n'-2)/2$.

Now, by Proposition~\ref{daisy:prop:ap}, $v$ is totally dominated in $G'$. Let $S_{1}$ be defined as before. Suppose $n_{1}\equiv 1\,(\mod\,5)$. Then, $S'\cup S_{1}$ is a DTD-set of $G$ and so $\dtd(G)\le (n'-2)/2+2(n_{1}-1)/5=(n-2)/2-(n_{1}+2)/10<(n-1)/2$. Suppose $n_{1}\equiv 3$ or $4\,(\mod\,5)$. Then, $S_{1}\cup S'\cup \{v\}$ is a DTD-set of $G$ and so $\dtd(G)\le 2(n_{1}-2)/5+1+(n'-2)/2=(n-2)/2+(2-n_{1})/10<(n-1)/2$, where the last inequality follows from the fact that $n_{1}\ge 2$. Since $F_{1}$ was chosen arbitrarily, it follows that $G'$ contains no cycle such that $n_{i}\equiv 1$, $2$ or $3\,(\mod\,5)$. By Proposition~\ref{cycle:ob}, $\dtd(F_{i})= 2n_{i}/5+1$ if $n_{i}\equiv 0\,(\mod\,5)$ and $\dtd(F_{i})= 2(n_{i}+1)/5$ if $n_{i}\equiv 4\,(\mod\,5)$. Therefore $|S'|\le \sum_{i=2}^{k}(\dtd(F_{i})-1)\le \sum_{i=2}^{k}(2n_{i}/5)=2(n-n_{1}-1)/5$. Let $S_{1}^*$ be a $\dtd(F_{1})$-set which includes $v$ if $n_{1}\equiv 4\,(\mod\,5)$ and let $S_{1}^{*}=S_{1}$ if $n_{1} \equiv 0 \,(\mod\,5)$. The set $S_{1}^*\cup S'$ is a DTD-set of $G$ and so $\dtd(G)\le 2(n_{1}+1)/5+2(n-n_{1}-1)/5=2n/5<(n-1)/2$, as required.~\qed

In order to prove our next result we require the value of the minimum disjunctive total domination number of a key. Let $G=L_{r,s}$ be the key formed by joining a vertex from a cycle on $r$ vertices to a leaf vertex of a path on $s$ vertices. Let $n=r+s$ be the order of $G$. Then, $\dtd(G)$ is given by Table~3 and the value of $\dtd(G)$ has been computed using Propositions~\ref{cycle:ob} and~\ref{path:prop:ap}. Further, let $y$ be the leaf vertex in $G$ and let $x$ be the neighbor of $y$. The vertex $x$ is contained in every $\dtd(G)$-set. The entries followed by a star $(*)$ denote that there is a $\dtd(G)$-set which contains $y$. The entries followed by a star $(*)$ and a dagger $(\dag)$ denote that there is a $\dtd(G)$-set, $S$, such that $y\in S$ and the set $S\setminus\{x\}$ DT-dominates the subgraph $G-y$ of $G$.  The entries followed by a dagger $(\dag)$ denote that there is a $\dtd(G)$-set, $S$, $y\notin S$, and the set $S\setminus\{x\}$ DT-dominates the subgraph $G-y$ of $G$. The double daggered entries denote that there is a $\dtd(G-\{x,y\})$-set with $\dtd(G-\{x,y\})=\dtd(G)-2$.

{
\[
\begin{array}{lcccccc}
\hline
\hline
 & \multicolumn{5}{c}{s} \\ \cline{2-6}
 r & \equiv 0\,(\mod\, 5) & \equiv 1\,(\mod\, 5) & \equiv 2\,(\mod\, 5) & \equiv 3\,(\mod\, 5) & \equiv 4\,(\mod\, 5)\\ \cline{1-6}
 \equiv 0\,(\mod\,5)  & 2n/5^{*} & 2(n-1)/5 & \left\lceil 2n/5 \right\rceil^{\dag} & \left\lceil 2n/5\right\rceil^{*\dag\ddag} & \left\lceil 2n/5\right\rceil^{*\dag} \\
 \equiv 1\,(\mod\, 5) & \left\lceil 2(n+1)/5\right\rceil^{*\dag\phantom{\ddag}} & \left\lceil 2n/5 \right\rceil^{*\phantom{\dag}} & \left\lceil 2(n-1)/5\right\rceil^{\phantom{*}} & \left\lceil 2n/5\right\rceil^{\dag\phantom{*\ddag}} & \left\lceil 2(n+1)/5\right\rceil^{*\dag\ddag}\\
 \equiv 2\,(\mod\, 5) & \left\lceil 2(n+1)/5\right\rceil^{*\dag\ddag} & \left\lceil 2n/5 \right\rceil^{*\dag} & \left\lceil 2(n-1)/5\right\rceil^{*} & \left\lceil 2n/5\right\rceil^{\phantom{*\dag\ddag}} & \left\lceil 2(n+1)/5\right\rceil^{\dag\phantom{*\ddag}}\\
 \equiv 3\,(\mod\, 5) & \left\lceil 2(n+1)/5\right\rceil^{*\dag\phantom{\ddag}} & \left\lceil 2n/5 \right\rceil^{*\phantom{\dag}} & \left\lceil 2(n-1)/5\right\rceil^{\phantom{*}} & \left\lceil 2n/5\right\rceil^{\dag\phantom{*\ddag}} & \left\lceil 2(n+1)/5\right\rceil^{*\dag\ddag}\\
\equiv 4\,(\mod\, 5) & \left\lceil 2(n+1)/5\right\rceil^{*\phantom{\dag\ddag}} & \left\lceil 2n/5 \right\rceil^{*\phantom{\dag}} & \left\lceil 2(n-1)/5\right\rceil^{\phantom{*}} &   \left\lceil 2n/5\right\rceil^{\dag\phantom{*\ddag}} & \left\lceil 2(n+1)/5\right\rceil^{*\dag\ddag} \\
\hline
\hline
\end{array}
\]}
\begin{center}
{Table~3:} The Disjunctive Total Domination Number of $L_{r,s}$ where $n=r+s$ and $r\ne 4$.
\end{center}

We are now in a position to prove Proposition~\ref{db:ob}. Recall its statement.

\begin{unnumbered}{Proposition~\ref{db:ob}} If $G$ is a dumb-bell of order $n$, then $\dtd(G)\le n/2$. Furthermore, $\dtd(G)\ge (n-1)/2$ if and only if $G \in \cD_{b} \cup \cG_{b}$.
\end{unnumbered}
\textbf{Proof of Proposition~\ref{db:ob}.} Let $G=D_{b}(n_{1}, n_{2}, \ell)$ be a dumb-bell of order $n=n_{1}+n_{2}+\ell$. We proceed with the following claim.

\begin{unnumbered}
{Claim P} $\ell\ge 1$.
\end{unnumbered}
\textbf{Proof of Claim P.} We show that if $\ell=0$, then the desired result follows.

Suppose $\ell=0$. Then, $G=D_{b}(n_{1},n_{2})$ and $n=n_{1}+n_{2}$. If $G\in\cD_{b}\cup\cG_{b}$, then $G\in\{D_{b}(3,4), D_{b}(4,4), D_{b}(4,5)\}$ and it may be verified that $\dtd(G)\ge (n-1)/2$. In order to prove the necessity, we proceed by induction on the order $n\ge 6$ of $G$. Let $F_{1}\cong C_{n_{1}}$ and $F_{2}\cong C_{n_{2}}$ be the two cycles of $G$. Let $F_{1}=uu_{1}\ldots u_{n_{1}-1}u$ and let $F_{2}=vv_{1}\ldots v_{n_{2}-1}v$. Let $uv$ be the edge joining the cycles $F_{1}$ and $F_{2}$ in $G$. Let $S_{1}$ be a $\dtd(F_{1})$-set which contains $u$. Let $S_{2}$ be a $\dtd(F_2)$-set which contains $v$. By Proposition~\ref{cycle:ob}, for $i=1,2$, $|S_{i}|=2n_{i}/5$ if $n_{i}\equiv 0\,(\mod \, 5)$ and $|S_{i}|=\lceil 2(n_{i}+1)/5\rceil$ otherwise. Let $D_{1}=\{u_{i} \mid i\equiv 0$ or $1\,(\mod\,5)\}$ and $D_{2}=\{v_{i} \mid i\equiv 0\,$ or $1\,(\mod\,5)\}$. Let $R_{1}=\{u_{i} \mid i\equiv 3\,$ or $4\,(\mod\,5)\}$ and $R_{2}=\{v_{i} \mid i\equiv 3\,$ or $4\,(\mod\,5)\}$. Let $T_{1}=\{u_{i} \mid i\equiv 0\,$ or $4\,(\mod\,5)\}$ and $T_{2}=\{v_{i} \mid i\equiv 0$ or $4\,(\mod\,5)\}$. If $n=6$, then $G=D_{b}(3,3)$ and $\dtd(G)=2<(n-1)/2$. If $n=7$, then $G=D_{b}(3,4)$ and $\dtd(G)=3=(n-1)/2$. In both cases the result holds, hence $n\ge 8$.

Suppose $n_{1}=4$. Then $n=n_{2}+4\ge 8$. If $n_{2}\equiv 0\,(\mod\,5)$, then $S_{1}\cup S_{2}$ is a DTD-set of $G$ and $\dtd(G)\le |S_{1}\cup S_{2}|=2n_{2}/5+2=2(n+1)/5\le (n-1)/2$. Further, we have equality if and only if $G=D_{b}(4,5)$. If $n_{2}\equiv 1\,(\mod\,5)$, then $S_{1}\cup D_{2}$ is a DTD-set of $G$ and $\dtd(G)\le |S_{1}\cup D_{2}|=2(n_{2}-1)/5+2=2(n-5)/5+2=2n/5<(n-1)/2$. If $n_{2}\equiv 2\,(\mod\,5)$, then $S_{1}\cup R_{2}$ is a DTD-set of $G$ and $\dtd(G)\le |S_{1}\cup R_{2}|=2(n_{2}-2)/5+2=2(n-6)/5+2=2(n-1)/5<(n-1)/2$. If $n_{2}=4$, then $G=D_{b}(4,4)$. Otherwise, if $n_{2}\equiv 3$ or $4\,(\mod\,5)$, then $n\ge 12$, and $S_{1}\cup T_{2}\cup\{v\}$ is a DTD-set of $G$ implying that $\dtd(G)\le |S_{1}\cup T_{2}\cup\{v\}|\le 2(n_{2}-3)/5+3=2(n-7)/5+3=(2n+1)/5<(n-1)/2$. Hence, $n_{i}\ne 4$ for $i=1,2$.

Suppose $n_{1}\equiv 0\,(\mod\,5)$. If $n_{2}\equiv 0\,(\mod\,5)$, then $S_{1}\cup S_{2}$ is a DTD-set of $G$ and $\dtd(G)\le |S_{1}\cup S_{2}|=2n_{1}/5+2n_{2}/5=2n/5<(n-1)/2$. If $n_{2}\equiv 1\,(\mod\,5)$, then $S_{1}\cup D_{2}$ is a DTD-set of $G$ and $\dtd(G)\le |S_{1}\cup D_{2}|=2n_{1}/5+2(n_{2}-1)/5=(2n-1)/5<(n-1)/2$. If $n_{2}\equiv 2\,(\mod \, 5)$, then $S_{1}\cup R_{2}$ is a DTD-set of $G$ and $\dtd(G)\le |S_{1}\cup R_{2}|=2n_{1}/5+2(n_{2}-2)/5=2(n-2)/5<(n-1)/2$. Otherwise, if $n_{2}\equiv\,3$ or $4\,(\mod\,5)$, $S_{1}\cup T_{2}\cup\{v\}$ is a DTD-set of $G$ and $\dtd(G)\le |S_{1}\cup T_{2}|+1\le 2n_{1}/5+2(n_{2}-3)/5+1=(2n-1)/5<(n-1)/2$. Hence, $n_{i}\not\equiv 0\,(\mod\,5)$ for $i=1,2$.

Suppose $n_{1}\equiv 1\,(\mod\,5)$. If $n_{2}\equiv 1\,(\mod\,5)$, then $S_{1}\cup D_{2}$ is a DTD-set of $G$ and $\dtd(G)\le |S_{1}\cup D_{2}|=2(n_{1}-1)/5+2(n_{2}-1)/5+1=(2n+1)/5<(n-1)/2$. If $n_{2}\equiv 2\,(\mod\,5)$, then $S_{1}\cup R_{2}$ is a DTD-set of $G$ and $\dtd(G)\le |S_{1}\cup R_{2}|=2(n_{1}-1)/5+2(n_{2}-2)/5+1=(2n-1)/5<(n-1)/2$. Otherwise, if $n_{2}\equiv 3$ or $4\,(\mod\,5)$, $D_{1}\cup T_{2}\cup\{v\}$ is a DTD-set of $G$ and $\dtd(G)\le |D_{1}\cup T_{2}\cup\{v\}|\le 2(n_{1}-1)/5+2(n_{2}-3)/5+1=(2n-3)/5<(n-1)/2$. Hence, $n_{i}\not\equiv 1\,(\mod\,5)$ for $i=1,2$ and, further, $n\ge 10$.

Suppose $n_{1}\equiv 2\,(\mod\,5)$. If $n_{2}\equiv 2\,(\mod\,5)$, then $R_{1}\cup R_{2}\cup\{u,v\}$ is a DTD-set of $G$ and $\dtd(G)\le |R_{1}\cup R_{2}\cup\{u,v\}|=2(n_{1}-2)/5+2(n_{2}-2)/5+2=2(n+1)/5<(n-1)/2$. Otherwise, if $n_{2}\equiv 3$ or $4\,(\mod\,5)$, $R_{1}\cup T_{2}\cup\{u,v\}$ is a DTD-set of $G$ and $\dtd(G)\le |R_{1}\cup T_{2}\cup \{u,v\}|\le 2(n_{1}-2)/5+2(n_{2}-3)/5+2=2n/5<(n-1)/2$. Hence, $n_{i}\not\equiv 2\,(\mod\,5)$ for $i=1,2$ and, further, $n\ge 11$.

Suppose $n_{1}\equiv 3\,(\mod\,5)$. If $n_{2}\equiv 3$ or $4\,(\mod\,5)$, then $T_{1}\cup T_{2}\cup\{u,v\}$ is a DTD-set of $G$ and $\dtd(G)\le |T_{1}\cup T_{2}\cup\{u,v\}|\le 2(n_{1}-3)/5+2(n_{2}-3)/5+2=2(n-1)/5<(n-1)/2$. Hence, $n_{i}\not\equiv 3\,(\mod\,5)$ for $i=1,2$ and, further, $n\ge 18$.

Suppose $n_{i}\equiv 4\,(\mod\,5)$ for $i=1,2$. Then, $T_{1}\cup T_{2}\cup\{u,v\}$ is a DTD-set of $G$ and $\dtd(G)\le |T_{1}\cup T_{2}\cup\{u,v\}|\le 2(n_{1}-4)/5+2(n_{2}-4)/5+2=2(n-3)/5<(n-1)/2$.~\qed

We now return to the proof of Proposition~\ref{db:ob}. If $G\in\cD_{b}\cup\cG_{b}$, then it is straightforward to check that $\dtd(G)\ge (n-1)/2$. In order to prove the necessity, let $G$ be such that $\dtd(G)\ge (n-1)/2$. By Claim P, $\ell\ge 1$. Hence $G=D_{b}(n_{1},n_{2},\ell)$ and $n\ge 7$.  Let $G_{1}$ and $G_{2}$ denote the two cycles in $G$ where $G_{i}\cong C_{n_{i}}$ for $i=1,2$. Let $uw_{1}\ldots w_{\ell}v$ be the path joining $G_{1}$ to $G_{2}$ where $u\in V(G_{1})$ and $v\in V(G_{2})$. Let $G_{2}=vv_{1}\ldots v_{n_{1}-1}v$. Let $F=G-V(G_{2})$. Then, $F$ is a key $L_{n_{1},\ell}$. Let $F$ have order $n_{F}$, and so $n_{F}=n_{1}+\ell=n-n_{2}$.

Let $S_{1}$ be a $\dtd(F)$-set and choose $S_{1}$ to contain the leaf vertex $w_{\ell}$ of $F$ if possible. Then, $|S_{1}|$ is given by Table~3. Further, $w_{\ell-1}\in S_{1}$. Let $R_{1}=S_{1}\setminus\{w_{\ell-1}\}$. Then, $|R_{1}|= |S_{1}|-1$. Let $D_{1}$ be a  $\dtd(F-\{w_{\ell-1},w_{\ell}\})$-set, possibly $|D_{1}|=|S_{1}|-1$ or $|D_{1}|=|S_{1}|-2$.  Let $S_{2}$ be a $\dtd(G_{2})$-set which contains $v$. By Proposition~\ref{cycle:ob}, $|S_{2}|=2n_{2}/5$ if $n\equiv 0\,(\mod\,5)$, and $|S_{2}|=\lceil 2(n+1)/5\rceil$ otherwise. Let $D_{2}=\{v_{i} \mid i\equiv 0$ or $1\,(\mod\,5)\}$. Let $R_{2}=\{v_{i} \mid i\equiv 3$ or $4\,(\mod\,5)\}$. Let $T_{2}=\{v_{i} \mid i\equiv 0$ or $4\,(\mod\,5)\}$. We proceed with the following series of claims.

\begin{unnumbered}{Claim~Q.1}
If $n_{1}\equiv 0\,(\mod\,5)$, then $G\in\{D_{b}(3,5,1)\}\subset G_{b}$.
\end{unnumbered}
\textbf{{Proof of Claim~Q.1}} Suppose $n_{1}\equiv 0\,(\mod\,5)$. We examine each possibility for $\ell$ in turn.

Suppose $\ell\equiv 0\,(\mod\,5)$. Then, $n\ge 9$. By Table~3, $|S_{1}|=2(n_{1}+\ell)/5$ and $S_{1}$ contains $w_{\ell}$. If $n_{2}\equiv 0\,(\mod\,5)$, then $S_{1}\cup S_{2}$ is DTD-set of $G$ and $\dtd(G)\le |S_{1}\cup S_{2}|=2n/5<(n-1)/2$, a contradiction. If $n_{2}\equiv 1\,(\mod\,5)$, then $S_{1}\cup D_{2}$ is a DTD-set of $G$, and $\dtd(G)\le |S_{1}\cup D_{2}|=2(n_{1}+\ell)/5+2(n_{2}-1)/5=2(n-1)/5<(n-1)/2$, a contradiction. If $n_{2}\equiv 2\,(\mod\,5)$, then $S_{1}\cup R_{2}$ is a DTD-set of $G$, and $\dtd(G)\le |S_{1}\cup R_{2}|=2(n_{1}+\ell)/5+2(n_{2}-2)/5=2(n-2)/5$, a contradiction. If $G_{2}=C_{4}$, then $S_{1}\cup S_{2}$ is a DTD-set of $G$, and $\dtd(G)\le |S_{1}\cup S_{2}|=2(n-4)/5+2=2(n+1)/5<(n-1)/2$, a contradiction. If $n_{2}\equiv 3$ or $4\,(\mod\,5)$, then $S_{1}\cup T_{2}\cup\{v\}$ is a DTD-set of $G$, and $\dtd(G)\le |S_{1}\cup T_{2}\cup\{v\}|\le 2(n_{1}+\ell)/5+2(n_{2}-3)/5+1=(2n-1)/5<(n-1)/2$, a contradiction.

Suppose $\ell\equiv 1\,(\mod\,5)$. By Table~3, $|S_{1}|=2(n_{1}+\ell-1)/5$. If $n_{2}\equiv 0\,(\mod\,5)$, then $S_{1}\cup S_{2}$ is a DTD-set of $G$, and $\dtd(G)\le |S_{1}\cup S_{2}|=2(n_{1}+\ell-1)/5+2n_{2}/5=2(n-1)/5<(n-1)/2$, a contradiction. If $n_{2}\equiv 1\,(\mod\,5)$, then $S_{1}\cup S_{2}$ is a DTD-set of $G$, and so $\dtd(G)\le |S_{1}\cup S_{2}|=2(n_{1}+\ell-1)/5+2(n_{2}-1)/5+1=2(n-2)/5+1=(2n+1)/5<(n-1)/2$, a contradiction. If $n_{2}\equiv 2\,(\mod\,5)$, then $S_{1}\cup R_{2}\cup\{w_{\ell}\}$ is a DTD-set of $G$, and $\dtd(G)\le |S_{1}\cup R_{2}\cup\{w_{\ell}\}|=2(n_{1}+\ell-1)/5+2(n_{2}-2)/5+1=2(n-3)/5+1=(2n-1)/5<(n-1)/2$, a contradiction. If $G_{2}=C_{4}$, then $S_{1}\cup S_{2}$ is a DTD-set of $G$, and $\dtd(G)\le |S_{1}\cup S_{2}|=2(n_{1}+\ell-1)/5+2=2(n-5)/5+2=2n/5<(n-1)/2$, a contradiction. Otherwise, if $n_{2}\equiv 3$ or $4(\mod\,5)$, then $S_{1}\cup S_{2}$ is a DTD-set of $G$, and $\dtd(G)\le |S_{1}\cup S_{2}|=2(n_{1}+\ell-1)/5+2(n_{2}-3)/5+2= 2(n+1)/5\le (n-1)/2$. We have equality if $n=9$ and $G=D_{b}(3,5,1)=G_{1}(1,0)\in \cG_{b}$ otherwise $\dtd(G)<(n-1)/2$, a contradiction.

Suppose $\ell\equiv 2\,(\mod\,5)$. Then, $n\ge 10$. By Table~3, $|S_{1}|=\lceil 2(n_{1}+\ell)/5\rceil =(2(n_{1}+\ell)+1)/5$. Further, $|R_{1}|=|S_{1}|-1$ and $R_{1}$ DT-dominates $F-w_{\ell}$. If $n_{2}\equiv 0\,(\mod\,5)$, then $R_{1}\cup S_{2}$ is a DTD-set of $G$, and $\dtd(G)\leq|R_{1}\cup S_{2}|=(2(n_{1}+\ell)+1)/5+2n_{2}/5-1=(2n-4)/5<(n-1)/2$, a contradiction. If $n_{2}\equiv 1\,(\mod\,5)$, then $R_{1}\cup S_{2}$ is a DTD-set of $G$, and $\dtd(G)\leq|R_{1}\cup S_{2}|=(2(n_{1}+\ell)+1)/5+2(n_{2}-1)/5=(2n-1)/5<(n-1)/2$, a contradiction. If $n_{2}\equiv 2\,(\mod\,5)$, then $S_{1}\cup R_{2}$ is a DTD-set of $G$, and $\dtd(G)\leq|S_{1}\cup R_{2}|=(2(n_{1}+\ell)+1)/5+2(n_{2}-2)/5-1=(2n-3)/5<(n-1)/2$, a contradiction. If $n_{2}\equiv 3$ or $4\,(\mod\,5)$, then $|R_{1}\cup S_{2}|$ is a DTD-set of $G$, and $\dtd(G)\le |R_{1}\cup S_{2}|\le (2(n_{1}+\ell)+1)/5+2(n_{2}-3)/5=2n/5<(n-1)/2$, a contradiction.

Suppose $\ell\equiv 3\,(\mod\,5)$. Then, $n\ge 11$. By Table~3 $|S_{1}|=2(n_{1}+\ell+2)/5$. Further, $|D_{1}|=|S_{1}|-2$. If $n_{2}\equiv 0$, $1$, $2$, $3$ or $4\,(\mod\,5)$, then $D_{1}\cup S_{2}$ is a DTD-set of $G$, and $\dtd(G)\le |D_{1}\cup S_{2}|= 2(n_{1}+\ell+2)/5+2(n_{2}+3)/5-2=2n/5<(n-1)/2$, a contradiction.

Suppose $\ell\equiv 4\,(\mod\,5)$. Then, $n\ge 12$. By Table~3 $|S_{1}|=2(n_{1}+\ell+1)/5$. Further, $|R_{1}|=|S_{1}|-1$, $w_{\ell}\in R_{1}$, and $R_{1}$ DT-dominates $F-w_{\ell}$. If $n_{2}\equiv 0\,(\mod\,5)$, then $R_{1}\cup S_{2}$ is a DTD-set of $G$, and $\dtd(G)\le |R_{1}\cup S_{2}|= 2(n_{1}+\ell+1)/5+2n_{2}/5-1=(2n-3)/5<(n-1)/2$, a contradiction. If $n_{2}\equiv 1\,(\mod\,5)$, then $R_{1}\cup D_{2}$ is a DTD-set of $G$, and $\dtd(G)\le |R_{1}\cup D_{2}|= 2(n_{1}+\ell+1)/5+2(n_{2}-1)/5=2n/5<(n-1)/2$, a contradiction. If $n_{2}\equiv 2\,(\mod\,5)$, then $R_{1}\cup R_{2}\cup\{v\}$ is a DTD-set of $G$, and $\dtd(G)\le |R_{1}\cup R_{2}\cup\{v\}|= 2(n_{1}+\ell+1)/5+2(n_{2}-2)/5=2(n-1)/5<(n-1)/2$, a contradiction. If $G_{2}=C_{4}$, then $R_{1}\cup S_{2}$ is a DTD-set of $G$, and $\dtd(G)\le |R_{1}\cup S_{2}|= 2(n_{1}+\ell+1)/5+1 = 2(n-3)/5+1=(2n-1)/5<(n-1)/2$, a contradiction. Otherwise, if $n_{2}\equiv\,3$ or $4\,(\mod\,5)$, then $R_{1}\cup T_{2}\cup \{v\}$ is a DTD-set of $G$, and $\dtd(G)\le |R_{1}\cup T_{2}\cup\{v\}|\le 2(n_{1}+\ell+1)/5+2(n_{2}-3)/5=2(n-2)/5<(n-1)/2$, a contradiction.~\smallqed

By Claim Q.1, it follows, relabelling indices if necessary, that $n_{i}\not\equiv 0\,(\mod\,5)$ for $i=1,2$.

\begin{unnumbered}{Claim Q.2} If $n_{1}\equiv 1\,(\mod\,5)$, then $G\in\{D_b(4,6,1), D_{b}(3,6,2)\}\subset\cG_{b}$.
\end{unnumbered}
\textbf{Proof of Claim Q.2} Suppose $n_{1}\equiv 1\,(\mod\,5)$. We examine each possibility for $\ell$ in turn.

Suppose $\ell\equiv 0\,(\mod\,5)$. Then, $n\ge 14$. By Table~3, $|S_{1}|=(2(n_{1}+\ell)+3)/5$ and $w_{\ell}\in S_{1}$. Further, we have $R_{1}$ with $|R_{1}|=|S_{1}|-1$, $w_{\ell}\in R_{1}$ and $R_{1}$ DT-dominates $F-w_{\ell}$. If $n_{2}\equiv 1\,(\mod\,5)$, then $S_{1}\cup D_{2}$ is a DTD-set of $G$, and $\dtd(G)\le |S_{1}\cup D_{2}|=(2(n_{1}+\ell)+3)/5+2(n_{2}-1)/5=(2n+1)/5<(n-1)/2$, a contradiction. If $n_{2}\equiv 2\,(\mod\,5)$, then $S_{1}\cup R_{2}$ is a DTD-set of $G$, and $\dtd(G)\le |S_{1}\cup R_{2}|=(2(n_{1}+\ell)+3)/5+2(n_{2}-2)/5=(2n-1)/5<(n-1)/2$, a contradiction. If $G_{2}=C_{4}$, then $R_{1}\cup S_{2}$ is a DTD-set of $G$, and $\dtd(G)\le |R_{1}\cup S_{2}|=(2(n_{1}+\ell)+3)/5+1=(2(n-4)+3)/5+1=2n/5<(n-1)/2$, a contradiction. Otherwise, if $n_{2}\equiv 3$ or $4\,(\mod\,5)$, then  $S_{1}\cup T_{2}\cup\{v\}$ is a DTD-set of $G$, and $\dtd(G)\le |S_{1}\cup T_{2}\cup \{v\}|\le (2(n_{1}+\ell)+3)/5+2(n_{2}-3)/5+1=2(n+1)/5<(n-1)/2$, a contradiction.

Suppose $\ell\equiv 1\,(\mod\,5)$. Then, $n\ge 10$. By Table~3, $|S_{1}|=((2n_{1}+\ell)+1)/5$ and $w_{\ell}\in S_{1}$.
If $n_{2}\equiv 1\,(\mod\,5)$, then $S_{1}\cup D_{2}$ is a DTD-set of $G$, and $\dtd(G)\le (2(n_{1}+\ell)+1)/5+2(n_{2}-1)/5=(2n-1)/5<(n-1)/2$, a contradiction. If $n_{2}\equiv 2\,(\mod\,5)$, then $S_{1}\cup R_{2}$ is a DTD-set of $G$, and $\dtd(G)\le |S_{1}\cup R_{2}|=((2n_{1}+\ell)+1)/5+2(n_{2}-2)/5=(2n-3)/5<(n-1)/2$. If $G_{2}=C_{4}$, then $n\ge 11$ and $S_{1}\cup S_{2}$ is a DTD-set of $G$, and $\dtd(G)\le (2(n_{1}+\ell)+1)/5+2=(2n+3)/5$. We have $\dtd(G)=(n-1)/2$ if $G=D_{b}(4,6,1)\cong G_{2}(0,1)\subset G_{b}$ for $n=11$, and, for $n\ge 11$, $\dtd(G)<(n-1)/2$, a contradiction. Otherwise, if $n_{2}\equiv 3$ or $4\,(\mod\,5)$, then  $S_{1}\cup T_{2}\cup\{v\}$ is a DTD-set of $G$, and $\dtd(G)\le |S_{1}\cup T_{2}\cup \{v\}|\le (2(n_{1}+\ell)+1)/5+2(n_{2}-3)/5+1=2n/5<(n-1)/2$, a contradiction.

Suppose $\ell\equiv 2\,(\mod\,5)$. Then, $n\ge 11$. By Table~3, $|S_{1}|=(2(n_{1}+\ell)-1)/5$. If $n_{2}\not\equiv 2\,(\mod\,5)$, then $S_{1}\cup S_{2}$ is a DTD-set of $G$, and $\dtd(G)\le |S_{1}\cup S_{2}|\le (2(n_{1}+\ell)-1)/5+2(n_{2}+2)/5=(2n+3)/5$. We have $\dtd(G)=(n-1)/2$ if $G=D_{b}(3,6,2)\cong G_{2}(1,0)\subset \cG_{b}$ and $n=11$, otherwise $\dtd(G)<(n-1)/2$, a contradiction. If $n_{2}\equiv 2\,(\mod\,5)$, then $S_{1}\cup R_{2}\cup\{w_{\ell}\}$ is a DTD-set of $G$, and $\dtd(G)\le |S_{1}\cup R_{2}\cup\{w_{\ell}\}|=(2(n_{1}+\ell)-1)/5+2(n_{2}-2)/5+1=2n/5<(n-1)/2$, a contradiction.

Suppose $\ell\equiv 3\,(\mod\,5)$. Then, $n\ge 12$. By Table~3, $|S_{1}|=2(n_{1}+\ell+1)/5$. Further, we have $R_{1}$ with $|R_{1}|=|S_{1}|-1$, and $R_{1}$ DT-dominates $F-w_{\ell}$. If $n_{2}\equiv\,1$, $2$, $3$ or $4\,(\mod\,5)$, then $R_{1}\cup S_{2}$ is a DTD-set of $G$, and $\dtd(G)\le |R_{1}\cup S_{2}|=2(n_{1}+\ell+1)/5+2(n_{2}+3)/5-1=(2n+3)/5<(n-1)/2$, a contradiction.

Suppose $\ell\equiv 4\,(\mod\,5)$. Then, $n\ge 13$. By Table~3, $|S_{1}|=(2(n_{1}+\ell)+3)/5$. Further, there is a set $D_{1}$ with $|D_{1}|=|S_{1}|-2$. If $n_{2}\equiv 1$, $2$, $3$ or $4\,(\mod\,5)$, then $D_{1}\cup S_{2}$ is a DTD-set of $G$, and $\dtd(G)\le |D_{1}\cup S_{2}|\le (2(n_{1}+\ell)+3)/5+2(n_{2}+3)/5-2=(2n-1)/5<(n-1)/2$, a contradiction.~\smallqed

By Claim Q.2, it follows, relabelling indices if necessary, that $n_{i}\not\equiv 1\,(\mod\,5)$ for $i=1,2$.

\begin{unnumbered}{Claim Q.3} If $n_{1}\equiv 2\,(\mod\,5)$, then $G\in\{D_{b}(4,7,2), D_{b}(3,7,3)\}\subset\cD_{b}$.
\end{unnumbered}
\textbf{Proof of Claim Q.3} Suppose $n_{1}\equiv\,2\,\mod\,5$. We examine each possibility for $\ell$ in turn.

Suppose $\ell\equiv 0\,(\mod\,5)$. Then, $n\ge 15$. By Table~3, $|S_{1}|=2(n_{1}+\ell+3)/5$. Further, $|D_{1}|=|S_{1}|-2$. If $n_{2}\equiv 2$, $3$, or $4\,(\mod\, 5)$, then $D_{1}\cup S_{2}$ is a DTD-set of $G$ and $\dtd(G)\le |D_{1}\cup S_{2}|\le 2(n_{1}+\ell+3)/5+2(n_{2}+3)/5-2=2(n+1)/5<(n-1)/2$.

Suppose $\ell\equiv 1\,(\mod\,5)$. Then, $n\ge 11$. By Table~3, $|S_{1}|=2(n_{1}+\ell+2)/5$ and $w_{\ell}\in S_{1}$. Further, we have $R_{1}$ with $|R_{1}|=|S_{1}|-1$, $w_{\ell}\in R_{1}$ and $R_{1}$ DT-dominates $F-w_{\ell}$. If $G_{2}=C_{4}$, then $n\ge 12$, $R_{1}\cup S_{2}$ is a DTD-set of $G$, and $\dtd(G)\le 2(n_{1}+\ell+2)/5+1=2(n-2)/5+1=(2n+1)/5<(n-1)/2$, a contradiction. Hence, $G_{2}\ne C_{4}$. Otherwise, if $n_{2}\equiv 2$, $3$, or $4\,(\mod\,5)$, then $R_{1}\cup T_{2}\cup\{v\}$ is a DTD-set of $G$, and $\dtd(G)\le |R_{1}\cup T_{2}\cup\{v\}|\le 2(n_{1}+\ell+2)/5+2(n_{2}-2)/5=2n/5<(n-1)/2$, a contradiction.

Suppose $\ell\equiv 2\,(\mod\,5)$. Then, $n\ge 12$. By Table~3, $|S_{1}|=2(n_{1}+\ell+1)/5$ and $w_{\ell}\in S_{1}$. If $G_{2}=C_{4}$, then $n\ge 13$, $S_{1}\cup S_{2}$ is a DTD-set of $G$, and $\dtd(G)\le|S_{1}\cup S_{2}|=2(n_{1}+\ell+1)/5+2=2(n+2)/5$. We have $\dtd(G)=(n-1)/2$ and $G=D_{b}(4,7,2)\in \cD_{b}$ if $n=13$, and, for $n\ge 14$, $\dtd(G)<(n-1)/2$, a contradiction. Otherwise, if $n_{2}\equiv 2$, $3$, or $4\,(\mod\,5)$, then $S_{1}\cup T_{2}\cup\{v\}$ is a DTD-set of $G$, and $\dtd(G)\le |S_{1}\cup T_{2}\cup\{v\}|\le 2(n_{1}+\ell+1)/5+2(n_{2}-2)/5+1=(2n+3)/5<(n-1)/2$.

Suppose $\ell\equiv 3\,(\mod\,5)$. Then, $n\ge 13$. By Table~3, $|S_{1}|=2(n_{1}+l)/5$. If $n_{2}\equiv 3$ or $4\,(\mod\,5)$, then $S_{1}\cup S_{2}$ is a DTD-set of $G$, and $\dtd(G)\le 2(n_{1}+\ell)/5+2(n_{2}+2)/5=(2n+4)/5$. We have $\dtd(G)=(n-1)/2$ if $G=D_{b}(3,7,3)\in \cD_{b}$ with $n=13$ and, for $n\ge 14$, $\dtd(G)<(n-1)/2$, a contradiction. If $n_{2}\equiv 2\,(\mod\,5)$, then $S_{1}\cup \{w_{\ell}\}\cup R_{2}$ is a DTD-set of $G$, and $\dtd(G)\le |S_{1}\cup R_{2}|+1= 2(n_{1}+\ell)/5+2(n_{2}-2)/5+1=(2n+1)/5<(n-1)/2$, a contradiction.

Suppose $\ell\equiv 4\,(\mod\,5)$. Then, $n\ge 14$. By Table~3, $|S_{1}|=(2(n_{1}+\ell)+1)/5$. Further, we have $R_{1}$ with $|R_{1}|=|S_{1}|-1$, and $R_{1}$ DT-dominates $F-w_{\ell}$. If $n_2\equiv 2$, $3$ or $4\,(\mod\,5)$, then $R_{1}\cup S_{2}$ is a DTD-set of $G$, and $\dtd(G)\le |R_{1}\cup S_{2}|\le (2(n_{1}+\ell)+1)/5+2(n_{2}+3)/5-1=2(n+1)/5<(n-1)/2$, a contradiction.~\smallqed

By Claim Q.3, it follows, relabelling indices if necessary, that $n_{i}\not\equiv 2\,(\mod\,5)$ for $i=1,2$.

\begin{unnumbered}{Claim Q.4}
If $n_1\equiv 3\,(\mod\,5)$, then $G\in\{D_{b}(3,3,1), D_{b}(3,4,1),$ $D_{b}(4,8,1),$ $D_{b}(3,3,2),$ \\$D_{b}(3,8,2),$ $D_{b}(3,4,6), D_{b}(3,3,7)\}\subset\cD_{b}$ or $G\in\{D_{b}(3,4,2),$ $D_{b}(3,3,3)\}\subset\cG_{b}$.
\end{unnumbered}
\textbf{Proof of Claim Q.4} Suppose $n_{1}\equiv 3\,(\mod\,5)$.  We examine each possibility for $\ell$ in turn.

Suppose $\ell\equiv 0\,(\mod\,5)$. Then, $n\ge 11$. By Table~3, $|S_{1}|=2(n_{1}+\ell+2)/5$ and $w_{\ell}\in S_{1}$. Further, $|R_{1}|=|S_{1}|-1$, $w_{\ell}\in R_{1}$ and $R_{1}$ DT-dominates $F-w_{\ell}$. If $G_{2}=C_{4}$, then $n\ge 12$, $R_{1}\cup S_{2}$ is a DTD-set of $G$, and $\dtd(G)\le |S_{1}\cup S_{2}|\le 2(n_{1}+\ell+2)/5+1=(2n+1)/5<(n-1)/2$, a contradiction. Hence, $G_{2}\ne C_{4}$. Otherwise, if $n_{2}\equiv 3$ or $4\,(\mod\,5)$, $R_{1}\cup T_{2}\cup\{v\}$ is a DTD-set of $G$, and $\dtd(G)\le |R_{1}\cup T_{2}\cup\{v\}|\le 2(n_{1}+\ell+2)/5+2(n_{2}-2)/5=2n/5<(n-1)/2$, a contradiction.

Suppose $\ell\equiv 1\,(\mod\,5)$. Then, $n\ge 7$. By Table~3, $|S_{1}|=2(n_{1}+\ell+1)/5$ and $w_{\ell}\in S_{1}$. If $G_{2}=C_{4}$, then $n\ge 8$, and $S_{1}\cup S_{2}$ is a DTD-set of $G$. We have $\dtd(G)\le |S_{1}\cup S_{2}|=2(n_{1}+\ell+1)/5+2=2(n+2)/5$. We have $\dtd(G)\ge (n-1)/2$ if $G\in\{D_{b}(3,4,1), D_{b}(4,8,1), D_{b}(3,4,6)\}\subset\cD_{b}$ with $n\le 13$ and, for $n\ge 14$, $\dtd(G)<(n-1)/2$, a contradiction. Hence, $G_{2}\ne C_{4}$. Otherwise, if $n_{2}\equiv 3$ or $4\,(\mod\,5)$, then $S_{1}\cup T_{2}\cup\{v\}$ is a DTD-set of $G$, and $\dtd(G)\le |S_{1}\cup T_{2}|\le 2(n_{1}+\ell+1)/5+2(n_{2}-2)/5+1=(2n+1)/5$. We have $\dtd(G)=(n-1)/2$ if $G=D_{b}(3,3,1)$ with $n=7$, otherwise $\dtd(G)<(n-1)/2$, a contradiction.

Suppose $\ell\equiv 2\,(\mod\,5)$. Then, $n\ge 8$. By Table~3, $|S_{1}|=2(n_{1}+\ell)/5$. If $n_{2}\equiv 3\,(\mod\,5)$, then $S_{1}\cup S_{2}$ is a DTD-set of $G$, and $\dtd(G)\le |S_{1}\cup S_{2}|\le 2(n_{1}+\ell)/5+2(n_{2}+2)/5=2(n+2)/5$. We have $\dtd(G)\ge (n-1)/2$ and $G\in\{D_{b}(3,3,2), D_{b}(3,8,2), D_{b}(3,3,7)\}\subset\cD_{b}$ for $n\le 13$ and, for $n\ge 14$, $\dtd(G)<(n-1)/2$, a contradiction. If $n_{2}\equiv 4\,(\mod\,5)$, then $n\ge 9$, and $S_{1}\cup S_{2}$ is a DTD-set of $G$. Further, $\dtd(G)\le |S_{1}\cup S_{2}|\le 2(n_{1}+\ell)/5+2(n_{2}+1)/5=2(n+1)/5$.  We have $\dtd(G)=(n-1)/2$ and $G=D_{b}(3,4,2)\cong G_{0}(1,1)\in\cG_{b}$ with $n=9$, otherwise $\dtd(G)<(n-1)/2$, a contradiction.

Suppose $\ell\equiv 3\,(\mod\,5)$. Then, $n\ge 9$. By Table~3, $|S_{1}|=(2(n_{1}+l)+3)/5$. Further, $|R_{1}|=|S_{1}|-1$ and $R_{1}$ DT-dominates $F-w_{\ell}$. If $n_{2}\equiv 3$ or $4\,(\mod\,5)$, then $R_{1}\cup S_{2}$ is a DTD-set of $G$, and $\dtd(G)\le |R_{1}\cup S_{2}|\le (2(n_{1}+\ell)+3)/5+2(n_{2}+2)/5-1=2(n+1)/5$. We have $\dtd(G)=(n-1)/2$, and $G=D_{b}(3,3,3)\cong G_{0}(2,0)\in\cG_{b}$ with $n=9$,  otherwise $\dtd(G)<(n-1)/2$, a contradiction.

Suppose $\ell\equiv 4\,(\mod\,5)$. Then, $n\ge 10$. By Table~3, $|S_{1}|=2(n_{1}+\ell+2)/5$. Further, we have the set $D_{1}$ and $|D_{1}|=|S_{1}|-2$. If $n_{2}\equiv 3$ or $4\,(\mod\,5)$, then $D_{1}\cup S_{2}$ is a DTD-set of $G$, and $\dtd(G)\le 2(n_{1}+\ell+3)/5+2(n_{2}+2)/5-2=2n/5<(n-1)/2$, a contradiction.~\smallqed

By Claim Q.4, it follows, relabelling indices if necessary, that $n_{i}\not\equiv 3\,(\mod\,5)$ for $i=1,2$.

\begin{unnumbered}{Claim Q.5}
If $n_{1}\equiv 4\,(\mod\,5)$, then $G\in\{D_{b}(4,4,5)\}\subset \cG_{b}$ or $G\in\{D_{b}(4,4,1)\}\subset\cG_{b}$.
\end{unnumbered}
\textbf{Proof of Claim Q.5} Suppose $G_{1}=G_{2}=C_{4}$. Then, $n=\ell+8$. Let $P:w_{1}\ldots w_{\ell}$ be the path joining $u$ to $v$. Let $G_{1}=uu_{1}u_{2}u_{3}u$ and $G_{2}=vv_{1}v_{2}v_{3}v$. Recall that $u$ is adjacent to $w_{1}$ and $v$ is adjacent to $w_{\ell}$. Let $T'=\{u, u_{1}, v, v_{1}\}$. Now, for $\ell\ge 1$ construct the set $T$ as follows: If $\ell=1,2,3$ then let $T=\emptyset$. If $\ell=4$, then let $T=\{w_{\ell}\}$. If $\ell\ge 5$, then let $T$ be a $\dtd(P-\{w_{1},w_{2},w_{\ell}\})$-set. Then, in this case, by Proposition~\ref{path:prop:ap}, $|T|=\lceil 2(\ell-2)/5\rceil+1$ if $\ell\equiv  4\,(\mod\,5)$ and $|T|=\lceil 2(\ell-2)/5\rceil$ otherwise. If $1\le \ell\le 4$, then $G\in\{D_{b}(4,4,1), D_{b}(4,4,2), D_{b}(4,4,3), D_{b}(4,4,4)\}$ and  $T'\cup T$ is a DTD-set of $G$. Thus $\dtd(G)\le |T'\cup T|$ and either $G=D_b(4,4,1)$ or $\dtd(G)<(n-1)/2$, a contradiction. Hence, assume $\ell\ge 5$, and so $n\ge 13$. Suppose $\ell\not\equiv 4\,(\mod\,5)$. Then $\dtd(G)\le |T'\cup T|\le 2\ell/5+4=2(n-8)/5+4=2(n+2)/5$. We have, for $n=13$, $G=D_{b}(4,4,5)$ and $\dtd(G)=(n-1)/2$, otherwise $\dtd(G)<(n-1)/2$. Suppose $\ell\equiv 4\,(\mod\,5)$. Then $\dtd(G)\le  |T'\cup T|= 2(\ell+1)/5+4=(2n+6)/5$. We have $n\ge 17$. If $n=17$, then $G=D_{b}(4,4,9)$ and it can be verified that $\dtd(G)=7<(n-1)/2$, otherwise $\dtd(G)\le |T'\cup T| = (2n+6)/5<(n-1)/2$, a contradiction. Thus at least one of $G_{1}$ and $G_{2}$, say $G_{1}\ne C_{4}$.

Suppose $n_{1}\equiv 4\,(\mod\,5)$, with $n_{1}\ge 9$. Suppose $\ell\not\equiv 4\,(\mod\,5)$ and $G_{2}=C_{4}$. By Table~3, $|S_{1}|\le 2(n_{1}+l+1)/5$. Since $G_{2}= C_{4}$, $n\ge 14$ and $S_{1}\cup S_{2}$ is a DTD-set of $G$. We have $\dtd(G)\le |S_{1}\cup S_{2}|\le 2(n_{1}+\ell+1)/5+2=2(n+2)/5<(n-1)/2$, a contradiction. Suppose $\ell\equiv 4\,(\mod\,5)$ and $G_{2}=C_{4}$. By Table~3, $|S_{1}|=2(n_{1}+\ell+2)/5$ and we have $|D_{1}|=|D_{1}|-2$. Then, $\dtd(G)\le |D_{1}\cup S_{2}|\le 2(n_{1}+\ell+2)/5=2(n-2)/5<(n-1)/2$, a contradiction. Hence, $G_{2}\ne C_{4}$ and it follows that $n\ge 19$. Further, by Table~3, $|S_{1}|\le 2(n_{1}+\ell+2)$. By our previous claims we have $n_{2}\equiv 4\,(\mod\,5)$, and so $\dtd(G)\le |S_{1}\cup S_{2}|=2(n_{1}+\ell+2)/5+2(n_{2}+1)/5=2(n+3)/5<(n-1)/2$.~\smallqed 

Proposition~\ref{db:ob} follows by Claims Q1--Q5.~\qed

Recall the statement of Observation~\ref{s:ob}.

\begin{unnumbered}{Observation~\ref{s:ob}}
A graph $G$ of order $3\le n\le 7$ is $\frac{1}{2}$-minimal if and only if $G\in\{B_{1},$ $B_{2}, D(3,3), D(4,4), D_b(3,4), D_{b}(3,3,1)\}\cup\{C_{3}, C_{4}, C_{5}, C_{6}, C_{7}\}$.
\end{unnumbered}
\textbf{Proof of Observation~\ref{s:ob}.} Suppose $G$ is a $\frac{1}{2}$-minimal graph of order $3\le n\le 7$ and let $\Delta(G)=\Delta$ be the maximum degree of a vertex in $G$. If $3\le n\le 6$, then $G\in\{B_{1}, D(3,3)\}\cup\{C_{3}, C_{4}, C_{5}, C_{6}\}$. Hence, $n=7$. Further, we may suppose $\dtd(G)\ge 3$ and $2\le \Delta\le 4$. If $\Delta=2$, then $G=C_{7}$ and $\dtd(G)=4$. Suppose $\Delta=3$. Let $u$ be a vertex of degree~3. If every vertex is within distance 2 of $u$, then $G=D_{b}(3,4)$ with $u$ the vertex of degree~3 in $G$ that belongs to the 4-cycle in $G$, or $\dtd(G)=2$, a contradiction. Hence, assume that there is vertex in $G-N[u]$ which is not adjacent to a vertex in $N(u)$. Then, $G=B_{2}$, $G=D_{b}(3,3,1)$, $G=D_{b}(3,4)$ with $u$ the vertex of degree~$3$ in $G$ that belongs to the 3-cycle in $G$, or $G$ is either not edge-minimal or has $\dtd(G)=2$, both of which give a contradiction. Suppose $\Delta=4$. Let $u$ be a vertex of degree $4$ in $G$, and let $V-N[u]=\{v,w\}$. If $v$ and $w$ have a common neighbor, then $\dtd(G)=2$, a contradiction. If $v$ and $w$ are adjacent and have no common neighbor, then $G$ contains a subgraph isomorphic to $C_{5}$ and $\dtd(G)=2$, a contradiction. It follows that $v$ and $w$ are at least distance 3 apart, and so $G=D(4,4)$.~\qed

\medskip

Recall the statement of Observation~\ref{sd:ob}.

\begin{unnumbered}
{Observation~\ref{sd:ob}}
Let $G$ be a graph with no isolated vertex and let $F$ be a $4$-subdivision of $G$. Then, $\dtd(F) \le \dtd(G)+2$.
\end{unnumbered}
\noindent{\bf Proof of Observation~\ref{sd:ob}.} Let $G=(V,E)$ be a graph with no isolated vertex. Let $e\in E$, where $e=uv$ be an edge of $G$. Let $F$ be obtained from $G$ by subdividing the edge $e$ four times, so that $F$ contains the path $P \colon uw_{1}w_{2}w_{3}w_{4}v$ which is not contained in $G$. Let $S$ be a $\dtd(G)$-set. We show that the claim is true by considering each of the following three cases: $u,v\in S$; $u\in S$ and $v\notin S$; $u,v\notin S$. Suppose $u,v$ in $S$. Then, $S\cup\{w_{1}, w_{4}\}$ is a DTD-set of $F$. Suppose $u\in S$ and $v\notin S$. Then, either $N_{G}(v)\cap S=\emptyset$ or $N_{G}(v)\cap S\neq\emptyset$. If $N_{G}(v)\cap S=\emptyset$, then $S\cup\{w_{3},w_{4}\}$ is a DTD-set of $F$. If $N_{G}(v)\cap S\ne \emptyset$, then $S\cup\{w_{2}, w_{4}\}$ is a DTD-set of $F$. Suppose $u,v\notin S$. Then, $S\cup\{w_{2},w_{3}\}$ is a DTD-set of $F$. Thus, in each case, $\dtd(F)\le |S|+2= \dtd(G)+2$. By symmetry, a similar result holds if $u\notin S$ and $v\in S$. This completes the proof.~\qed

\end{document}